\documentclass{article}
\usepackage{graphicx} 
\usepackage{rotating} 
\usepackage[a4paper, total={6in, 9in}]{geometry}
\usepackage{tikz}
\usepackage{tikz-cd}
\usepackage{amsfonts}
\usepackage{amsmath}
\usepackage{amsthm}
\usepackage{amssymb}
\usepackage{stmaryrd} 
\usepackage{comment}
\usepackage{mathrsfs}
\usepackage{xcolor}
\usepackage{enumitem}

\usepackage{hyperref}

\usepackage{placeins} 

\usepackage{subcaption}

\usepackage{diagbox}

\usepackage{standalone}

\usetikzlibrary{decorations.pathmorphing}
\usepackage{algorithm}
\usepackage[noend]{algpseudocode}

\tikzset{
  red_black_dash/.style={
    draw,
    line width=2pt,
    dash pattern=on 15pt off 15pt on 15pt off 15pt, 
    preaction={draw, red, line width=2pt}, 
    postaction={draw, black, line width=2pt, dash pattern=on 15pt off 15pt on 15pt off 15pt}  
  }
}
\tikzset{
  blue_black_dash/.style={
    draw,
    line width=2pt,
    dash pattern=on 15pt off 15pt on 15pt off 15pt, 
    preaction={draw, line width=2pt,cyan}, 
    postaction={draw, black, line width=2pt, dash pattern=on 15pt off 15pt on 15pt off 15pt}  
  }
}

\newcommand{\Tfourpiece}[1]{
  \begin{scope}[shift={(#1)}]
        \draw[line width=1mm, fill=blue!30] 
        (0,0) 
        -- (0,3) 
        -- (1,3) 
        -- (1,4) 
        -- (2,4) 
        -- (2,5)
        -- (3,5) 
        -- (3,6) 
        -- (4,6) 
        -- (4,3)
        -- (3,3) 
        -- (3,2) 
        -- (2,2) 
        -- (2,1)
        -- (1,1) 
        -- (1,0) 
        -- cycle;
  \end{scope}
}

\newcommand{\lineUp}[1]{
  \begin{scope}[shift={(#1)}]
    \draw (0, 0) -- (0, 1);   
  \end{scope}
}

\newcommand{\lineUpWArrow}[1]{
  \begin{scope}[shift={(#1)}]
    \draw[->] (0, 0) -- (0, 1);   
  \end{scope}
}

\newcommand{\lineUpBendRight}[1]{
  \begin{scope}[shift={(#1)}]
    \draw (0,0) to[bend right=50] (0,1);
  \end{scope}
}
\newcommand{\lineUpBendLeft}[1]{
  \begin{scope}[shift={(#1)}]
    \draw (0,0) to[bend left=50] (0,1);
  \end{scope}
}

\newcommand{\lineUpBendRightMore}[1]{
  \begin{scope}[shift={(#1)}]
    \draw (0,0) to[bend right=90] (0,1);
  \end{scope}
}
\newcommand{\lineUpBendLeftMore}[1]{
  \begin{scope}[shift={(#1)}]
    \draw (0,0) to[bend left=90] (0,1);
  \end{scope}
}

\newcommand{\lineRightBendUp}[1]{
  \begin{scope}[shift={(#1)}]
    \draw (0,0) to[bend left=50] (1,0);
  \end{scope}
}

\newcommand{\lineRightBendDown}[1]{
  \begin{scope}[shift={(#1)}]
    \draw (0,0) to[bend right=50] (1,0);
  \end{scope}
}

\newcommand{\lineRightBendUpMore}[1]{
  \begin{scope}[shift={(#1)}]
    \draw (0,0) to[bend left=90] (1,0);
  \end{scope}
}

\newcommand{\lineRightBendDownMore}[1]{
  \begin{scope}[shift={(#1)}]
    \draw (0,0) to[bend right=90] (1,0);
  \end{scope}
}

\newcommand{\zeroSmoothing}[1]{
  \begin{scope}[shift={(#1)}]
  \lineRightBendUpMore{0,0}
  \lineRightBendDownMore{0,1}
  \end{scope}
}

\newcommand{\oneSmoothing}[1]{
  \begin{scope}[shift={(#1)}]
  \lineUpBendRightMore{0,0}
  \lineUpBendLeftMore{1,0}
  \end{scope}
}

\newcommand{\negCrossing}[1]{
  \begin{scope}[shift={(#1)}]
    \draw (0,0)  -- (0.35,0.35);
    \draw (0.65,0.65)  -- (1,1);
    \draw (1,0) -- (0,1);
  \end{scope}
}
\newcommand{\posCrossing}[1]{
  \begin{scope}[shift={(#1)}]
    \draw (0,1)  -- (0.35,0.65);
    \draw (0.65,0.35)  -- (1,0);
    \draw (0,0) -- (1,1);
  \end{scope}
}
\newcommand{\negCrossingWArrows}[1]{
  \begin{scope}[shift={(#1)}]
    \draw (0,0)  -- (0.35,0.35);
    \draw[->] (0.65,0.65)  -- (1,1);
    \draw[->] (1,0) -- (0,1);
  \end{scope}
}

\newcommand{\Cob}{\operatorname{Cob}_{\bullet}^3}
\newcommand{\CCob}{\operatorname{CCob}^3}

\newcommand{\term}[1]{\textit{#1}}
\newcommand{\Mgr}{M_{\operatorname{gr}}}
\newcommand{\Mlex}{M_{\operatorname{lex}}}
\newcommand{\Mlexgr}{M_{\operatorname{lex/gr}}}

\newcommand{\bzero}{\textcolor{blue}{\boldsymbol 0}}
\newcommand{\gzero}{\textcolor{green}{\boldsymbol 0}}
\newcommand{\bone}{\textcolor{blue}{\boldsymbol 1}}
\newcommand{\rzero}{\textcolor{red}{\boldsymbol 0}}
\newcommand{\rone}{\textcolor{red}{\boldsymbol 1}}
\newcommand{\rbullet}{\textcolor{red}{\bullet}}
\newcommand{\bbullet}{\textcolor{blue}{\bullet}}
\newcommand{\bA}{\boldsymbol{\mathcal{A}}}
\newcommand{\bB}{\boldsymbol{\mathcal{B}}}
\newcommand{\bC}{\boldsymbol{\mathcal{C}}}
\newcommand{\fA}{\mathfrak{A}}
\newcommand{\fB}{\mathfrak{B}}
\newcommand{\fC}{\mathfrak{C}}
\newcommand{\barfatslash}{\overline{ \fatslash \,}}

\newcommand{\dimcob}{\operatorname{dim}_{\operatorname{Cob}}}
\newcommand{\dimkom}{\operatorname{dim}_{\operatorname{Kom}}}

\SetSymbolFont{stmry}{bold}{U}{stmry}{m}{n}

\newtheorem{theorem}{Theorem}[section]
\newtheorem{corollary}[theorem]{Corollary}
\newtheorem{lemma}[theorem]{Lemma}
\newtheorem{conjecture}[theorem]{Conjecture}

\newtheorem{numericalresult}[theorem]{Numerical Result}

\usepackage{amsthm}
\usepackage{thmtools, thm-restate}

\declaretheorem[name=Proposition,numberlike=theorem]{proposition}

\theoremstyle{definition}

\theoremstyle{remark}
\newtheorem{example}[theorem]{Example}
\newtheorem{remark}[theorem]{remark}

\usepackage[nottoc,numbib]{tocbibind}
\usepackage{csquotes}
\usepackage[doi=false,isbn=false,url=false,style=alphabetic,sorting=nyt]{biblatex}

\title{Morse matchings and Khovanov homology of 4-strand torus links}
\author{Tuomas Kelomäki}

\begin{document}

\maketitle

\begin{abstract}
Given a link or a tangle diagram, we define algorithmic Morse theoretic simplifications on their Khovanov homology. In contrast to Bar-Natan's scanning algorithm, the cancellations are postponed until the end and performed in one go. Although our novel approach is computationally inferior to Bar-Natan's algorithm, it side-steps the need for a large amount of iterations, making it more fitting for theoretical analysis. Our main application is towards integral Khovanov homology of 4-strand torus links, for which we compute non-trivial Khovanov homology groups in all homological degrees and find an abundance of $4$-torsion. At the limit $T(4,\infty)$, our computations agree with a conjecture of Gorsky, Oblomkov and Rasmussen. For finite $n$, we use the  $\lambda$-invariant of Lewark, Marino and Zibrowius to derive lower bounds on proper rational Gordian distances from $T(4,n)$. 
\end{abstract}

\tableofcontents

\newpage

\section{Introduction}

Khovanov homology is a powerful link invariant which categorifies the Jones polynomial \cite{khovanov1999categorification}. Due to its combinatorial nature, it is readily computable for every link diagram and nowadays there are a number of computer programs which can compute Khovanov homology of diagrams with up to 50-100 crossings within minutes \cite{javakh-v2} \cite{Khocasoftwarearticle} \cite{knotjob} \cite{khtpp}.
The good computability of Khovanov homology is not just an algorithmic feat: a computer calculation played a crucial role in proving that the Conway knot is not smoothly slice \cite{Piccirillo2020} and the efficient  programs open up possibilities of discovering a counterexample to the smooth 4-dimensional Poincaré conjecture \cite{Freedman2010}. 

Despite the fact that for any `small' link $L$ its Khovanov homology, denoted $\operatorname{Kh}(L)$, is easily obtainable, the list of infinite link families  whose Khovanov homology is known remains short. Noteworthy infinite families for which Khovanov homology is well understood include alternating and quasi-alternating links \cite{LeeSpectralSequence}, \cite{manolescu2008khovanovknotfloerhomologies}, $T(2,m)$ and $T(3,m)$ torus links \cite{khovanov1999categorification},      \cite{TurnerSpectralSequenceQcoef}, \cite{Gillam2012},\cite{benheddi2017khovanov}, \cite{chandler_lowrance_sazdanović_summers_2022},  3-strand pretzel links \cite{KhovanovHomologyOfPrezelsRevisited} and closures of $3$-braids \cite{schuetz2025khovanovhomology3braids}. Especially significant attention has been directed to torus links $T(n,m)$. Stošić computed some of the lowest and the highest nontrivial Khovanov homology groups of certain torus links \cite{Stosic2007}, \cite{Stoi2009}.  Furthermore, he showed that Khovanov homologies admit a stable limit $\lim_{m\to \infty} \operatorname{Kh} (T(n,m))$ when renormalized correctly. These limits $\operatorname{Kh} (T(n,\infty))$ were given conjectural descriptions as homologies of explicit Koszul complexes by Gorsky, Oblomkov and Rasmussen \cite{Gorsky2013}.



While a complete picture of Khovanov homology of torus links remains elusive, some related invariants have been determined for all $T(n,m)$. The decategorification, Jones polynomial, admits a closed form \cite{Isidro1993} and the triply graded Khovanov-Rozansky homology is also known for positive $n,m$ \cite{hogancamp2019toruslinkhomology}. Recently, the Lee degeneration of Khovanov homology was also computed for all $T(n,m)$ \cite{Ren2024}. The triply graded homology and the Lee degeneration are connected to Khovanov homology via a spectral sequences $\operatorname{HHH}(L) \Rightarrow\operatorname{Kh}(L)$ \cite{Rasmussen2015} and $\operatorname{Kh}(L) \Rightarrow \mathcal H_{\operatorname{Lee}}(L)$ \cite{Rasmussen2010}. 

In this article we do not aim to compute Khovanov homology using any of the spectral sequences involving it. Instead, we use a more elementary, combinatorial tool called discrete Morse theory. Originally, Forman developed a discretization of the classical Morse theory as a method to simplify cell complexes up to homotopy \cite{FormanCellComplexMorseTheory}. Later Sköldberg formulated it purely in algebraic terms \cite{SkoldbergAlgebraicMorseTheory} which is more suitable for our purposes. Roughly speaking, Sköldberg's algebraic discrete Morse theory is Gaussian elimination for chain complexes employed at a larger scale. Assuming that certain combinatorial conditions are met, one can use the theory to simultaneously cancel many pairs of direct summands, each connected by an isomorphism. This way of packaging the data enables one to avoid a large number of iterations, which would be present, if one were to perform Gaussian elimination naively.

Algebraic discrete Morse theory has been previously employed for knot homologies only very recently. In  \cite{Maltoni2024} Maltoni used discrete Morse theory to reduce the size of Rouquier complexes of braids which define the triply graded homology. Banerjee, Chakraborty and Das wielded it on a model of Khovanov homology \cite{WEHRLI2008} which is based on the Tait graph of the link diagram \cite{banerjee2025spanningtreemodelkhovanov}. The author also applied it in his previous work to compute Khovanov homologies of certain families of links \cite{kelomaki2024discretemorsetheorykhovanov} and this article aims to unify those ideas.

In the original article, Khovanov homology was obtained from a large, explicit chain complex of $\mathbb Z$-modules. A few year later, partly motivated by the search for an efficient algorithm, Bar-Natan formulated a local theory for Khovanov homology of tangles \cite{BarNatanKhovanovTangles}. In this paper, we employ discrete Morse theory on Bar-Natan's local complexes over cobordism categories with a special focus on braids. Given a Morse presentation of a tangle diagram, we introduce two algorithimically defined ansätze of Morse matchings on the tangle complex: \term{a lexicographic matching} $\Mlex$ and \term{a greedy matching} $\Mgr$. For every tangle diagram, we show that $\Mlex$ is Morse matching but strangely, there is exactly one counter-example in a database of 2977 braid diagrams, for which $\Mgr$ fails to generate a meaningful Morse complex. Nevertheless, $\Mgr$ seems to be a more useful tool than $\Mlex$ when performing concrete homology calculations for infinite link families.

\begin{figure}[ht]
    \centering
    \input{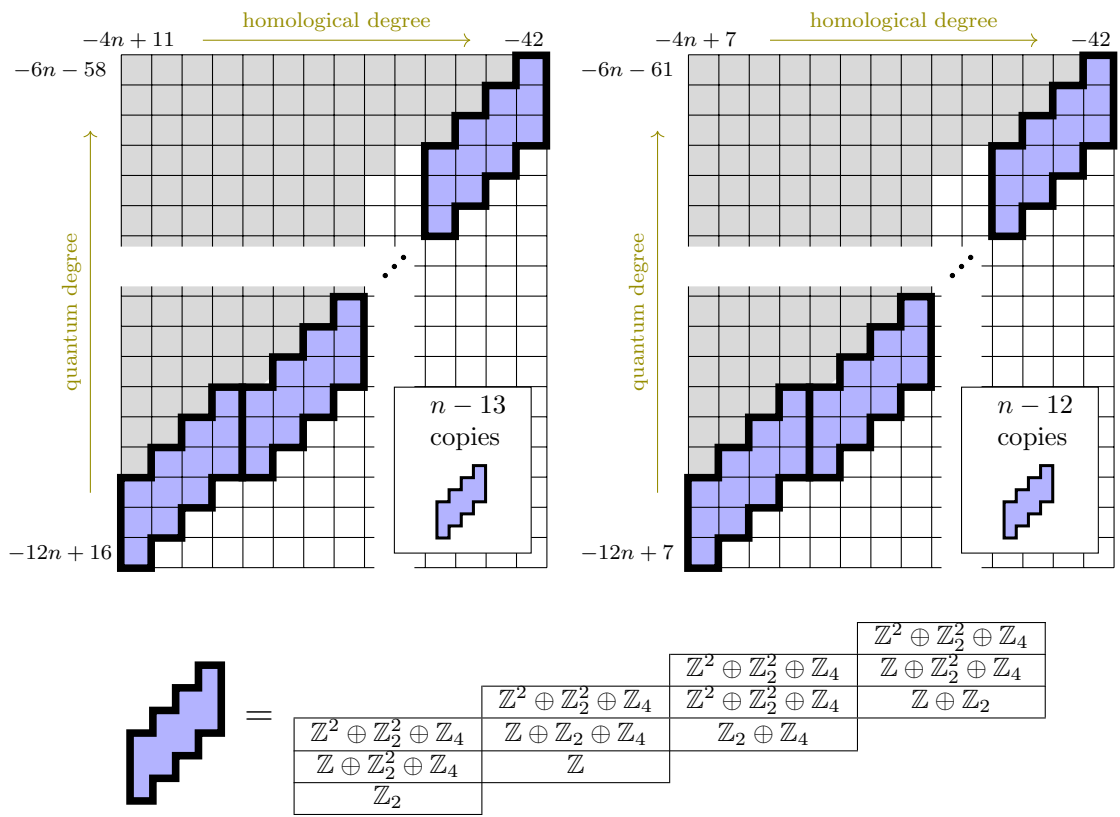}
    \caption{For $n\geq 14$, the unreduced integer Khovanov homology $\operatorname{Kh}^{i,j}(T(4,-2n))$ (top left) and  $\operatorname{Kh}^{i,j}(T(4,-2n-1))$ (top right) and the reoccurring block (bottom). The homology remains unknown in the gray area above the purple blocks and it vanishes below them. More precisely, the top left table describes $\operatorname{Kh}^{i,j}(T(4,-2n))$ only when $i\in [-4n+11,-42]$ and $j\leq \frac{3}{2}i-6n+5$. The top right table describes $\operatorname{Kh}^{i,j}(T(4,-2n-1))$ when $i\in [-4n+7,-42]$ and $j\leq \frac{3}{2}i-6n+2$. Outside these homological degrees $\operatorname{Kh}^{i,j}(T(4,-2n))$ and $\operatorname{Kh}^{i,j}(T(4,-2n-1))$ are described completely in Figures \ref{Figure: T4 homology lowest degrees} and \ref{Figure: T4 homology highest degrees} in Appendix A.}
    \label{Figure: T4 middle}
\end{figure}

Our main application is towards Khovanov homology of negative\footnote{
Negative torus links are mirror images of positive ones and Khovanov homology of a mirror link can easily be taken from the dual chain complex. We follow the conventions of \cite{chandler_lowrance_sazdanović_summers_2022} except that we enforce that $T(4,-n)$ has negative crossings and thus its Khovanov homology is supported on the non-positive homological degrees for $n>0$.
}
4-strand torus links $T(4,-n)$ for which we show that $\Mgr$ is a Morse matching and thus can be applied. Working explicitly through graph theory and combinatorics, we obtain repetitions on the level of open braid Morse complexes which, when taking the braid closure, yield:
\begin{restatable}{theorem}{homologyRecThm}\label{Theorem: recursions for T4 homology}
    The unreduced and reduced Khovanov homology groups of negative 4-strand torus links admit recursions: 
    $$
    \operatorname{Kh}^{i,j}\big(T(4,-n)\big) \cong \operatorname{Kh}^{i,j-12}\big(T(4,-n-4)\big), \quad  \overline{\operatorname{Kh}}^{i,j}\big(T(4,-n)\big) \cong \overline{\operatorname{Kh}}^{i,j-12}\big(T(4,-n-4)\big)
    $$
    for all $n,i,j\in \mathbb Z$ such that $n\geq 0$ and $-2n+2i-j\geq 14$.
    They also admit recursions:
    $$
    \operatorname{Kh}^{i,j}\big(T(4,-n)\big) \cong \operatorname{Kh}^{i-8,j-24}\big(T(4,-n-4)\big), \quad  \overline{\operatorname{Kh}}^{i,j}\big(T(4,-n)\big) \cong \overline{\operatorname{Kh}}^{i-8,j-24}\big(T(4,-n-4)\big)
    $$
    for all $n,i,j\in \mathbb Z$ such that $n\geq 0$ and  $-9n+4i-3j\geq 41$.
\end{restatable}
Combining Theorem \ref{Theorem: recursions for T4 homology} with some vanishing results and a hefty amount of computer base cases for induction, we obtain:
\begin{theorem}\label{Theorem: T4 homology figures}
    For $n\geq 28$, the integral Khovanov homology of negative 4-strand torus links $\operatorname{Kh}^{i,j}(T(4,-n))$ of a certain range  are displayed in Figures \ref{Figure: T4 middle}, \ref{Figure: T4 homology lowest degrees} and \ref{Figure: T4 homology highest degrees}. Figure \ref{Figure: T4 middle} partially describes the middle homological degrees whereas Figures \ref{Figure: T4 homology lowest degrees} and \ref{Figure: T4 homology highest degrees} fully describe the lowest and highest homological degrees respectively.
\end{theorem}

Our results in Figure \ref{Figure: T4 middle} also give out  infinitely many $\mathbb Z /4 \mathbb Z$ torsion groups in the limit $\operatorname{Kh} (T(4,\infty))$ and arbitrarily many $\mathbb Z /4 \mathbb Z$  given  $T(4,n)$ with large $n$. Over the integers, similar $4$-torsion is also present in the conjecture of Gorsky, Oblomkov and Rasmussen and when taken with $\mathbb F_2$ coefficients our calculations for $\operatorname{Kh} (T(4,\infty); \mathbb F_2)$ agree with the conjecture. In addition to the two recursion of Theorem \ref{Theorem: recursions for T4 homology}, our computations also suggest a third one which would be able to completely describe $\operatorname{Kh}^{i,j} (T(4,n))$ for all $i,j,n$. While we do have a rigorous construction for the third correspondence on the level of chain spaces of the open braid complexes, we are unable to prove that it commutes with the differentials. Hence we do not manage to construct a third recursion on the level of homology groups.


Not long after Khovanov's original work, Lee showed that Khovanov homology carries an endomorphism which is also a link invariant \cite{LeeSpectralSequence}. The endomorphism takes slightly different forms depending on the coefficients and whether one considers the unreduced or the reduced Khovanov homology \cite{Rasmussen2010}, \cite{BarNatanKhovanovTangles}, \cite{doi:10.1142/S0218216506005111}. A unified approach to this array of homology theories can be taken by associating to a knot $K$ a homotopy class of chain complexes $\operatorname{BN}(K)$ of graded $\mathbb Z [G]$-modules  where $G$ is a formal variable \cite{zbMATH05118580}. From this $H(\operatorname{BN}(K) \otimes_{\mathbb Z[G]} \mathbb Z[G] / (G) )$ recovers the reduced Khovanov homology $\overline{ \operatorname{Kh}}(K)$ and the unreduced  one can be also retrieved as well. 
Although  we are unable to completely compute the whole Khovanov homology of 4-strand torus knots, let alone  $\operatorname{BN}(T(4,2n+1))$, we manage to show in Proposition \ref{Proposition: wierd piece in ZG homology of T4} that an indecomposable subcomplex of the form:    
\begin{equation}\label{Diagram: intro ZG piece}
\begin{tikzcd}[row sep=0.2em, column sep=3em]
                                                  & {\mathbb Z[G]} \arrow[rd, "2"]  &                \\
{\mathbb Z[G]} \arrow[ru, "G^2"] \arrow[rd, "2G"'] & \oplus                          & {\mathbb Z[G]} \\
                                                  & {\mathbb Z[G]} \arrow[ru, "-G"'] &               
\end{tikzcd}
\end{equation}
splits from $\operatorname{BN}(T(4,2n+1))$ in the highest non-trivial homological degrees when $n\geq 2$.

The majority of geometric applications of Khovanov homology stem from the Rasmussen $s$-invariants. These integers $s_{\mathbb F}(K)$ can be recovered from quantum grading of the free part of $\operatorname{BN}(K)\otimes_{\mathbb Z [G]} \mathbb F [G]$. In Rasmussen's original paper, $s_{\mathbb Q}(K)$ was computed for all positive knots $K$  which gave a combinatorial proof for the slice genus and unknotting numbers of all torus knots $T(n,m)$.

Another source of geometric applications of Khovanov homology come from the $G$-torsion part of $\operatorname{BN}(K)$. In the first edition of these $G$-torsion invariants, Alishahi and Dowlin showed that $\mathfrak u_{\mathbb F_2}(K)  \leq u(K)$ where $\mathfrak u_{\mathbb F_2}(K)$ is the maximum order of $G$-torsion in $\operatorname{BN}(K)\otimes_{\mathbb Z [G]} \mathbb F_2 [G]$ and $u(K)$ is the unknotting number of $K$  \cite{zbMATH07178864}, \cite{zbMATH07005602}. Iltgen, Lewark and Marino refined  $\mathfrak u_{\mathbb F_2}(K)$ to $\lambda(K)$ which works over the integers and includes the information carried by the chain homotopy type of $\operatorname{BN}(K)$. Moreover, they showed that $\lambda$ defines pseudo metric on knots and that $\lambda(K_1,K_2)\leq u_q(K_1,K_2)$ where $u_q$ denotes the rational Gordian distance, that is, the minimum number of proper rational tangle replacements needed to go from $K_1$ to $K_2$. The invariant $\lambda$ was further refined in \cite{lewark2024khovanovhomologyrefinedbounds} to accommodate for grading shifts. 

Since $s_{\mathbb F}(K)$ has already been computed for torus knots, we will rely on the $\lambda$ to make a modest new geometric application of our homology computations: The following proposition, which mimics Theorem 1.4 and Proposition 4.6 from \cite{lewark2024khovanovhomologyrefinedbounds}, can be obtained as a corollary of (\ref{Diagram: intro ZG piece}).
\begin{restatable}{proposition}{ratDistProp}\label{Proposition: rational distance 2 if not positive enough} 
    Suppose $K$ is a knot which admits a diagram with $c_+$ positive crossings. If $n\geq 2$ and $c_+\leq 4n$  then $u_q(T(4,2n+1),K)\geq 2$.  Distinct $4$-strand torus knots are at least $2$ rational replacements from each other: $u_q(T(4,n),T(4,m))\geq 2$ whenever $n,m\in 2\mathbb Z +1$, $n\neq m$ and $|n|,|m|\geq 5$.   
\end{restatable}
To prove (\ref{Diagram: intro ZG piece}) and thus Proposition \ref{Proposition: rational distance 2 if not positive enough} we first computed the reduced non-equivariant Khovanov homology $\overline{\operatorname{Kh}} (T(4,2n+1))$ using Theorem \ref{Theorem: recursions for T4 homology}. Since the homology is reasonably `thin' in the highest degrees, we can guess the $G$-action on the homology, essentially by using the fact that the reduced Lee spectral sequence converges to $0$ in those degrees. It is noteworthy that this tactic would not have sufficed had we done our main Khovanov homology computations over any one field, instead of over the integers. With field coefficients, the necessary $G$-action for $\lambda$ could not be determined from the non-equivariant homology alone.








\textbf{Structure of the paper.} In Section \ref{Section: Preliminaries} we introduce discrete Morse theory and our encoding of tangles and their Khovanov complexes. Section \ref{Section: lex and morse matchings} lays out our protagonists $\Mlex$ and $\Mgr$. Almost all results of this paper are contained in Section \ref{Section: Main results} where homological recursions are derived from the open braid complex ones. The proofs of these braid recursions take up the next two sections. 
Finally in Section \ref{Section: numerical evidence}, we numerically examine the effectiveness of our Morse theoretic simplifications and display the curious braid diagram, where $\Mgr$ breaks down.

\textbf{Computer assistance in this paper.}
The following computer tools have been used in deriving the main mathematical results:
\begin{itemize}
    \item \texttt{Khoca}\cite{Khocasoftwarearticle}, was used to compute the reduced and unreduced integer Khovanov homologies which played the role of base cases in the proofs of Theorem \ref{Theorem: T4 homology figures} and Proposition \ref{Proposition: wierd piece in ZG homology of T4}.
    
    \item Simple custom made \texttt{Python} program \texttt{basecases.py} was used to verify base cases for inductions in Lemmas \ref{Lemma: minimal tA,tC}, \ref{Lemma: t value implies W} and \ref{Lemma: q goes up at most 3 in a path}. It evaluates rational valued functions on specific sets of formal strings. \texttt{Python} was also used to verify that Theorem \ref{Theorem: T4 homology figures} holds for $n=28,\dots, 82$ before the induction kicks in.
    
    \item \texttt{Lean}'s linarith tactic \cite{Moura2021} was used to deduce linear inequalities from sets of linear inequalities in the proofs of Theorem \ref{Theorem: T4 homology figures}, Corollary \ref{Corollary: bounds on Cn dimensions}, Proposition \ref{Proposition: wierd piece in ZG homology of T4} and Lemma \ref{Lemma: q goes up at most 3 in a path}.

    \item \texttt{Mathematica} was used in Proposition \ref{Proposition: GOR-comparison over F_2} to extract power series coefficients from an expansion of a rational function related to the the GOR-conjecture.
    \end{itemize}
Additionally to obtain numerical results in Section \ref{Section: numerical evidence}, a software written by the author \texttt{braidalgo.py}, which computes the greedy matching for any braid diagram, and \texttt{kht++} \cite{khtpp} were utilized. The data of the base cases and the aforementioned code can be found in \cite{computerCodeAndDataFor4BraidPaper}. 


\textbf{Acknowledgments.} The author was supported by the Väisälä Fund of the Finnish Academy of Science and Letters.
The author would like to thank Gregory Arone, Jouko Kelomäki, Oscar Kivinen, Kalle Kytölä, Lukas Lewark, Milo Orlich, Dirk Schütz and Claudius Zibrowius for useful discussions and technical and moral support.

\section{Preliminaries} \label{Section: Preliminaries}
In this section we lay out the main tool, discrete Morse theory, and our conventions for local Khovanov complexes.

\subsection{Algebraic discrete Morse theory}
Let $\mathbf C$ be an additive category. A \term{based} chain complex $C$ is a chain complex over $\mathbf C$ with a fixed \term{basis}, that is, fixed direct sum decomposition $C^i= \bigoplus_{ j\in J_i}C^i_j$  on every \term{chain space} $C^i$. \term{A matrix element} of a based complex $C$ is 
$$
d_{C^{i+1}_k,C^{i}_j}\colon C^i_j \to C^{i+1}_k, \quad d_{C^{i+1}_k,C^{i}_j}= \pi d^i \iota
$$
where $\iota $ and $\pi$ are the canonical inclusion and projection morphisms associated to the direct sums. A based complex induces a simple directed graph $G(C)=(V,E)$ whose vertices are the summands $C_j^i$ and directed edges are non-zero matrix elements. Choosing a subset of edges $M\subset E$ generates another directed graph $G(C,M)=(V',E')$ by reversing the edges of $M$, that is, $V'=V$ and
$$
E'=(E\setminus M) \cup \{ b\to a \mid (a\to b)\in M \}. 
$$
We call the set of edges $M$ a \term{Morse matching} on $G(C)$, if 
\begin{enumerate}
    \item $M$ is finite.\label{Condition: finiteness condition of Morse}
    \item $M$ is a matching, i.e., its edges are pairwise non-adjacent. \label{Condition: matching condition of Morse}
    \item For every edge $f\in M$, the corresponding matrix element $f$ is an isomorphism in $\mathbf C$.\label{Condition: isomorphism condition of Morse}
    \item $G(C, M )$ has no directed cycles. \label{Condition: acyclicity condition of Morse}
\end{enumerate}

We often do not distinguish between vertices of $G(C)$ and $G(C,M)$ and the summands of $C$ and refer to all of them as \term{cells} and denote the collection of all of them by $\operatorname{cells}(C)$. A cell $x$ is called \term{unmatched} if $x$ is not the domain or the codomain of any $f\in M$. Abusing the notation slightly, we can consider the graph $G(C,M)$ as a category, whose morphisms are paths and identity morphisms are paths of length $0$. This allows one to define \term{a remembering functor} $R \colon G(C,M) \to \mathbf C$ which maps an edge $f\in E\setminus M$ to the corresponding matrix element $f$ in $\mathbf C$ and a reversed edge $g$ of $M$ to $-g^{-1}$ in $\mathbf C$. (Naturally, the functor $R$ depends heavily on $C$ and $M$.) 

We can now define the based \term{Morse complex} $(M C, \partial)$ whose basis is given by the unmatched cells of $C$. The matrix elements of $MC$ are given by
$$
\partial_{b,a}= \sum_{p\in \{ \text{paths: } a\to b\} } R(p)
$$
which completely define the differentials $\partial^i$. As a sanity check, the empty matching $M=\emptyset$ recovers the original complex $C=MC$ and matching of one edge $M=\{f\}$ recovers Gaussian elimination of chain complexes, Lemma 3.2 from \cite{FastKhovavnovComputations}. A consequence of Condition \ref{Condition: matching condition of Morse} is that only those paths which zig-zag between two adjacent homological degrees will end up contributing to $\partial$.

\begin{theorem}[Sköldberg \cite{SkoldbergAlgebraicMorseTheory}]\label{Theorem: discrete Morse theory}
If $M$ is a Morse matching on $G(C)$, then $C \simeq MC$.    
\end{theorem}

Although we will not need them in this article, the explicit chain maps of Theorem \ref{Theorem: discrete Morse theory} can also be written as sums of paths: the maps $f\colon C\to MC$ and $g\colon MC \to C$ are defined by their matrix elements with 
$$
\pi_b f \iota_a=  \sum_{p\in \{ \text{paths: } a \to b\} } R(p), \qquad \pi_{d} g \iota_{c}=  \sum_{p\in \{ \text{paths: } c\to d\} } R(p) 
$$
where $\iota_a, \iota_{c}, \pi_{b}$ and $\pi_d$ are again the canonical projection and inclusion morphisms with respect to $C$ and $MC$. The map $g$ could be used to track back a homology cycle of $MC$ to a cycle of $C$. Explicit homology cycles are important in Khovanov homology since they are needed to evaluate the induced homomorphisms which in turn can be used to distinguish exotically embedded smooth surfaces in $B^4$, see \cite{HaydenSundberg+2024+217+246}. Recently this strategy was used in \cite{banerjee2025spanningtreemodelkhovanov}, where discrete Morse theory was employed on a graph theoretic model of Khovanov homology to obtain explicit homology cycles which were used to obtain a new family of exotic surfaces.

\subsection{The Bar-Natan hypercube complex and delooping}

In \cite{BarNatanKhovanovTangles} Bar-Natan constructed from a tangle diagram $T$ with $2n$ endpoints a local Khovanov complex $\llbracket T \rrbracket$ over the category $\operatorname{Mat} (\Cob(2n))$ which we will review here. The building blocks of this category are of the form $c\{m\}$ where $c$ is a crossingless planar diagrams  and $m\in \mathbb Z$ is a quantum grading shift. Furthermore, one adds finite formal direct sums to the mix, so that a generic object in $\operatorname{Mat} (\Cob(2n))$ is written as $\bigoplus_i c_i\{m_i\}$. The morphisms of $\operatorname{Mat} (\Cob(2n))$ are  matrices whose entries are of formal sums of dotted cobordisms. These formal sums of cobordisms are considered up to boundary preserving isotopy, moving the dots around and the ``dotted" relations depicted in Figure \ref{Figure: relations of Cob}. The morphisms are graded, so that if $f\colon c_1\{m_1\} \to c_2\{m_2\}$ is a cobordism, then  $\operatorname{deg}(f)= \chi (f)-2\cdot\#(\text{dots in } f)+n+m_2-m_1$ where $\chi(f)$ is the Euler characteristic of the underlying cobordism $f$. We will only consider morphisms of degree 0.

\begin{figure}[ht]
    \centering    
\begin{tikzpicture}[scale=0.4]
  
  \def\xradius{1};
  \def\yradius{0.5};
  \def\ballsxcoor{-15}
  \def\rectanglexcoor{-9}

      \draw (\ballsxcoor,4) ellipse (\xradius cm and \xradius cm);
    \draw[dashed] plot[domain=-\xradius:\xradius] ({\x +\ballsxcoor},{4+1*\yradius * sqrt(1 - (\x / \xradius)^2)});
  \draw plot[domain=-\xradius:\xradius] ({\x +\ballsxcoor},{4-1*\yradius * sqrt(1 - (\x / \xradius)^2)});
\draw (\ballsxcoor+2.5,4) node { = 0};

      \draw (\ballsxcoor,0) ellipse (\xradius cm and \xradius cm);
    \draw[dashed] plot[domain=-\xradius:\xradius] ({\x +\ballsxcoor},{1*\yradius * sqrt(1 - (\x / \xradius)^2)});
  \draw plot[domain=-\xradius:\xradius] ({\x +\ballsxcoor},{-1*\yradius * sqrt(1 - (\x / \xradius)^2)});
\draw (\ballsxcoor+2.5,0) node { = 1};

  \coordinate (A) at (\rectanglexcoor,0.5);
  \coordinate (B) at (\rectanglexcoor+3,1);
  \coordinate (C) at (\rectanglexcoor+3,4);
  \coordinate (D) at (\rectanglexcoor,3.5);
  
  \draw (A) -- (B) -- (C) -- (D) -- cycle;

    \filldraw (\rectanglexcoor+1.5,1.75) circle (3pt);
    \filldraw (\rectanglexcoor+1.5,2.5) circle (3pt);
    \draw (\rectanglexcoor+4.5,2.5) node {= 0};

    \filldraw (\ballsxcoor,0) circle (3pt);

  \draw (0,4) ellipse (\xradius cm and \yradius cm);
  \draw (4,4) ellipse (\xradius cm and \yradius cm);

\draw[dashed] plot[domain=-\xradius:\xradius] ({\x},{1*\yradius * sqrt(1 - (\x / \xradius)^2)});
\draw plot[domain=-\xradius:\xradius] ({\x},{-1*\yradius * sqrt(1 - (\x / \xradius)^2)});

  \draw plot[domain=-\xradius:\xradius] ({\x},{4 - 2.5*\yradius * sqrt(1 - (\x / \xradius)^2)});
  \draw plot[domain=-\xradius:\xradius] ({\x + 4},{4 - 2.5*\yradius * sqrt(1 - (\x / \xradius)^2)});

  \draw plot[domain=-\xradius:\xradius] ({\x},{2.5*\yradius * sqrt(1 - (\x / \xradius)^2)});
  \draw plot[domain=-\xradius:\xradius] ({\x + 4},{2.5*\yradius * sqrt(1 - (\x / \xradius)^2)});
  \draw[dashed] plot[domain=-\xradius:\xradius] ({\x + 4},{1*\yradius * sqrt(1 - (\x / \xradius)^2)});
  \draw plot[domain=-\xradius:\xradius] ({\x + 4},{-1*\yradius * sqrt(1 - (\x / \xradius)^2)});
  
  \draw (2,2) node { +};

    \draw (6,2) node { =};

    \draw (8.5,4) ellipse (\xradius cm and \yradius cm);
  \draw[dashed] plot[domain=-\xradius:\xradius] ({\x + 8.5},{1*\yradius * sqrt(1 - (\x / \xradius)^2)});
  \draw plot[domain=-\xradius:\xradius] ({\x + 8.5},{-1*\yradius * sqrt(1 - (\x / \xradius)^2)});

    \draw (7.5,0) -- (7.5,4);
    \draw (9.5,0) -- (9.5,4);
    
    \filldraw (4,3) circle (3pt);
    \filldraw (0,1) circle (3pt);
  
\end{tikzpicture}

    \caption{The dotted relations of $\operatorname{Mat} (\Cob(2n))$.}
    \label{Figure: relations of Cob}
\end{figure}

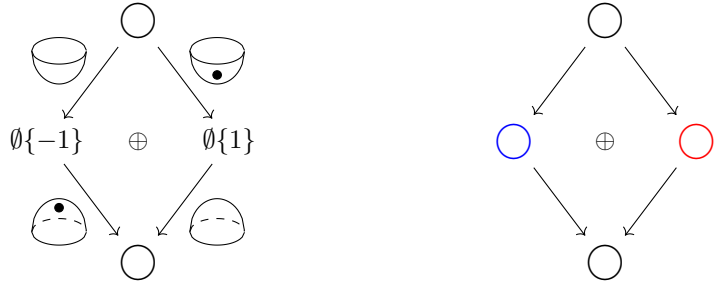
\begin{figure}[ht]
  \centering
  \begin{subfigure}[b]{0.4\textwidth}
  \centering
    \begin{tikzpicture}[scale=0.8]

  \node (A) at (0, 4) {\Large$\bigcirc$};
\node (B) at (1.5, 2) { $\emptyset \{1\}$};
  \node (C) at (-1.5, 2) { $\emptyset \{-1\}$};
  \node (D) at (0, 0) {\Large$\bigcirc$};

\node (E) at (0, 2) {$\oplus$};

  \draw[->] (A) -- (B);
  \draw[->] (A) -- (C);
  \draw[->] (B) -- (D);
  \draw[->] (C) -- (D);






  \def\xradius{0.45};
  \def\yradius{0.22};

\def\a{1.3}
\def\b{3.5}

\def\c{-1.3}
\def\d{3.5}

\def\e{1.3}
\def\f{0.5}

\def\g{-1.3}
\def\h{0.5}

\draw plot[domain=-\xradius:\xradius] ({\a+\x},{\b+1*\yradius * sqrt(1 - (\x / \xradius)^2)});
\draw plot[domain=-\xradius:\xradius] ({\a+\x},{\b-1*\yradius * sqrt(1 - (\x / \xradius)^2)});
\draw plot[domain=-\xradius:\xradius] ({\a+\x},{\b-2.5*\yradius * sqrt(1 - (\x / \xradius)^2)});
\filldraw (\a,\b-0.4) circle (2pt);

\draw plot[domain=-\xradius:\xradius] ({\c+\x},{\d+1*\yradius * sqrt(1 - (\x / \xradius)^2)});
\draw plot[domain=-\xradius:\xradius] ({\c+\x},{\d-1*\yradius * sqrt(1 - (\x / \xradius)^2)});
\draw plot[domain=-\xradius:\xradius] ({\c+\x},{\d-2.5*\yradius * sqrt(1 - (\x / \xradius)^2)});

\draw[dashed] plot[domain=-\xradius:\xradius] ({\e+\x},{\f+1*\yradius * sqrt(1 - (\x / \xradius)^2)});
\draw plot[domain=-\xradius:\xradius] ({\e+\x},{\f-1*\yradius * sqrt(1 - (\x / \xradius)^2)});
\draw plot[domain=-\xradius:\xradius] ({\e+\x},{\f+2.5*\yradius * sqrt(1 - (\x / \xradius)^2)});

\draw[dashed] plot[domain=-\xradius:\xradius] ({\g+\x},{\h+1*\yradius * sqrt(1 - (\x / \xradius)^2)});
\draw plot[domain=-\xradius:\xradius] ({\g+\x},{\h-1*\yradius * sqrt(1 - (\x / \xradius)^2)});
\draw plot[domain=-\xradius:\xradius] ({\g+\x},{\h+2.5*\yradius * sqrt(1 - (\x / \xradius)^2)});
\filldraw (\g,\h+0.4) circle (2pt);

\end{tikzpicture}

  \end{subfigure}
  \begin{subfigure}[b]{0.4\textwidth}
  \centering
    \begin{tikzpicture}[scale=0.8]

  \node (A) at (0, 4) {\Large$\bigcirc$};
\node (B) at (1.5, 2) { \textcolor{red}{\Large$\bigcirc$}};
  \node (C) at (-1.5, 2) { \textcolor{blue}{\Large$\bigcirc$}};
  \node (D) at (0, 0) {\Large$\bigcirc$};

\node (E) at (0, 2) {$\oplus$};

  \draw[->] (A) -- (B);
  \draw[->] (A) -- (C);
  \draw[->] (B) -- (D);
  \draw[->] (C) -- (D);

\end{tikzpicture}

  \end{subfigure}
    \caption{For a diagram which consists of only a single circle: delooping isomorphism $\Psi$ and its inverse $\Psi^{-1}$.  The cobordism of the matrix elements of $\Psi$ and $\Psi^{-1}$ are displayed on the left and  our diagrammatic convention for the delooped objects is shown on the right. }
    \label{Figure: delooping isomorphism}
\end{figure}

Due to the relations imposed on cobordisms, the category $\operatorname{Mat} (\Cob(2n))$ admits a local delooping isomorphism $\Psi$ depicted in Figure \ref{Figure: delooping isomorphism}. Iteratively removing loops with $\Psi$, one obtains a global isomorphism 
$$
\Psi_{c\{m\}} \colon c\{m\} \to \bigoplus_{K\subset  L_c} c'\{m+2\cdot\# K -\# L_c\}
$$
where $L_c$ is the set of loops in $c$ and $c'$ is the planar diagram of $c$ with all of the loops removed. We denote $\Psi c$ for the codomain of $\Psi_c$ and graphically represent a direct summand of $\Psi c$ by coloring the loops contained in $K$ with red and those in $L_c \setminus K$ blue. By acting diagonally, $\Psi$ can be extended to formal direct sums and for any object $a$.


The standard Bar-Natan hypercube complex $(\llbracket T \rrbracket,d)$ is generated from the tangle diagram $T$ by resolving every crossing \begin{tikzpicture}[baseline={(0,0.05)},scale=0.3]
    \negCrossing{0,0} 
\end{tikzpicture} with either a 0-smoothing \begin{tikzpicture}[baseline={(0,0.05)},scale=0.3]
    \lineRightBendUp{0,0}
    \lineRightBendDown{0,1}
\end{tikzpicture} or a 1-smoothing \begin{tikzpicture}[baseline={(0,0.05)},scale=0.3]
    \lineUpBendRight{0,0}
    \lineUpBendLeft{1,0}
\end{tikzpicture}. 
If a single crossing is $0$-smoothed in $c_0$ and $1$-smoothed in $c_1$ and the diagrams $c_0$ and $c_1$ coincide otherwise, then the matrix element $d_{c_1,c_0}$ is the simplest saddle cobordism between $c_0$ and $c_1$ with a sign $s\in \{-,+\}$ which is explained later. Suppose that $T$ has $n_+$ and $n_-$ positive and negative crossings, see Figure \ref{Morse presentation of tangle} for conventions, and $c$ is a diagram with $k$ number of 1-smoothings, then $c$ will have homological degree $k-n_-$ and quantum degree shift $k+n_+-2n_-$.



In the hypercube complex $(\llbracket T \rrbracket, d)$, all of the morphisms are saddles which in particular means that there are no isomorphisms and discrete Morse theory will not make a dent on $\llbracket T \rrbracket$. Instead we work with $\Psi \llbracket T \rrbracket$ whose the chain spaces are defined as $\Psi(\llbracket T \rrbracket^i)$ and the differential are pulled back from $d$ and defined as $\Psi d^i \Psi^{-1}$.   
The basis of $\Psi \llbracket T \rrbracket$ consists of diagrams with black line segments from the boundary to the boundary accompanied with red and blue colored circles. The following lemma allows us to classify the matrix elements of $\Psi \llbracket T\rrbracket$ in a local manner. Its proof is a case-by-case application of the relations in Figure \ref{Figure: relations of Cob} and we omit it.

\begin{figure} 
    \centering
\input{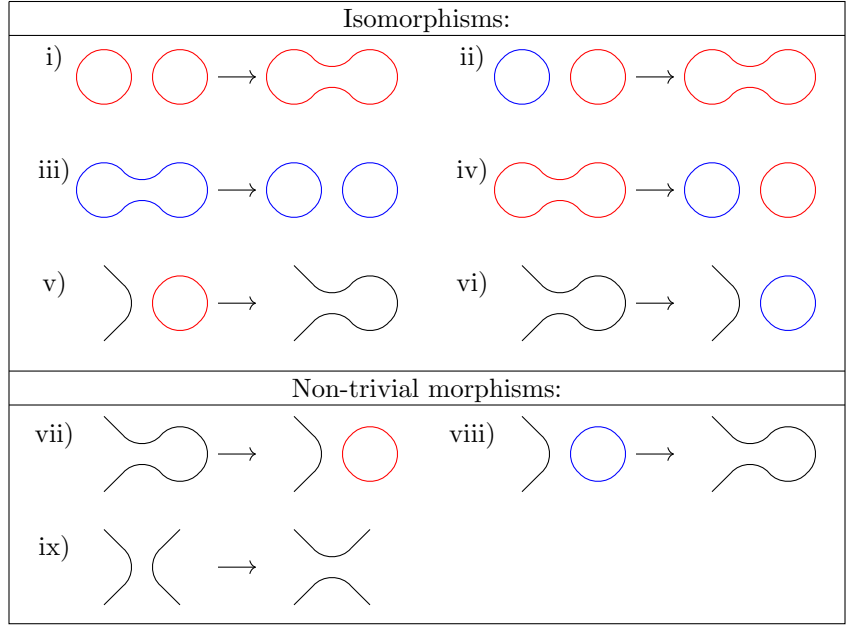}

    \caption{A local dictionary for matrix elements of $\Psi \llbracket T \rrbracket$. The pictures i)-vi) represent identity cobordisms. The pictures vii) and viii) correspond to adding a dot to the black component  whereas ix) is a saddle.}
    \label{Figure: 2d local dictionary}
\end{figure}

\begin{lemma}\label{Lemma: local classification of delooped matrix elements}
    Let $T$ be tangle diagram and suppose $d_{b,a}$ is a matrix element of $\Psi \llbracket T \rrbracket$ where $a,b$  are completely smoothed pictures with red/blue colored loops. If $a$ and $b$ agree everywhere except for a local change displayed in Figure \ref{Figure: 2d local dictionary}, pictures i)-vi), then $d_{b,a}$ is an isomorphism. If $a$ and $b$ agree everywhere except for a local change of pictures vii)-ix), then $d_{b,a}$ is neither a zero morphism nor an isomorphism. Otherwise $d_{b,a}=0$.
\end{lemma}

\begin{figure}
    \centering
    \input{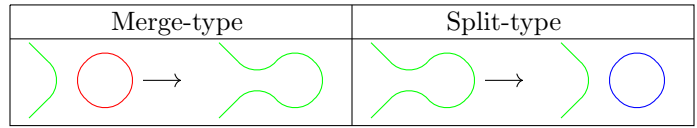}
    \caption{Using green as a placeholder color for red, blue or black, we can divide the isomorphisms of Figure \ref{Figure: 2d local dictionary} into two classes. Types i), ii) and v) amount to the merge-type, whereas iii),  iv) and vi) make up the split-type.}
    \label{Figure: merge type and split type}
\end{figure}

The cells of $\Psi \llbracket T \rrbracket$ are large pictures, but often we only care about some of their local aspects. In particular, when drawing local diagrams we are typically very interested in the colors of some loops, but less concerned whether other components in the global picture form red or blue loops, or whether they are connected to the boundary and are thus colored black. In these cases, we draw the less important components with green as to say that a local phenomenon can be lifted into any of the global pictures of $\Psi \llbracket T \rrbracket$ that are obtained by replacing the green color with red, blue or black individually in each component. As an example, the isomorphisms of Figure \ref{Figure: 2d local dictionary} can be divided into two classes: \term{merge-type} and \term{split-type}, see Figure \ref{Figure: merge type and split type}. 

\subsection{Morse presentations of tangles}
\begin{figure}
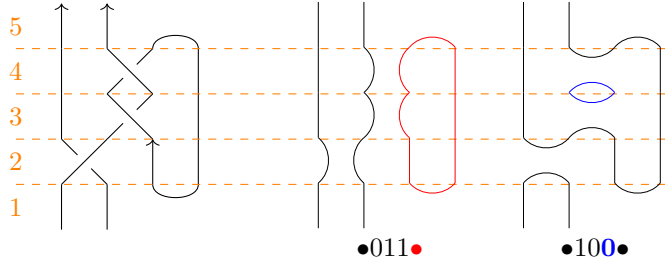

    \centering
\begin{tikzpicture}[scale=0.6]
\lineUp{0,0}
\lineUp{1,0}
\lineRightBendDownMore{2,1}

\posCrossing{0,1}
\lineUpWArrow{2,1}
\lineUp{3,1}

\lineUp{0,2}
\negCrossing{1,2}
\lineUp{3,2}

\lineUp{0,3}
\negCrossing{1,3}
\lineUp{3,3}

\lineUpWArrow{0,4}
\lineUpWArrow{1,4}
\lineRightBendUpMore{2,4}

\draw[orange, dashed] (-1, 1) -- (13.5, 1);

\draw[orange, dashed] (-1, 2) -- (13.5, 2);

\draw[orange, dashed] (-1, 3) -- (13.5, 3);
\draw[orange, dashed] (-1, 4) -- (13.5, 4);

\node at (-1,0.5) {\textcolor{orange}{1}};

\node at (-1,1.5) {\textcolor{orange}{2}};
\node at (-1,2.5) {\textcolor{orange}{3}};
\node at (-1,3.5) {\textcolor{orange}{4}};
\node at (-1,4.5) {\textcolor{orange}{5}};

\node at (5.5,-0.1)[scale=0.6, anchor=south west] {\input{cell_examples/cAbbc0011y}
};
\node at (10,-0.1)[scale=0.6, anchor=south west] {\input{cell_examples/cAbbc010x0}
};
\node at (7.2,-0.4) {$\bullet 011 \rbullet$};

\node at (11.7,-0.4) {$\bullet 10 \bzero \bullet$};

\end{tikzpicture}
    
    \caption{Morse presentation of an oriented tangle with 5 Morse layers. On layer 2 we see a positive crossing and on layers 3 and 4 we see negative ones. On the right we see 2 examples of cells in $\Psi \llbracket T \rrbracket$ with their symbolic encodings.}
    \label{Morse presentation of tangle}
\end{figure}

Topologically, a tangle $T$ is an embedded smooth 1-manifold in $B^3$ with possibly a non-empty boundary which intersects $S^2=\partial B^3$ transversally. Two tangles are equivalent, if they are related by a boundary fixing ambient isotopy.  In this paper, we will work with (highly non-unique) combinatorial descriptions of tangles as \term{Morse presentations}. They consist of $a$ strings at top and $b$ at the bottom with $a+b\in 2\mathbb Z_{\geq 0}$. In between there is a finite number of \term{Morse layers} consisting of either of the two crossings
$\begin{tikzpicture}[baseline={(0,0.05)},scale=0.3]
    \posCrossing{0,0} 
\end{tikzpicture}$,
$\begin{tikzpicture}[baseline={(0,0.05)},scale=0.3]
    \negCrossing{0,0} 
\end{tikzpicture}$, a cap
$\begin{tikzpicture}[baseline={(0,0.05)},scale=0.3]
    \lineRightBendUpMore{0,0} 
\end{tikzpicture}$ or a cup
$\begin{tikzpicture}[baseline={(0,-0.15)},scale=0.3]
    \lineRightBendDownMore{1,0}
\end{tikzpicture}$ 
accompanied with an appropriate number of vertical bars 
$\begin{tikzpicture}[baseline={(0,0.05)},scale=0.3]
    \lineUp{0,0}
\end{tikzpicture}$
to match the connectivities. An oriented tangle diagram and our conventions for positive and negative crossings can be seen in Figure \ref{Morse presentation of tangle}.

To proceed efficiently, we will introduce a symbolic way of encoding summands of $\llbracket T \rrbracket$ and $ \Psi \llbracket T \rrbracket$, where $T$ is a Morse presentation of a tangle. For $\llbracket T \rrbracket$, the Morse layers from bottom to top are written from left to right with symbols $0$ and $1$ denoting smoothings of crossings and $\bullet$ placed whenever the Morse layer is a cap or a cup. A matrix element $d_{b,a}$ of $\llbracket T \rrbracket$ changes $0$-symbol of $a$ to $1$-symbol of $b$ at some index $1\leq i \leq n$. The sign $s$ of $d_{b,a}$ is defined as $(-1)^k$ where $k$ is the number of symbols $1$ in $a$ with indices $j<i$.

For $ \Psi \llbracket T \rrbracket$ the colors of loops need to be stored as well: given a loop $l$ in a cell $a$ of $ \Psi \llbracket T \rrbracket$, we denote 
the highest Morse layer which contains a piece of $l$, by $\max(l)$. The color of the loop $l$ is engraved to symbol of $a$ at index $\max(l)$. Thus for an $n$-layered Morse presentation of a tangle $T$, the summands of $\llbracket T \rrbracket$ and $\Psi \llbracket T \rrbracket$ are \term{symbolically encoded} with words of $\{ 0, 1, \bullet  \}^n$ and $\{ 0,\bzero, \rzero,1,\bone,\rone,\bullet, \bbullet,\rbullet \}^n$ respectively, see Figure \ref{Morse presentation of tangle} for examples. Due to technical reasons, we allow for degenerate tangles with zero Morse layers and $n$ degenerate vertical lines.  For a degenerate tangle $T$ the unique summand of $\llbracket T \rrbracket$ and $\Psi \llbracket T \rrbracket$ is denoted with $e$.



\begin{figure}
\centering
\scalebox{0.9}{
\input{deloopedhypercube}
}
    
    \caption{A morse presentation of a tangle $T$ (left) and its delooped hypercube complex $\Psi \llbracket T \rrbracket $ (right). The complex has two isomorphisms, $f\colon 00\bullet \to 10\bbullet $ and $g\colon 10\rbullet \to 11 \bullet$ which satisfy $i(f)=2$, $i(g)=1$ and $u(f)=u(g)=3$. Both isomorphism get matched in $\Mgr$ and hence $\Mgr\Psi\llbracket T  \rrbracket$ consists of 
    only the picture $10\bullet$ which in turn proves invariance of Khovanov homology under Reidemeister 2 move.
    }    
\end{figure}

\section{Algorithmic matchings for Khovanov complexes}\label{Section: lex and morse matchings}

In this section, we will define a greedy matching on Khovanov complexes of tangles with fixed Morse presentations. In it, the following two maps play a crucial role:
\begin{align*}
i \colon &\{\text{matrix elements of } \Psi \llbracket T \rrbracket\} \to \mathbb Z_{>0} \\
u \colon &\{\text{isomorphism matrix elements of } \Psi \llbracket T \rrbracket\} \to \mathbb Z_{>0}.   
\end{align*}
For a matrix element $a\to b$ of $\Psi \llbracket T \rrbracket$ we define $i(a\to b)$ to be the index of the Morse layer for which a character $0,\bzero$ or $\rzero$ of $a$ changes to a character $1, \bone$ or $\rone$ of $b$.  Denote $L_{a\to b}$ to be the set of red loops in $a$ merged by $a\to b$ and blue loops in $b$ split out by $a\to b$. Notice that by Lemma \ref{Lemma: local classification of delooped matrix elements} we have $L_{a\to b}\neq \emptyset$ if and only if $a\to b$ is an isomorphism. Hence we can define
$
u(a\to b)=\min \{\max(l) \mid l\in L_{a \to b} \}
$ and one should note that $i(a\to b) \leq u(a\to b)$ always holds. Given a (candidate) Morse matching $M$ and an edge $(a\to b)\in G(\Psi \llbracket T \rrbracket,M)$ we are often somewhat loose with the orientations and write $i(a\to b)$ and $u(a\to b)$ regardless of whether $a\to b$ or $b\to a$ actually represents a matrix element of $\Psi \llbracket T \rrbracket$. 

For a fixed tangle $T$ and $j,k\in \mathbb Z_{>0} $ we define $A_{j,k}=i^{-1}\{j\} \cap u^{-1}\{k\}$. Furthermore, we set $N_1=M_{1,1}=A_{1,1}$ and recursively assign for any $j \in \mathbb Z_{>0}$ 
$$
N_{j}=\left\{ (a\to b) \in A_{j,j} \mid
      a,b  \text{ are not matched in any } N_l  \text{ with } l<j   
 \right\} 
$$
and for any pair $j,k \in \mathbb Z_{>0}$
$$
M_{j,k}=\left\{ (a\to b) \in A_{j,k} \middle\vert \begin{array}{c}
      a,b \text{ are not matched in any } M_{l,m}  \\
      \text{with } m<k      \\
     \text{ or with } m=k  \text{ and } j>l
\end{array} \right\}. 
$$
Finally we define \term{the lexicographic matching} to be $\Mlex=\bigsqcup_j N_j$ and \term{the greedy matching}\footnote{
The matching $\Mlex$ is named \textit{lexicographic}, as it seems to move $\rzero$ and $\bone$ characters to the lower Morse layers in an orderly fashion. The matching $\Mgr$ is labeled greedy, since it is a maximal partial matching on $G(\Psi\llbracket T \rrbracket)$ that is, there does not exist any strictly larger partial matching of isomorphisms containing $\Mgr$.
}
to be $\Mgr=\bigsqcup_{j,k} M_{i,j}$.  

It is immediate that both $\Mlex$ and $\Mgr$ satisfy Conditions \ref{Condition: finiteness condition of Morse}-\ref{Condition: isomorphism condition of Morse} of being a Morse matching. In Proposition \ref{Proposition: lexicographic matching is Morse} we prove that $G(\Psi \llbracket T \rrbracket, \Mlex )$ is always acyclic and thus $\Mlex$ a Morse matching for all tangle diagrams. Curiously, the matching $\Mgr$ is not always acyclic, although it often appears to be; in Section \ref{Section: numerical evidence} we find that in a test set of 2977 braid diagrams for exactly one $B$ the graph $G(\Psi \llbracket B \rrbracket, \Mgr )$ contains directed cycles.

Although, $\Mgr$ has the obvious aforementioned caveat that it does not always work, for concrete homology calculations it seems to be a more useful tool than $\Mlex$. In previous work of the author \cite{kelomaki2024discretemorsetheorykhovanov}, explicit Morse matchings were imposed on 2- and 3-strand torus braids. These ad hoc matchings can now be seen as instances of $\Mgr$ and hence their acyclicity proofs (Lemmas 3.2 and 4.1 from \cite{kelomaki2024discretemorsetheorykhovanov}) show that graphs $G(\Psi \llbracket \sigma_1^n\rrbracket,\Mgr )$ and $G(\Psi \llbracket (\sigma_1\sigma_2)^n\rrbracket,\Mgr )$ contain no directed cycles. For 4-strand torus braids, $\Mgr$ also gives an acyclic matching which is proven in Proposition \ref{Proposition: Mgr Morse on T4}. This leads to the following conjecture:
\begin{conjecture}\label{Conjecture: greedy matching is Morse for T(n,m)}
    The greedy matching is a Morse matching for all torus braids, that is, the graphs $G(\Psi \llbracket (\sigma_1 \dots \sigma_{m-1})^n\rrbracket,\Mgr )$ contain no directed cycles for all $n,m$.
\end{conjecture}



Both $\Mlex$ and $\Mgr$ are defined on the whole graph $G(\Psi \llbracket T \rrbracket)$, but in practice one can ignore large portions of $G(\Psi \llbracket T \rrbracket)$ when evaluating $\Mlex$ and $\Mgr$. In this sense they bear resemblance to Bar-Natan's scanning algorithm \cite{FastKhovavnovComputations}, where small ``local" cancellations end up massively reducing the large ``global" complex. To efficiently\footnote{In terms of \textit{computational efficiency}, Bar-Natan's algorithm seems to greatly outperform our use of discrete Morse theory. Hence we aim towards \textit{mathematical efficiency}, by which we mean that we hope to construct an easy-to-use theory and prove theorems using it. 
} construct $\Mlex \Psi \llbracket T \rrbracket$ and  $\Mgr \Psi \llbracket T \rrbracket$, two questions need to  answered:
\vspace{-0.1 cm}
\begin{equation}
\text{\textit{Given a tangle diagram $T$, what are the unmatched cells of $ M_{\operatorname{lex/gr}} \Psi \llbracket T \rrbracket$?}} \label{Question: what are unmatched cells of Mlex/Mgr} 
\end{equation}
\vspace{-0.6 cm}
\begin{equation}
    \text{\textit{Given unmatched cells $a$ and $b$, what are all paths from $a$ to $b$ in $G(\Psi \llbracket T \rrbracket, M_{\operatorname{lex/gr}})$?}} \label{Question: what are all paths from a to b in Mlex/Mgr}
\end{equation}


The cells of $\Psi \llbracket T \rrbracket$ can be represented as formal words of symbols which is why we can also manipulate them as such. To do so, we denote $\mathcal{S}$ as the free monoid generated by symbols $0,\bzero, \rzero,1,\bone,\rone, \bullet, \bbullet,\rbullet$, that is, as a set 
$$
\mathcal S=\bigcup_{n=0}^{\infty} \{ 0,\bzero, \rzero,1,\bone,\rone, \bullet, \bbullet,\rbullet \}^n.
$$
The multiplication in $\mathcal S$ is the concatenation of words and marked with $[.]$, e.g. $0\bone\rbullet. \rone=0\bone\rbullet\rone$ and the neutral element is denoted with $e$. The number of characters in a word $w$ is denoted with $|w|$. 

Given a tangle diagram $T$, we define $\operatorname{cut}_m(T)$ as the tangle consisting of the first $m$ Morse layers of $T$. For a formal set of words $W\subset \operatorname{cells}(\operatorname{cut}_m(T)) \subset \mathcal S$ we denote $\operatorname{expand}(W)$ as the set of all sensible continuations of $W$ under $T$, that is, 
$$
\operatorname{expand}(W)=\big\{w.s\in \operatorname{cells}(\operatorname{cut}_{m+1}(T)) \, \big \vert \  w\in W, \  s\in \{ 0,\bzero, \rzero,1,\bone,\rone,\bullet, \bbullet,\rbullet \} \big\}.
$$
It is easy to see that for any tangle diagram $T$, we have $\operatorname{cells}(\Psi \llbracket T\rrbracket)= \operatorname{expand}\circ \dots \circ \operatorname{expand}(\{e\})$. We can now algorithmically answer Question \ref{Question: what are unmatched cells of Mlex/Mgr} for $\Mlex$ and for $\Mgr$.

\begin{algorithm*}
\caption{Obtaining unmatched cells of $\Mlex \Psi \llbracket T \rrbracket$ and answering Question \ref{Question: what are unmatched cells of Mlex/Mgr}, case $\Mlex$.} \label{Algorithm: unmatched cells of Mlex}
 \textbf{Input} $T$, a Morse presentation of a tangle diagram.\\
  \textbf{Output} $U$, the set of unmatched cells of $\Mlex \Psi\llbracket T \rrbracket$. 

   \begin{algorithmic}[1]
   \State $U\gets \{e \}$ 
  
\For{k=1,\dots, \#\text{Morse layers of $T$}}
    \State $U\gets \operatorname{expand}(U)$ 
    \State Remove all pairs of cells $(a,b)$ from $U$ with $i(a\to b)=u(a\to b)=k$   \label{Mlex algoline remove pairs}
\EndFor
\State \textbf{return} $U$
\end{algorithmic}
\end{algorithm*}
\begin{algorithm*}
\caption{Obtaining unmatched cells of $\Mgr \Psi \llbracket T \rrbracket$ and answering Question \ref{Question: what are unmatched cells of Mlex/Mgr}, case $\Mgr$.} \label{Algorithm: unmatched cells of Mgr}
 \textbf{Input} $T$, a Morse presentation of a tangle diagram.\\
  \textbf{Output} $U$, the set of unmatched cells of $\Mgr \Psi \llbracket T \rrbracket$. 

   \begin{algorithmic}[1]
   \State $U\gets \{e \}$ 
  
\For{k=1,\dots, \#\text{Morse layers of $T$}}
    \State  $U\gets \operatorname{expand}(U)$
    \For{j=k,\dots,1}    
        \State Remove all pairs of cells $(a,b)$ from $U$ with $i(a\to b)=j$ and $u(a\to b)=k$ \label{Mgr algoline remove pairs}
    \EndFor
\EndFor
\State \textbf{return} $U$
\end{algorithmic}
\end{algorithm*}

In a directed acyclic graph, the set of all paths from $a$ to $b$ can be obtained by using an exhaustive depth-first search. 
For a cell $v$ in the graph $G(\Psi \llbracket T \rrbracket )$, the set of edges arriving to $v$ or departing from $v$ can easily be attained, but it is not obvious which ones are contained in $\Mlexgr$, that is, which edges are reversed. Naturally, this information is crucial in constructing the set of all paths from $a$ to $b$. It follows that we can use an exhaustive search to answer Question \ref{Question: what are all paths from a to b in Mlex/Mgr} assuming that we first give an answer to the following question:
\begin{equation}
\text{\textit{Given $a\in \operatorname{cells} (\Psi \llbracket T \rrbracket)$, what is the cell $a$ is matched to in $\Mlexgr$ or is $a$ is left unmatched?}} \label{Question: where is a single vertex matched}    
\end{equation}

In addition to existence, the properties of a single edge $f$ in $G(\Psi \llbracket T \rrbracket)$ are also not hard to investigate; Lemma \ref{Lemma: local classification of delooped matrix elements} tells us whether $f$ is an isomorphism and $i(f),u(f)$ are easy to compute. Combining these one can straightforwardly decide whether in $G(\Psi \llbracket T \rrbracket )$ there exists an isomorphism edge $f$ to, or from, $a$ with specific values $i(f)$ and $u(f)$. This is to say that given a cell $a$ and integers $j,k$ one can evaluate the function:
    $$
    \operatorname{isopair}(a,j,k)=
    \begin{cases}
    b, & \text{if } \exists f\colon a\to b, \ f \text{ isomorphism edge of }G(\Psi \llbracket T \rrbracket ), \  i(f)=j,\ u(f)=k\\
    b, & \text{if }\exists g\colon b\to a, \ g \text{ isomorphism edge of }G(\Psi \llbracket T \rrbracket ), \ i(g)=j,\ u(g)=k\\
    \star, & \text{otherwise}
    \end{cases}
    $$
Using this, we can answer Question \ref{Question: where is a single vertex matched} for $\Mlex$, where the algorithm is simple, and for $\Mgr$, where a possibly branching recursion is needed.
\begin{algorithm*}
\caption{Matching a vertex in $\Mlex \Psi\llbracket T \rrbracket$ and answering Question \ref{Question: where is a single vertex matched}, case $\Mlex$.} \label{Algorithm: matching a single cell in Mlex}
 \textbf{Input} $a\in \operatorname{cells}(\Psi\llbracket T \rrbracket)$, tangle diagram $T$ is implicit  \\
  \textbf{Output}  $a$ is unmatched or $a$ is matched to $b\in \operatorname{cells} \llbracket T \rrbracket$
    in $\Mlex$  
   \begin{algorithmic}[1]
\For{k=1,\dots, \#\text{Morse layers of $T$}}
    \State $b\gets \operatorname{isopair}(a,k,k)$
    \If{$b\neq \star$}
        \State \Return $a$ matched to $b$ 
    \EndIf
\EndFor
\State \Return $a$ unmatched
\end{algorithmic}
\end{algorithm*}

\begin{algorithm*}
\caption{Matching a vertex in $\Mgr \Psi\llbracket T \rrbracket$ and answering Question \ref{Question: where is a single vertex matched}, case $\Mgr$.} \label{Algorithm: matching a single cell in Mgr}
 \textbf{Input} $a\in \operatorname{cells}\Psi\llbracket T \rrbracket$, tangle diagram $T$ is implicit  \\
  \textbf{Output}  $a$ is unmatched or $a$ is matched to $b\in \operatorname{cells} \llbracket T \rrbracket$
    in $\Mgr$  
   \begin{algorithmic}[1]
\Function{MatchesBack}{$x,y$} 
  \If{$x=\star$}
    \State \Return False
  \EndIf
\For{k=1,\dots, \#\text{Morse layers of $T$}}
    \For{j=k,\dots,1}    
        \State $z\gets \operatorname{isopair}(x,j,k)$  \label{algoline: isopair 1}
            \If{z=y}
                \State \textbf{return} True
            \ElsIf{\Call{MatchesBack}{$z,x$}} \label{algoline: MatchesBack 1}
                \State \Return False
            \EndIf
    \EndFor
\EndFor
\EndFunction

\State
   
\For{k=1,\dots, \#\text{Morse layers of $T$}}
    \For{j=k,\dots,1}    
        \State $b\gets \operatorname{isopair}(a,j,k)$ \label{algoline: isopair 2}
        \If{\Call{MatchesBack}{$b,a$}} \label{algoline: MatchesBack 2}
            \State \Return $a$ matched to $b$ 
        \EndIf
    \EndFor
\EndFor
\State \Return $a$ unmatched
\end{algorithmic}
\end{algorithm*}

\begin{lemma} \label{Lemma: equivalence of matching based on first characters}
    Let $x,y$ be cells of $\Psi \llbracket T_1 \rrbracket$ and $\Psi \llbracket T_2 \rrbracket$ respectively. Assume that first $m$ Morse layers of $T_1$ and $T_2$ agree and the first $m$ characters on the symbolic encodings of $x$ and $y$ also coincide. The following equivalences hold: 
    \begin{alignat}{7}
    \exists a_1&\colon  \ (a_1\to x)&&\in \Mlex, \ u(a_1\to x)&&\leq m \quad &&\iff \quad 
    \exists a_2&&\colon  \ (a_2\to y)&&\in \Mlex, \ u(a_2\to y)&&\leq m \\
    \exists b_1&\colon  \ (x\to b_1)&&\in \Mlex, \ u(x\to b_1)&&\leq m \quad &&\iff \quad 
    \exists b_2&&\colon  \ (x\to b_2)&&\in \Mlex, \ u(x\to b_2)&&\leq m \\
    \exists c_1&\colon  \ (c_1\to x)&&\in \Mgr, \ u(c_1\to x)&&\leq m \quad &&\iff \quad 
    \exists c_2&&\colon  \ (c_2\to y)&&\in \Mgr, \ u(c_2\to y)&&\leq m \\
    \exists d_1&\colon  \ (x\to d_1)&&\in \Mgr, \ u(x\to d_1)&&\leq m \quad &&\iff \quad 
    \exists d_2&&\colon  \ (x\to d_2)&&\in \Mgr, \ u(x\to d_2)&&\leq m.
    \end{alignat}
\end{lemma}
\begin{proof}
    We first claim the following auxiliary result: \textit{For cells $w.r_1$ of $\Psi \llbracket T_1 \rrbracket$, $w.r_2$ of $\Psi \llbracket T_2 \rrbracket$ and $j,k\leq |w|$, we have $\operatorname{isopair}((w.r_1),j,k)\in \{\star\} \cup \{v.r_1 \mid v \in \mathcal{S}\}$ and}
    \begin{alignat*}{2}
     &\operatorname{isopair}((w.r_1),j,k)=\star \quad &&\implies \quad\operatorname{isopair}((w.r_2),j,k)=\star   \\
    &\operatorname{isopair}((w.r_1),j,k)=v.r_1, \quad |v|=|w| \quad &&\implies \quad \operatorname{isopair}((w.r_2),j,k)=v.r_2.
    \end{alignat*}
    The auxiliary result holds due to the fact that a loop in the vertex $w.r_1$, which realizes the value of $u$, has to be completely contained in Morse layers $1,\dots ,|w|$ and remains unaffected when swapping $r_1$ with $r_2$. 

    In the lemma, the existence of $a_1,a_2,b_1$ and $b_2$ can be obtained by running Algorithm \ref{Algorithm: matching a single cell in Mlex} for $(x,T_1)$ and $(y,T_2)$ and stopping when $k=m+1$. Similarly, the existence of $c_1,c_2,d_1$ and $d_2$ is attained from Algorithm \ref{Algorithm: matching a single cell in Mgr}. The statement of the lemma follows now from analyzing the algorithms separately with the auxiliary result: for both inputs $x$ and $y$ the same number of orders is executed and line-by-line their action can be controlled with the auxiliary result.
\end{proof}

\begin{lemma}\label{Lemma: u/i}
   Let $a\to b\to c$ be a path in $G(\Psi \llbracket T \rrbracket, \Mlexgr )$ with $(b \to a) \in \Mlexgr$ and $(c\to d) \notin \Mlexgr$ for all $d$. Then $i(b\to c)\leq u(a\to b)$. Dually, let $e\to f \to g$ be a path in $G(\Psi \llbracket T \rrbracket, \Mlexgr )$ with $(g \to f) \in \Mlexgr$ and $(h\to e) \notin \Mlexgr$ for all $h$. Then $i(e\to f)\leq u(f\to g)$.
\end{lemma}
Notice that if $a\to b \to c$ or $e\to f \to g$ are subpaths of either a cycle or a path contributing towards the differential $\partial$ of $\Mlexgr \Psi \llbracket T \rrbracket$, then the conditions $(c\to d), (h\to e) \notin \Mlexgr$ automatically hold. This follows from Condition \ref{Condition: matching condition of Morse} of Morse matching which $\Mlexgr$ satisfy by construction.

\begin{proof}
    If $i(b\to c)>u(a\to b)$, then $b$ and $c$ agree on the first $u(a\to b)$ characters. Hence by Lemma \ref{Lemma: equivalence of matching based on first characters} there has to exist $(c\to d)\in \Mlexgr$ contradicting our assumptions. The case $e\to f \to g$ is analogous.
\end{proof}

\begin{proposition}\label{Proposition: lexicographic matching is Morse}
    For any tangle diagram $T$, the lexicographic matching $\Mlex$ is a Morse matching on $\Psi \llbracket T \rrbracket$.
\end{proposition}


\begin{proof}
    Assume towards contradiction that there exists a directed cycle $L$ in $G(\Psi \llbracket T \rrbracket, \Mlex)$. Let $a\to b\to c \to  d$ be a subpath of $L$ such that $(c\to b)\in \Mlex$ and $i(b\to c)\leq i(f)$ for all edges $f$ of $L$. Denote the $i(b\to c)$:th Morse layer of $T$ with $F$.

    \textbf{Case I:} $F=\vert \cdots \vert \ \begin{tikzpicture}[baseline={(0,0.05)},scale=0.3]
    \negCrossing{0,0} 
\end{tikzpicture} \ \vert \cdots \vert $. Since $(c\to b) \in \Mlex$ we have $i(b\to c)=u(b\to c)$. The $1$-smoothing of $F$ in $b$ must be $\color{green}\vert \color{black}\cdots  \color{green} \vert \ \begin{tikzpicture}[baseline={(0,0.05)},scale=0.3, color=green]
    \lineUpBendLeftMore{1,0}
    \lineUpBendRightMore{0,0}
\end{tikzpicture} \ \color{green}\vert \color{black}\cdots \color{green} \vert $ 
 and $(c\to b) \in G(\Psi \llbracket T \rrbracket)$ is merge-type, see Figure \ref{Figure: merge type and split type}. Lemma \ref{Lemma: u/i} and the $i$-minimality of $b\to c$ yield
$$
i(b\to c)\leq i(c\to d) \leq u(b\to c) =i(b\to c)
$$
giving $i(c\to d)=i(b\to c)$. Now each possible type i), ii) or v) of $c\to b$, yields $b=d$ which implies a contradiction: $(b\to c)\notin \Mlex$ and $(b\to c)\in \Mlex$.

\textbf{Case II:} $F=\vert \cdots \vert \ \begin{tikzpicture}[baseline={(0,0.05)},scale=0.3]
    \posCrossing{0,0} 
\end{tikzpicture} \ \vert \cdots \vert $. Now $0$-smoothing of $F$ in $c$ must be $\color{green}\vert \color{black}\cdots  \color{green} \vert \ \begin{tikzpicture}[baseline={(0,0.05)},scale=0.3, color=green]
    \lineUpBendLeftMore{1,0}
    \lineUpBendRightMore{0,0}
\end{tikzpicture} \ \color{green}\vert \color{black}\cdots \color{green} \vert $ and  hence $(c\to b) \in G(\Psi \llbracket T \rrbracket)$ is split-type. As in Case I, one can here that $a=c$ yielding a contradiction.

\textbf{Case III:} $F=\vert \cdots \vert \ \begin{tikzpicture}[baseline={(0,0.05)},scale=0.3]
    \lineRightBendUpMore{0,0} 
\end{tikzpicture} \ \vert \cdots \vert $ or $F=\vert \cdots \vert \ \begin{tikzpicture}[baseline={(0,0.05)},scale=0.3]
    \lineRightBendDownMore{0,1}
\end{tikzpicture} \ \vert \cdots \vert $. This is impossible; the function $i$ can only take values in the indices of the Morse layers which contain crossings in them.
\end{proof}

A Morse presentation of a tangle diagram is called \term{a negative braid word}, if its every Morse layer is of the form $\vert \cdots \vert \ \begin{tikzpicture}[baseline={(0,0.05)},scale=0.3]
    \negCrossing{0,0} 
\end{tikzpicture} \ \vert \cdots \vert $ and the strands are oriented from bottom to top. The canonical generators $\sigma_i$ of the braid groups correspond to the crossings of the tangle diagram in a standard way, so that
$$
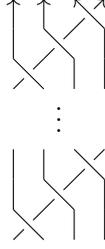
\begin{tikzpicture}
    \begin{scope}[scale=0.4]
        \negCrossing{0,0} \lineUp{2,0} \lineUp{3,0}
        \lineUp{0,1} \negCrossing{1,1} \lineUp{3,1}
        \lineUp{0,2} \lineUp{1,2} \negCrossing{2,2}

        \negCrossing{0,5} \lineUp{2,5} \lineUp{3,5}
        \lineUp{0,6} \negCrossing{1,6} \lineUp{3,6}
        \lineUpWArrow{0,7} \lineUpWArrow{1,7} \negCrossingWArrows{2,7}
    \end{scope}
    \node at (0.6,1.7) {$\vdots$};
\end{tikzpicture}
$$
represents the negative braid word of the 4-strand torus braid $(\sigma_1 \sigma_2 \sigma_3)^n$.

\begin{lemma}\label{lemma: red zeros matched}
    Let $B$ be a negative braid word and let $a\to b$ be a merge-type isomorphism in the graph $G(\Psi \llbracket B \rrbracket , \Mgr)$.
    \begin{enumerate}
        \item If $i(a\to b)=u(a\to b)$, then $(a\to b)\in \Mgr$ or $a,b$ are both otherwise matched in $\Mgr$. \label{Lemma subitem: negbraid matching 1}
        \item If $i(a\to b)< u(a\to b)$, then $(a\to b) \notin \Mgr$ since $a$ is otherwise matched in $\Mgr$.\label{Lemma subitem: negbraid matching 2}
        \item Unmatched cells of $G(\Psi \llbracket B \rrbracket, \Mgr )$ contain only characters $0,\bzero,1$.
    \end{enumerate}
\end{lemma}
\begin{proof}
    Suppose $i(a\to b)=u(a\to b)$. Then $a$ and $b$ agree on the first $i(a\to b) -1$ characters and thus by Lemma \ref{Lemma: equivalence of matching based on first characters} either both $a$ and $b$ or neither of them is matched with arrows whose $u$ value is at most $i(a\to b)-1$. In case neither is, it follows from Algorithm \ref{Algorithm: matching a single cell in Mgr} that $(a\to b)\in \Mgr$. 

    Now suppose $i(a\to b) < u(a\to b)$. Consider the highest point of the red circle in $a$ from which $u$ is measured. This Morse layer must be $0$-smoothed as $\color{green}\vert \color{black}\cdots  \color{green} \vert \ \begin{tikzpicture}[baseline={(0,0.05)},scale=0.3, color=green]
    \begin{scope}[green]
        \lineRightBendDownMore{0,1}
    \end{scope}
    \begin{scope}[red]
        \lineRightBendUpMore{0,0}
    \end{scope}
\end{tikzpicture} \ \color{green}\vert \color{black}\cdots \color{green} \vert $. By the paragraph above, $a$ is matched to a cell whose same Morse layer is $\color{green}\vert \color{black}\cdots  \color{green} \vert \ \begin{tikzpicture}[baseline={(0,0.05)},scale=0.3, color=green]
    \lineUpBendLeftMore{1,0}
    \lineUpBendRightMore{0,0}
\end{tikzpicture} \ \color{green}\vert \color{black}\cdots \color{green} \vert $ or $a$ is matched to some other vertex which is not $b$. 

Let $c$ be an unmatched cell of $\Psi \llbracket B \rrbracket \subset \{0,\bzero, \rzero,1,\bone,\rone, \bullet, \bbullet,\rbullet \}^n$. Since $B$ is a braid, all of its Morse layers are crossings and thus the characters $\bullet,\rbullet,\bbullet$ cannot appear in $c$. The highest points of loops are contained in $0$-smoothed Morse layers  $\color{green}\vert \color{black}\cdots  \color{green} \vert \ \begin{tikzpicture}[baseline={(0,0.05)},scale=0.3, color=green]
    \lineRightBendDownMore{0,1}
    \lineRightBendUpMore{0,0}
\end{tikzpicture} \ \color{green}\vert \color{black}\cdots \color{green} \vert $ and hence characters $\rone, \bone$  are excluded. By Claim \ref{Lemma subitem: negbraid matching 1}, every cell containing a red loop is matched and thus $\rzero$ cannot be contained in $c$. From what remains, we conclude that $c\in \{0,\bzero,1\}^n$.
\end{proof}

\begin{proposition}\label{Proposition: Mgr Morse on T4}
    For $n\in \mathbb Z_{\geq 0}$, the greedy matching $\Mgr$ is a Morse matching on $\Psi\llbracket(\sigma_1\sigma_2\sigma_3)^n\rrbracket$.
\end{proposition}
\begin{proof}
    Again, we assume towards contradiction that there exists a directed cycle $L$ in the graph $G(\Psi \llbracket (\sigma_1\sigma_2\sigma_3)^n \rrbracket, \Mgr)$. Let $z\to a\to b\to c \to  d$ be a subpath of $L$ such that $(a\to z),(c\to b)\in \Mgr$ and $i(b\to c)\leq i(f)$ for all edges $f$ of $L$. We denote $I=i(b\to c)$.

    If $c\to b$ is a merge-type isomorphism, it follows from Lemma \ref{lemma: red zeros matched} that $I=u(b\to c)$. From this one can reach a contradiction $b=d$ similar to Case I of the proof of Proposition \ref{Proposition: lexicographic matching is Morse}. We can thus assume that $c \to b $ is split-type isomorphism and divide into cases with $I$.

    \textbf{Case I:} $I\in 3\mathbb Z +1$. The $I$:th Morse layer  is $1$-smoothed in $b$ and $0$-smoothed in $c$. In the Morse layers $1,\dots, I-1$ of $b$ and $c$ the $\ast$-marked component is connected to either to $A$, to $B$ or in the bottom boundary to $C$: 
        \input{Acyclicity_proof/I1mod3fig1}
    The only way $b\to c$ is a split-type isomorphism is if $\ast$ is connected to $B$. In $(\sigma_1\sigma_2\sigma_3)^n$ the $I+1$:th Morse layer is $\begin{tikzpicture}[baseline={(0,0.05)},scale=0.3]
    \lineUp{0,0}
    \negCrossing{1,0}
    \lineUp{3,0}
\end{tikzpicture}$ and by Lemma \ref{lemma: red zeros matched} it has to be smoothed in $b$ as 
 \begin{tikzpicture}[baseline={(0,0.05)},scale=0.3, color=green]
    \begin{scope}[green]
        \lineUp{-1,0}
        \lineRightBendDownMore{0,1}
        \lineUp{2,0}
    \end{scope}
    \begin{scope}[blue]
        \lineRightBendUpMore{0,0}
    \end{scope}
\end{tikzpicture}
    and so $b$ and $c$ are of the form:
    \begin{center}
    \input{Acyclicity_proof/I1mod3fig2}    
    \end{center}
    Hence $u(b\to c)=I+1$ and by $i$ mimimality if $b\to c$ and Lemma \ref{Lemma: u/i} we get 
    $$
    I \leq i(a\to b) \leq u(b\to c) =I +1.
    $$
    The $I+1$:th Morse layer of $b$ is $0$-smoothed rendering $i(a\to b)=I+1$ impossible. Thus $i(a\to b)=I$ but since $b\to c$ is split type, we get that $a=c$ which is a contradiction.
    
    \textbf{Case II:} $I\in 3\mathbb Z +2$. We similarly divide into cases based on connectivity of $\ast$ in the Morse layers $1,\dots , I-1$ of $b$ and $c$:
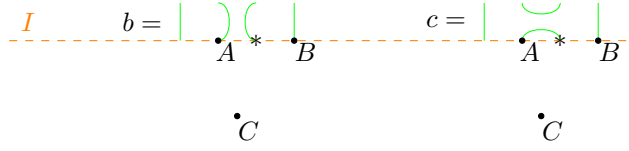
    \begin{center}
        \input{Acyclicity_proof/I2mod3fig1}
    \end{center}
    The connection $\ast$ to $B$ is almost identical with the similar connectivity in the case $I\in 3\mathbb Z +1$. The only difference in arguments can be handled by moving the left-most line to the right-most position in all of the pictures. Since $c\to b$ is split-type, $\ast$ cannot be connected to $A$ and we can proceed by assuming that $\ast$ is connected to the bottom boundary $C$. In order for $b$ to be matched to $c$, they have to be of the form:
    \begin{center}
        \input{Acyclicity_proof/I2mod3fig2}
    \end{center}
    Hence $u(b\to c)=I+2$ and it follows that either $i(a\to b)\in \{ I,I+1,I+2\}$. From these, $i(a\to b)=I+2$ is impossible, since $I+2$:th Morse layer of $b$ is $0$-smoothed. The value $i(a\to b)=I$ yields a contradiction with $a=c$, so we proceed with $i(a\to b)=I+1$.

    We denote $x=\operatorname{isopair(a,I,I+2)}$, so the cells at our disposal are pictorially:
    \begin{center}
        \input{Acyclicity_proof/I2mod3fig3}
    \end{center}
    Since $(c \to b) \in \Mgr$, we aim to show that $(x\to a) \in \Mgr$ as well. This will be a contradiction, as $(a\to z)\in \Mgr$ and $\Mgr$ is a matching. Now $(x\to a) \in A_{I,I+2}$, so in order for $(x\to a) \notin M_{I,I+2}$ either $a$ or $x$ has to be matched in $M_{I+1,I+2}$, $M_{I+2,I+2}$ or in $M_{J,K}$ for some $J,K\leq I+1$. We can exclude these possibilities one-by-one. 
    \begin{itemize}
        \item By staring at the diagrams of $a$ and $x$, one can conclude that $\operatorname{isopair}(y,j,k)= \star$, whenever
    \begin{multline*}
    (y,j,k)\in \{(a,I+1,I+2),(a,I+2,I+2),(a,I,I+1),(a,I+1,I+1), \\(a,I,I),(x,I+1,I+2),  (x,I+2,J+2)\} \cup \{(x,J,I+1) \mid J\leq I+1 \}.        
    \end{multline*}
    Hence $a$ and $x$ cannot be matched in the corresponding $M_{j,k}$ sets.
        \item Since $(c\to b)\in \Mgr$ and $u(c\to b)=I+2$, the cell $c$ cannot be matched in $M_{J,K}$ for $K\leq I$. The first $I$ characters of $c$ and $x$ coincide, so by Lemma \ref{Lemma: equivalence of matching based on first characters}, $x$ is not matched in any $M_{J,K}$ for $J,K\leq I$ either.
        \item We know that $i(a\to z)\geq I$ and $(a\to z)\in \Mgr$ which is why $a$ cannot be matched in any $M_{J,K}$ with $J<I$.  
    \end{itemize}
    The three arguments above eliminate all opportunities for $a$ and $x$ to be matched before $M_{I,I+2}$. Thus $(x\to a)\in M_{I,I+2}\subset \Mgr$ and we have reached a contradiction.

    \textbf{Case III:} $I\in 3\mathbb Z$. In order for $c\to b$ to be a split-type isomorphism, the connectivity below $I$:th Morse layer and the smoothing at $I+1$:th Morse layer have to be look like
    \begin{center}
        \input{Acyclicity_proof/I3mod3fig1}
    \end{center}
    Now it can be argued with Lemma \ref{lemma: red zeros matched} that $I+2$:th Morse layer of $b$ has to be  \begin{tikzpicture}[baseline={(0,0.05)},scale=0.3, color=green]
    \begin{scope}[green]
        \lineUp{-1,0}
        \lineRightBendDownMore{0,1}
        \lineUp{2,0}
    \end{scope}
    \begin{scope}[blue]
        \lineRightBendUpMore{0,0}
    \end{scope}
\end{tikzpicture}. Hence $u(b\to c)=I+2$ and with the help of \ref{Lemma: u/i} we can conclude that $i(a\to b)=I+1$ so that diagrammatically we acquire:
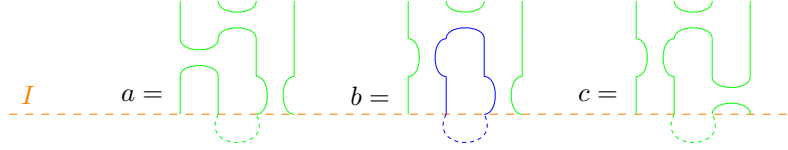
\begin{center}
    \input{Acyclicity_proof/I3mod3fig2}
\end{center}
Since $i(c\to b)\in \Mgr$ instead of $(a\to c)\in \Mgr$ the arrow $a\to z$ was prioritized in $\Mgr$ over the arrow $a\to c$, that is,
$$
(a\to z) \in M_{I+2,I+2}\cup M_{I,I+1} \cup M_{I+1,I+1}\cup M_{I,I} \cup \bigcup_{J<I} M_{J,K}.
$$
The first four cases can be deemed impossible as $(a\to z) \notin A_{j,k}$ for the relevant $j,k$ whereas last case is eliminated by $i(a\to z)\geq I$. Again, we have reached a contradiction and thus proven that $\Mgr$ is a Morse matching on $\Psi\llbracket(\sigma_1\sigma_2\sigma_3)^n\rrbracket$.
\end{proof}

\section{Torus braids with 4 strands}\label{Section: Main results}
In this section, we derive the main results from properties of the Morse complexes $\Mgr \Psi (\llbracket (\sigma_1\sigma_2 \sigma_3)^n \rrbracket)$ whose proofs are postponed to Sections \ref{Section: proof of vertex correspondences} and  \ref{Section: proof of morphism correspondences}.

\subsection{Upper bounds on the homological dimensions}
For brevity, throughout Sections \ref{Section: Main results}, \ref{Section: proof of vertex correspondences} and \ref{Section: proof of morphism correspondences} we will write $C_n= \Mgr \Psi (\llbracket (\sigma_1\sigma_2 \sigma_3)^n \rrbracket)$ and $U_n$ for the unmatched cells of $G(\Psi \llbracket (\sigma_1\sigma_2 \sigma_3)^n\rrbracket, \Mgr )$, so that $U_n$ forms a basis for $C_n$. For a cell $w$ we denote $h(w)$ and $q(w)$ for its the homological and quantum degrees respectively. Let $X$ denote the set of all matched and unmatched cells, $X=\bigcup_{n=0}^\infty \operatorname{cells}(\Psi \llbracket (\sigma_1\sigma_2 \sigma_3)^n\rrbracket)$ and define functions
\begin{align*}
    t_{\bA} &\colon  X \to \mathbb R, \quad t_{\bA}(w)=-\frac{1}{2}n(w) +\frac{1}{2}h(w) -\frac{1}{4}q(w) -1 \\    
    t_{\bB} &\colon X \to \mathbb R, \quad t_{\bB}(w)=3n(w) -\frac{3}{2}h(w) +q(w) -\frac{5}{2} \\
    t_{\bC} &\colon X \to \mathbb R, \quad t_{\bC}(w)=-\frac{9}{4}n(w) +h(w) -\frac{3}{4}q(w) -\frac{1}{4}. 
\end{align*}
where  $n(w)$ denotes the unique index $n$ for which $w\in \operatorname{cells}(\Psi \llbracket (\sigma_1\sigma_2 \sigma_3)^n\rrbracket)$. In order to compare the complexes $C_n$ for different $n$, grading shifts are required: For an individual cell $w$ and or a whole complex $A$ we write $w[a]\{b\}$ or $A[a]\{b\}$ to shift the homological degree by $a$ and the quantum degree by $b$.

\begin{proposition}\label{Proposition: abc bijections on cells}
    For all $n\geq 0$ there are bijections
    \begin{align*}
        \mathfrak{a}_n\colon& \{ w\in U_n \mid t_{\bA}(w)\geq 1 \} \to \{ w\in U_{n+4} \mid t_{\bA}(w)\geq 2 \} \\ 
        \mathfrak{b}_n\colon& \{ w\in U_n \mid t_{\bB}(w)\geq 1 \} \to \{ w\in U_{n+4} \mid t_{\bB}(w)\geq 2 \} \\ 
        \mathfrak{c}_n\colon& \{ w\in U_n \mid t_{\bC}(w)\geq 1 \} \to \{ w\in U_{n+4} \mid t_{\bC}(w)\geq 2 \}.  
    \end{align*}
    The bijections shift the homological and quantum degrees and maintain the connectivity of the underlying Temperley-Lieb diagrams, so that
    $$
    w_1=(\mathfrak{a}_n(w_1) )\{12\}, \quad w_2= (\mathfrak{b}_n(w_2) )[6] \{20\} , \quad w_3=(\mathfrak{c}_n(w_3) )[8] \{24\} 
    $$
    as objects in the category $\operatorname{Mat}(\Cob(8))$.
  \end{proposition}

\begin{figure}
\centering
\begin{tikzpicture}[scale=0.5]
\begin{scope}[shift={(0,0)}]
    \oneSmoothing{0,0} 
    \lineUp{2,0}
    \lineUp{3,0}
    \lineUp{0,1}
    \oneSmoothing{1,1}
    \lineUp{3,1}
    \lineUp{0,2}
    \lineUp{1,2}
    \oneSmoothing{2,2}
    
    \node at (4,0.5) {1};
    \node at (4,1.5) {1};
    \node at (4,2.5) {1};
    \node at (4,3.5) {1};
    \node at (4,4.5) {1};
    \node at (4,5.5) {1};
    \node at (4,6.5) {1};
    \node at (4,7.5) {1};
    \node at (4,8.5) {1};
    \node at (4,9.5) {1};
    \node at (4,10.5) {1};
    \node at (4,11.5) {1};
\end{scope}

\begin{scope}[shift={(0,3)}]
    \oneSmoothing{0,0}
    \lineUp{2,0}
    \lineUp{3,0}
    \lineUp{0,1}
    \oneSmoothing{1,1}
    \lineUp{3,1}
    \lineUp{0,2}
    \lineUp{1,2}
    \oneSmoothing{2,2}
\end{scope}

\begin{scope}[shift={(0,6)}]
    \oneSmoothing{0,0}
    \lineUp{2,0}
    \lineUp{3,0}
    \lineUp{0,1}
    \oneSmoothing{1,1}
    \lineUp{3,1}
    \lineUp{0,2}
    \lineUp{1,2}
    \oneSmoothing{2,2}
\end{scope}

\begin{scope}[shift={(0,9)}]
    \oneSmoothing{0,0}
    \lineUp{2,0}
    \lineUp{3,0}
    \lineUp{0,1}
    \oneSmoothing{1,1}
    \lineUp{3,1}
    \lineUp{0,2}
    \lineUp{1,2}
    \oneSmoothing{2,2}
\end{scope}

\begin{scope}[shift={(7,0)}]
    \oneSmoothing{0,0}
    \lineUp{2,0}
    \lineUp{3,0}
    \lineUp{0,1}
    \zeroSmoothing{1,1}
    \lineUp{3,1}
    \lineUp{0,2}
    \lineUp{1,2}
    \oneSmoothing{2,2}
    \lineRightBendUpMore{0,3}
    \lineUp{2,3}
    \lineUp{3,3}
    \lineUpBendLeftMore{2,4}
    \lineUp{3,4}
    \begin{scope}[color=blue]
        \lineRightBendDownMore{0,4}
        \lineUp{0,4}
        \lineUpBendRightMore{1,4}
        \lineUp{0,5}
        \lineUp{1,5}
        \lineRightBendUpMore{0,6}
    \end{scope}
    \oneSmoothing{2,5}
    \lineUp{2,6}
    \lineUp{3,6}
    \lineRightBendDownMore{0,7}
    \lineUp{0,7}
    \zeroSmoothing{1,7}
    \lineUp{3,7}
    \lineUp{0,8}
    \lineUp{1,8}
    \lineRightBendUpMore{2,8}
    \oneSmoothing{0,9}
    \begin{scope}[color=blue]
        \lineRightBendDownMore{2,9}
        \lineUp{2,9}
        \lineUp{3,9}
        \lineUpBendLeftMore{2,10}
        \lineUp{3,10}
        \lineRightBendUpMore{2,11}
    \end{scope}
    \lineUp{0,10}
    \lineUpBendRightMore{1,10}
    \lineUp{0,11}
    \lineUp{1,11}
    \begin{scope}[green]
        \lineRightBendDownMore{2,12}
    \end{scope}
    
    \node at (4,0.5) {1};
    \node at (4,1.5) {0};
    \node at (4,2.5) {1};
    \node at (4,3.5) {0};
    \node at (4,4.5) {1};
    \node at (4,5.5) {1};
    \node at (4,6.5) {$\bzero$};
    \node at (4,7.5) {0};
    \node at (4,8.5) {0};
    \node at (4,9.5) {1};
    \node at (4,10.5) {1};
    \node at (4,11.5) {$\bzero$};
\end{scope}

\begin{scope}[shift={(14,0)}]
    \begin{scope}[color=blue]
    \zeroSmoothing{0,0}
    \lineUp{2,0}
    \lineUp{3,0}
    \lineUp{0,1}
    \oneSmoothing{1,1}
    \lineUp{3,1}
    \lineUp{0,2}
    \lineUp{1,2}
    \zeroSmoothing{2,2}
    \zeroSmoothing{0,3}
    \lineUp{2,3}
    \lineUp{3,3}
    \lineUp{0,4}
    \zeroSmoothing{1,4}
    \lineUp{3,4}
    \lineUp{0,5}
    \lineUp{1,5}
    \oneSmoothing{2,5}
    \zeroSmoothing{0,6}
    \lineUp{2,6}
    \lineUp{3,6}
    \lineUp{0,7}
    \oneSmoothing{1,7}
    \lineUp{3,7}
    \lineUp{0,8}
    \lineUp{1,8}
    \lineRightBendUpMore{2,8}
    \lineRightBendUpMore{0,9}
    \end{scope}

    \begin{scope}[color=green]
    \lineRightBendDownMore{2,9}
    \lineRightBendDownMore{0,10}
    \lineUp{2,9}
    \lineUp{3,9}
    \lineUp{0,10}
    \zeroSmoothing{1,10}
    \lineUp{3,10}
    \lineUp{0,11}
    \lineUp{1,11}
    \oneSmoothing{2,11}
    
    \end{scope}
    
    \node at (4,0.5) {0};
    \node at (4,1.5) {1};
    \node at (4,2.5) {$\bzero$};
    \node at (4,3.5) {$\bzero$};
    \node at (4,4.5) {0};
    \node at (4,5.5) {1};
    \node at (4,6.5) {0};
    \node at (4,7.5) {1};
    \node at (4,8.5) {$\bzero$};
    \node at (4,9.5) {$\bzero$};
    \node at (4,10.5) {0};
    \node at (4,11.5) {1};

\end{scope}

\end{tikzpicture}
    \caption{The three periodic building blocks $\bA, \bB$ and $\bC $ of unmatched cells of $G(\Psi \llbracket (\sigma_1\sigma_2 \sigma_3)^n\rrbracket, \Mgr )$. The maps $\mathfrak a_n$, $\mathfrak b_n$ and  $\mathfrak c_n$ from Proposition \ref{Proposition: abc bijections on cells} take a diagram $w$ of $U_n$ which contains at least one of the periodic building blocks $\bA, \bB$ and $\bC $ and duplicates the first instance of this block, making the tangle diagram 12 Morse layers taller and an unmatched cell in $U_{n+4}$. 
    }
    \label{Figure: periodic building blocks}
\end{figure}

    Diagrammatically the correspondences $\mathfrak a_n$, $\mathfrak b_n$ and $\mathfrak c_n$ are described in Figure \ref{Figure: periodic building blocks} and a more formal description will be given in Section  \ref{Section: proof of vertex correspondences}.

    The $i$th chain space of a complex $A$ over $\operatorname{Mat}(\Cob(2m))$ is denoted by $A^i$. To further split $A^i$ to quantum gradings, denote the direct sum of those Temperley-Lieb diagrams with $\{j\}$ quantum grading shift by $A^{i,j}$, so that
    $
    \bigoplus_{j\in \mathbb Z} A^{i,j} = A^i.
    $
    To measure the size an object $B\in\operatorname{Cob}(2m)$, we define 
    $\dimcob(B)\in \mathbb Z_{\geq 0}$ to be the number of loopless Temperley-Lieb diagrams in $\Psi(B)$. Similarly we define size of the whole chain complex with 
    $$
    \dimkom (A)=\sum_{i\in \mathbb Z} \dimcob(A^i). 
    $$  
    
  \begin{corollary}\label{Corollary: bounds on Cn dimensions}
      The dimensions of $C_n^{i,j}$ and $C_n$ admit constant and quadratic bounds respectively:
      \begin{align}
          \max \{ &\dimcob(C_n^{i,j}) \mid i,j\in \mathbb Z,\ n\geq 0\} < \infty \\
          &\dimkom(C_n) =O(n^2).
      \end{align}
  \end{corollary}
  \begin{proof}
    With \texttt{Lean}'s linarith tactic, one can show that over $\mathbb Q$ the system of inequalities 
    \begin{align*}
        2&>-\frac{1}{2}n +\frac{1}{2}h -\frac{1}{4}q -1 \\    
    2&>3n -\frac{3}{2}h +q -\frac{5}{2} \\
    2&>-\frac{9}{4}n +h -\frac{3}{4}q -\frac{1}{4} 
    \end{align*}
    implies $n<39$. Hence,  either  $t_{\bA}(w) \geq 2$, or $t_{\bB}(w) \geq 2$ or $t_{\bC}(w) \geq 2$ holds for any $w\in U$ with $n(w)\geq 39$. It follows from Proposition \ref{Proposition: abc bijections on cells} that 
      $$
      \max \{ \dimcob(C_n^{i,j}) \mid i,j\in \mathbb Z,\ 0\leq n\} = 
      \max \{ \dimcob(C_n^{i,j}) \mid i,j\in \mathbb Z,\  0\leq n\leq 38\} < \infty.
      $$
      Morse complex $C_n$ consists of unmatched cells of the complex $\Psi \llbracket (\sigma_1 \sigma_2 \sigma_3)^n \rrbracket$ which yields
      $$
      \dimcob (C_n^{i,j}) \leq \dimcob ((\Psi \llbracket (\sigma_1 \sigma_2 \sigma_3)^n \rrbracket)^{i,j}).
      $$
      On the other hand the complex $\Psi \llbracket (\sigma_1 \sigma_2 \sigma_3)^n \rrbracket$ is supported in bigradings $[-3n,0] \times [ -9n,-3n]$, hence
      $$
          \dimkom(C_n) \leq 18n^2 \cdot \max \{ \dimcob(C_n^{i,j}) \mid i,j\in \mathbb Z,\ n\geq 0\} =O(n^2).
      $$
  \end{proof}

\begin{figure}
    \centering
    \input{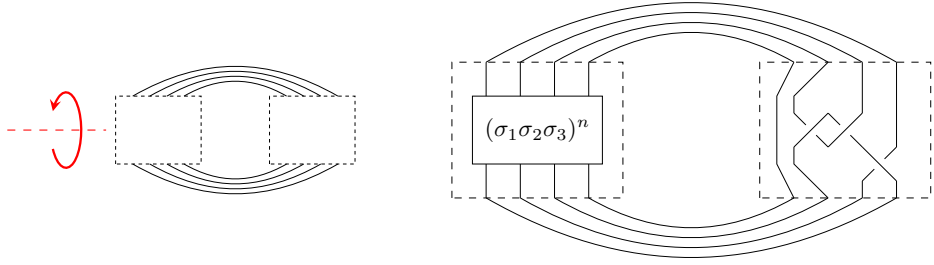}
    \caption{Our convention for composing 8-ended tangles (left) and an example link $(\sigma_1\sigma_2\sigma_3)^n \otimes_{\operatorname{Link}} T$ (right), for which Theorem \ref{Theorem: homology bound for compositions of T4} applies. The red rotation axis is applied in the proof of Theorem \ref{Theorem: homology bound for compositions of T4}. }
    \label{Figure: composition of tangles into unoriented link}
\end{figure}

Given two 8-ended tangles $T_1$ and $T_2$ we define an unoriented link  $T_1\otimes_{\operatorname{Link}} T_2$ by connecting the strands together as in Figure \ref{Figure: composition of tangles into unoriented link}. Choosing an orientation $o$ on $T_1\otimes_{\operatorname{Link}} T_2$ yields an oriented link which is denoted with  $T_1\otimes_{\operatorname{Link},o} T_2$. From the wonderful work of Bar-Natan, we know that not only tangles can be composed with this planar algebra operation, but also two objects $O_1$, $O_2$ of $\operatorname{Cob}(8)$ can be fused as an object $O_1\otimes_{\operatorname{Cob}} O_2$ in $\operatorname{Cob}(0)$ (join the pictures and add up their grading shifts).  Even complexes $A,B\in \operatorname{Mat}(\Cob(8))$ can be merged as a complex $A\otimes_{\operatorname{Kom}}B $ where
$$
(A\otimes_{\operatorname{Kom}}B)^{i,j}= \bigoplus_{a+c=i,\ b+d=j} A^{a,b}\otimes_{\operatorname{Cob}}B^{c,d}.
$$
Assuming that orientations agree throughout, the Bar-Natan bracket complex commutes with the composition of tangles: 
$$
\llbracket T_1\otimes_{\operatorname{Link},o} T_2 \rrbracket \simeq \llbracket T_1 \rrbracket \otimes_{\operatorname{Kom}} \llbracket T_2 \rrbracket.
$$

The standard Khovanov TQFT is an additive functor $F$ from $\Cob(0)$ to graded $\mathbb Z$-modules which maps a diagram with $k$ loops to $(\mathbb Z \{-1 \} \oplus \mathbb Z\{ 1\})^{\otimes k}$. Naturally, $F$ can be applied to a whole chain complex $A$ over $ \operatorname{Mat}(\Cob(0))$ to obtain a bigraded complex $F(A)$ and to access the $(i,j)$-graded $\mathbb Z$-module of $F(A)$, we write $F^{i,j}(A)$. For a bigraded complex of $\mathbb Z$-modules $B$, we write $H(B)$ for the bigraded homology module and $H^{i,j}(B)$ for the individual homology groups. Finally, for an oriented link diagram $L$ the even Khovanov homology $\operatorname{Kh}^{i,j}(L)=H^{i,j}(F\llbracket L \rrbracket)$ is the (co)homology obtained from
\begin{equation}\label{Equation: definition of Khovanov homology}
F^{i-1,j} \llbracket L \rrbracket  \to 
F^{i,j} \llbracket L \rrbracket  \to
F^{i+1,j} \llbracket L \rrbracket . 
\end{equation}

One can add a basepoint to a link diagram $L$, and from there to the diagrams of $\Cob(0)$ that generate the complex $\llbracket L \rrbracket$. With the basepoints in hand, $F$ gets replaced by the reduced version $\overline{F}$  where a diagram $D$ of $k$ loops and one basepoint gets mapped $(\mathbb Z \{-1 \} \oplus \mathbb Z\{ 1\})^{\otimes k-1}$. 
This yields the reduced Khovanov homology $\overline{\operatorname{Kh}} (L)=H(\overline{F} \llbracket L \rrbracket)$. The even Khovanov homology also has a strange cousin, odd Khovanov homology $\operatorname{Kh}_{\text{odd}}^{i,j}(L)$, which is also obtained from (\ref{Equation: definition of Khovanov homology}) although with different sign convention on the maps \cite{Ozsvth2013}. The odd theory also has a reduced version, $\overline{\operatorname{Kh}}_{\text{odd}}$, which satisfies 
\begin{equation}
\operatorname{Kh}_{\text{odd}}^{i,j}(L) \cong  \overline{\operatorname{Kh}}_{\text{odd}}^{i,j-1}(L) \oplus \overline{\operatorname{Kh}}_{\text{odd}}^{i,j+1}(L). \label{Equation: reduced unreduced Khovanov splitting}
\end{equation}
The corresponding relation does not hold for the even theory. The unreduced even and odd Khovanov homologies with coefficients $\operatorname{Kh}_{\text{odd/even}}^{i,j}(L; \mathbb F)$ are taken from (\ref{Equation: definition of Khovanov homology}) by first tensoring with $\mathbb F$ before taking the homology. 

\begin{theorem}\label{Theorem: homology bound for compositions of T4}
    Let $T$ be any $8$-ended tangle. The dimensions of even and odd Khovanov homologies of links $(\sigma_1\sigma_2\sigma_3)^n \otimes_{\operatorname{Link},o} T$ admit constant and quadratic bounds:
    \begin{align}
        \max_{i,j,n,o,\mathbb F, \mathcal H} \left\{ \operatorname{dim}_{\mathbb F} \left( \mathcal H ^{i,j}((\sigma_1\sigma_2\sigma_3)^n \otimes_{\operatorname{Link},o} T; \mathbb F) \right) \right \} &<\infty \label{Equation: tangle homology dim bound 1}\\
        \max_{o,\mathbb F, \mathcal{H}} \Big\{ \operatorname{dim}_{\mathbb F} \Big( \bigoplus_{i,j\in \mathbb Z}\mathcal H^{i,j} ((\sigma_1\sigma_2\sigma_3)^n \otimes_{\operatorname{Link},o} T; \mathbb F) \Big) \Big \} &=O(n^2). \label{Equation: tangle homology dim bound 2}
    \end{align}
    In (\ref{Equation: tangle homology dim bound 1}) the maximum is taken over all integers $i,j,n$, all orientations $o$, all fields $\mathbb F$ and all theories: $\mathcal{H}\in \{ \operatorname{Kh}, \overline{\operatorname{Kh}}, \operatorname{Kh}_{\text{odd}} , \overline{\operatorname{Kh}}_{\text{odd}} \}$ with every choice of basepoint. In (\ref{Equation: tangle homology dim bound 2}) the maximum is taken similarly over all $o, \mathbb F, \mathcal H$ and the left-hand side is to be interpreted as a function of $n$. 
\end{theorem}
In particular, Theorem \ref{Theorem: homology bound for compositions of T4} yields constant and quadratic bounds on $\operatorname{Kh} (T(4,n))$ by choosing $T$ to be the trivial tangle with 4 vertical lines.
\begin{proof}
    Let us first look at the even unreduced Khovanov homology and $n\geq 0$. Fixing $n$, a diagram $T$ and an orientation $o$, we have 
    $$
    \llbracket (\sigma_1\sigma_2\sigma_3)^n \otimes_{\operatorname{Link},o} T \rrbracket \simeq C_n \otimes_{\operatorname{Kom}} D
    $$
    where $D=\Psi \llbracket T \rrbracket [s_1]\{s_2\}$ and shifts $s_1,s_2\in \mathbb Z$ depend on $o$. It follows that  
    $$
    \operatorname{dim}_{\mathbb F} \operatorname{Kh}^{i,j}( (\sigma_1\sigma_2\sigma_3)^n \otimes_{\operatorname{Link},o} T; \mathbb F ) \leq \operatorname{rank} F^{i,j} (C_n \otimes_{\operatorname{Kom}} D ).
    $$
    By opening up the definition of tensor product and using additivity of $F^{i,j}$ we get
    $$
    F^{i,j} (C_n \otimes_{\operatorname{Kom}} D ) =\bigoplus_{a,b,c,d\in \mathbb Z} F^{i,j}(C_n^{a,b} \otimes_{\operatorname{Cob}} D^{c,d}).
    $$
    
    The complex $D$ consists of finitely many loopless Temperley-Lieb diagrams with 8 boundary points. Hence there exist $x_1,x_2,y_1,y_2$ such that $D$ is supported in degrees $[x_1,x_2]\times [y_1\times y_2]$ and 
    $\max_{i,j}\{\dimcob D^{i,j}\}$ is finite. Any pair of loopless diagrams $A,B$ can in total produce at most $4$ loops in $A\otimes_{\operatorname{Cob}}B$ which is why
    $$
    F^{i,j}(C_n^{a,b} \otimes_{\operatorname{Cob}} D^{c,d}) \not \cong 0
    $$
    implies that 
    \begin{align}
        i&=a+c \label{Equation: first eq of index bounds for rank} \\ 
        j& \in [b+d-4,b+d+4] \\
        c& \in [x_1,x_2] \\
        d&\in [y_1,y_2].\label{Equation: last eq of index bounds for rank}
    \end{align}
    Denote $K$ as the set of 4-tuples $(a,b,c,d)$ satisfying (\ref{Equation: first eq of index bounds for rank})-(\ref{Equation: last eq of index bounds for rank}). It is easy to see that 
    $$
    \# K=9(x_2-x_1+1)(y_2-y_1+1)
    $$
    and that $\# K$ does not depend on $i,j,n,o,\mathbb F$. Combining the above and Corollary \ref{Corollary: bounds on Cn dimensions} yields
    \begin{align*}
        \operatorname{dim}_{\mathbb F} \operatorname{Kh}^{i,j}((\sigma_1\sigma_2\sigma_3)^n \otimes_{\operatorname{Link},o} T; \mathbb F ) & \leq \operatorname{rank}\Big ( \bigoplus_{(a,b,c,d)\in  K} F^{i,j}(C_n^{a,b} \otimes_{\operatorname{Cob}} D^{c,d}) \Big ) \\ 
        &\leq 16 \cdot \# K \cdot \max_{n,k,l} \{ \dimcob C_n^{k,l} \}  \cdot \max_{k,l}\{\dimcob D^{k,l}\} < \infty
    \end{align*}
    proving (\ref{Equation: tangle homology dim bound 1}) for $n\geq 0$ and unreduced even Khovanov homology. Replacing $F$ with $\overline{F}$ in the proof above, yields the constant bound (\ref{Equation: tangle homology dim bound 1}) for $n\geq 0$ for the reduced even Khovanov homology with any choice of basepoint.

    As previously mentioned, odd Khovanov homology is defined similarly to even Khovanov homology, but with different sign conventions. The reduction of discrete Morse theory only cares whether the morphisms are zero, nonzero or isomorphisms --- properties which are not affected by signs. Let $A$ denote the odd unreduced Khovanov complex of $\mathbb Z$-modules
    $$
     \to F^{i-1,\ast }\llbracket (\sigma_1\sigma_2\sigma_3)^n \otimes_{\operatorname{Link},o} T \rrbracket
    \to F^{i,\ast }\llbracket (\sigma_1\sigma_2\sigma_3)^n \otimes_{\operatorname{Link},o} T \rrbracket
    \to F^{i+1,\ast }\llbracket (\sigma_1\sigma_2\sigma_3)^n \otimes_{\operatorname{Link},o} T \rrbracket
    \to
    $$
    with the odd differentials. We can use discrete Morse theory on $A$ and reverse all the edges which correspond to edges reversed in $G((\sigma_1 \sigma_2 \sigma_3)^n, \Mgr)$ and call the resulting Morse complex $MA$. This matching is acyclic since the corresponding graph is a Cartesian product of two acyclic graphs.  By observing carefully, one can see that there are isomorphisms of chain spaces 
    $$
    (MA)^{i,\ast} \cong F^{i,\ast} (C_n\otimes_{\operatorname{Kom}} D).
    $$
    It follows that the constant bound (\ref{Equation: tangle homology dim bound 1}) holds for unreduced odd Khovanov homology, $n\geq 0$. The reduced odd Khovanov homology case follows from (\ref{Equation: reduced unreduced Khovanov splitting}). 
    
    In case $n <0$, we can rotate the link around the red axis displayed in Figure \ref{Figure: composition of tangles into unoriented link} and take the mirror image. This gives out a new link diagram for which $n > 0$. Rotation does not affect the invariants and mirroring a link dualizes the homology: $\mathcal H^{i,j}(L,\mathbb F)\cong \mathcal H^{-i,-j}(L^{!},\mathbb F) $ (See \cite{khovanov1999categorification} for even and \cite{PutyraLobawskiOddMirror} for odd Khovanov homology). Hence (\ref{Equation: tangle homology dim bound 1}) holds for negative $n$ as well.
    
    For all $n$ and all theories, the quadratic bound (\ref{Equation: tangle homology dim bound 2}) can be deduced from the constant bound (\ref{Equation: tangle homology dim bound 1}) by finding a quadratic bound on the support of the whole complex.
\end{proof}

\subsection{Homology results: recursions and vanishing of \texorpdfstring{$\operatorname{Kh}(T(4,-n))$}{H(T(4,-n))}}




The following powerful proposition states that $\mathfrak a_n$ and $\mathfrak c_n$ commute with the differentials and are thus in some sense chain maps in their respective degrees. If similar commutation would hold for $\mathfrak b_n$, one could use it to obtain all Khovanov homology of $T(4,n)$. While we do expect $\mathfrak b_n$ to follow along, rigorously proving this would likely require some new ideas compared to our brute force approach in Section \ref{Section: proof of morphism correspondences}.

\begin{proposition}\label{Proposition: ac matrix element equations}
    Given $t_{\bA}(w),t_{\bA}(u) \geq 1$, the matrix elements of $C_n$ and $C_{n+4}$ are equal:
    $$\partial_{u,w}=\partial_{\mathfrak a_n(u),\mathfrak a_n(w)}.$$
    Given $t_{\bC}(x)\geq 6$ and $t_{\bC}(y)\geq 1$, the matrix elements of $C_n$ and $C_{n+4}$ are equal: 
    $$\partial_{y,x}=\partial_{\mathfrak c_n(y),\mathfrak c_n(x)}.$$

\end{proposition}

Let us now restate and prove Theorem \ref{Theorem: recursions for T4 homology} from the introduction.
\homologyRecThm*
If necessary, Theorem \ref{Theorem: recursions for T4 homology} could be extended to account for all possible tangle compositions $(\sigma_1 \sigma_2 \sigma_3)^n \otimes_{\operatorname{Link},o} T$ as Theorem \ref{Theorem: homology bound for compositions of T4} did. For the sake of simplicity we only consider the trivial decompositions: 4-strand torus links.
\begin{proof}
Let $I$ be the trivial 8-ended tangle diagram with  4 vertical lines and 0 crossings and let $o$ be the orientation so that 
$$
T(4,-n)=(\sigma_1\sigma_2\sigma_3)^n \otimes_{\text{Link,o}} I
$$
as links. There is a homotopy equivalence 
$$
\llbracket T(4,-n) \rrbracket \simeq C_n \otimes_{\operatorname{Kom}} I'
$$
where $I'$ is the complex with one object $I$ in the homological degree $0$ and quantum shift $0$. It follows that $\operatorname{Kh}^{i,j}\big(T(4,-n) \big)$ can be computed from the sequence
$$
\bigoplus_{a,b\in \mathbb Z} F^{i-1,j} \big(C_n^{a,b} \otimes_{\operatorname{Cob}} I \big)\to 
\bigoplus_{a,b\in \mathbb Z} F^{i,j} \big(C_n^{a,b} \otimes_{\operatorname{Cob}} I \big)\to 
\bigoplus_{a,b\in \mathbb Z} F^{i+1,j} \big( C_n^{a,b} \otimes_{\operatorname{Cob}} I\big)  
$$
which is equivalent to the sequence
\begin{equation} \label{Torus link homology domain sequence}
\bigoplus_{b=j-4}^{j+4} F^{i-1,j} \big(C_n^{i-1,b} \otimes_{\operatorname{Cob}} I \big)\to 
\bigoplus_{b=j-4}^{j+4} F^{i,j} \big(C_n^{i,b} \otimes_{\operatorname{Cob}} I \big)\to 
\bigoplus_{b=j-4}^{j+4} F^{i+1,j} \big( C_n^{i+1,b} \otimes_{\operatorname{Cob}} I \big).      
\end{equation}

Now, let $n,i,j\in \mathbb Z$ such that $n\geq 0$ and $-2n+2i-j\geq 14$. It follows that $t_{\bA}(w)\geq 1$ for all $w\in U_n$ with $(h(w),q(w))\in [i-1,i+1] \times [j-4,j+4]$. Thus, by Proposition \ref{Proposition: ac matrix element equations} there is a commutative diagram in the category $\operatorname{Mat}(\Cob(8))$:
\begin{equation} \label{commutative diagram: a_n diagram in cob}
\begin{tikzcd}[row sep=2em, column sep=3em]
\bigoplus_{b=j-4}^{j+4}C_n^{i-1,b}\arrow[r] \arrow[d] & \bigoplus_{b=j-4}^{j+4}C_n^{i,b} \arrow[r] \arrow[d] & \bigoplus_{b=j-4}^{j+4}C_n^{i+1,b} \arrow[d] \\
\bigoplus_{b=j-4}^{j+4}C_{n+4}^{i-1,b}\{12\} \arrow[r]           & \bigoplus_{b=j-4}^{j+4}C_{n+4}^{i,b}\{12\} \arrow[r]           & \bigoplus_{b=j-4}^{j+4}C_{n+4}^{i+1,b}\{12\}          
\end{tikzcd}    
\end{equation}
where the vertical arrows are isomorphisms induced by $\mathfrak a_n$. Applying $F(- \otimes_{\operatorname{Cob}} I)$ to (\ref{commutative diagram: a_n diagram in cob}) yields commutating isomorphisms between (\ref{Torus link homology domain sequence}) and 
\begin{equation*}\label{Torus link homology codomain sequence}
\bigoplus_{b=j-4}^{j+4} F^{i-1,j} \big(C_{n+4}^{i-1,b} \{12\}\otimes_{\operatorname{Cob}} I \big)\to 
\bigoplus_{b=j-4}^{j+4} F^{i,j} \big(C_{n+4}^{i,b} \{12\}\otimes_{\operatorname{Cob}} I \big)\to 
\bigoplus_{b=j-4}^{j+4} F^{i+1,j} \big( C_{n+4}^{i+1,b} \{12\}\otimes_{\operatorname{Cob}} I \big).      
\end{equation*}
This sequence gives out the homology group $H^{i,j}(F(C_{n+4} \{12\} \otimes_{\operatorname{Kom}} I'))$ which can be seen isomorphic to $\operatorname{Kh}^{i,j-12}\big( T(4,-n-4)\big)$. Hence  $\operatorname{Kh}^{i,j}\big(T(4,-n)\big) \cong \operatorname{Kh}^{i,j-12}\big(T(4,-n-4)\big)$. The same proof works for the reduced Khovanov homology: simply replace $F$ with $\overline{F}$ and $\operatorname{Kh}$ with $\overline{\operatorname{Kh}}$. 

The second recursion is deduced from Proposition \ref{Proposition: ac matrix element equations} analogously.
\end{proof}

To prove a vanishing result for $T(4,-n)$, we need a quick lemma whose proof is postponed to Section \ref{Section: proof of vertex correspondences}.
\begin{lemma}\label{Lemma: minimal tA,tC}
    For all $w\in U$ it holds that $t_{\bA}(w),t_{\bC}(w) \geq -\frac{3}{2}$ .
\end{lemma}

\begin{proposition}\label{Proposition: vanishing bounds for T4 homology}
    The Khovanov homology groups of $4$-strand torus link vanishes at $(i,j)$-bidegree:
    $$
    \operatorname{Kh}^{i,j} (T(4,-n)) \cong 0, \qquad \overline{\operatorname{Kh}}^{i,j} (T(4,-n)) \cong 0
    $$
    whenever $n\geq 0$ and either $-2n+2i-j<-6$ or $-9n+4i-3j<-17$.
\end{proposition}
\begin{proof}
    Let $n\geq 0$ and take the contraposition that $\operatorname{Kh}^{i,j}(T(4,-n)) \not \cong 0$. By the proof of Theorem \ref{Theorem: recursions for T4 homology} this homology group is computed from Sequence \ref{Torus link homology domain sequence} and hence
    $$
    \bigoplus_{b=j-4}^{j+4} F^{i,j} \big(C_n^{i,b} \otimes_{\operatorname{Cob}} I \big) \not \cong 0.
    $$
    One of the $C_{n}^{i,j}$ must be non-zero and so there exist $w\in U_n$ with $h(w)=i$ and $q(w)\in [j-4,j+4]$. Lemma \ref{Lemma: minimal tA,tC} now yields:
    \begin{align*}
    -\frac{3}{2} &\leq t_{\bA}(w)\leq -\frac{1}{2} n+\frac{1}{2}i -\frac{1}{4}j \\
    -\frac{3}{2} &\leq t_{\bC}(w)\leq -\frac{9}{4} n+i -\frac{3}{4}j + \frac{11}{4}
    \end{align*}
    which can be massaged to $-2n+2i-j\geq-6$ and $-9n+4i-3j \geq -17$ proving the claim. The reduced case is proven analogously.   
\end{proof}

Theorem \ref{Theorem: T4 homology figures}, which states that Figures \ref{Figure: T4 middle}, \ref{Figure: T4 homology lowest degrees} and \ref{Figure: T4 homology highest degrees} describe $\operatorname{Kh}^{i,j} (T(4,-n))$, can now be obtained as a corollary of Theorem \ref{Theorem: recursions for T4 homology}, Proposition \ref{Proposition: vanishing bounds for T4 homology} and some computer data.

\begin{proof}[Proof of Theorem \ref{Theorem: T4 homology figures}]
Given $n\geq 83$, at least one of the recursions of Theorem  \ref{Theorem: recursions for T4 homology} or one of the vanishing results of Proposition \ref{Proposition: vanishing bounds for T4 homology} applies for every $n,i,j$ for which Figures \ref{Figure: T4 middle}, \ref{Figure: T4 homology lowest degrees} and  \ref{Figure: T4 homology highest degrees} claim $\operatorname{Kh}^{i,j}(T(4,-n))$. This can be verified with the help of \texttt{Lean}'s linarith tactic, by proving case-by-case that the inequalities of the figures  imply the inequalities of Theorem  \ref{Theorem: recursions for T4 homology} and Proposition \ref{Proposition: vanishing bounds for T4 homology}. Hence, it is straightforward to see that if the figures are accurate for $n=28,\dots,82$, then by induction they are accurate for $n\geq 83$. Luckily, Khovanov homology of small links is computable: we use Lewark's program \texttt{Khoca} \cite{Khocasoftwarearticle} to compute $\operatorname{Kh}(T(4,-28)),\dots,\operatorname{Kh}(T(4,-82))$ and check (with simple \texttt{Python} scripts) that the results match with Figures \ref{Figure: T4 middle}, \ref{Figure: T4 homology lowest degrees} and  \ref{Figure: T4 homology highest degrees}. 
\end{proof}

\subsection{Comparison with GOR-conjecture}
In order to match our conventions to the Gorsky-Oblomkov-Rasmussen conjecture and to the rational Gordian bounds in \cite{lewark2024khovanovhomologyrefinedbounds} we will consider positive torus links in this subsection and the following one.

The stable Khovanov homology of $n$-strand torus links (not to be confused with the stable Khovanov homotopy type \cite{Lipshitz2014}) is defined as 
$$
\operatorname{Kh} (T(n,\infty)) = \lim_{m\to \infty} (\operatorname{Kh}(T(n,m))\{-(n-1)(m-1)+1\} ). 
$$
Stošić showed that this is a well-defined limit as for every bigrading $(i,j)$ the shifted homology groups are eventually isomorphic \cite{Stosic2007}. 

\begin{conjecture}[\cite{Gorsky2013}]
    The stable Khovanov homology $\operatorname{Kh} (T(n,\infty))$ is the dual to the homology of the Koszul complex from the sequence $s_0,\dots , s_{n-1} \in \mathbb Z[x_0,\dots, x_{n-1}]$ where
    $$ 
    s_k = \sum_{i=0}^k x_ix_{k-i}. 
    $$
    The monomials $x_i$ are bigraded with $h(x_i)=2i$ and $q(x_i)=2i+2$ and the Koszul differential respects the quantum grading and lowers the homological grading by $1$.
\end{conjecture}

Following \cite{Gorsky2013}, we denote the dual Koszul homology as $\operatorname{Kh}_{\operatorname{alg}} (T(n,\infty))$. The next two propositions verify that our computations for $\operatorname{Kh}(T(4,\infty))$ coincide with $\operatorname{Kh}_{\operatorname{alg}}(T(4,\infty))$ as conjectured. 

\begin{proposition}\label{Proposition: GOR-comparison 4 torsion}
    Both $\operatorname{Kh} (T(4,\infty))$ and $\operatorname{Kh}_{\operatorname{alg}} (T(4,\infty))$ have 4-torsion at bidegrees $(9+4k,14+6k)$ for $k\geq 0$.
\end{proposition}
\begin{proof}
    Let $m\geq 2$ and consider the monomials $x_0x_2^m$ in $\operatorname{Kh}_{\operatorname{alg}} (T(4,\infty))$ at Koszul degree $0$. 
    The homology at Koszul degree $0$ is defined as
    $$\mathbb Z[ x_0,x_1,x_2,x_3]/ (x_0^2, \ 2x_0x_1, \  x_1^2+2x_0x_2, \  2x_0x_3+2x_1x_2).$$
    Therefore in order to prove that $x_0x_2^m$  represents $4$-torsion, it suffices to show that for all $z\in \mathbb Z $, there exist $p_0,p_1,p_2,p_3\in \mathbb Z[ x_0,x_1,x_2,x_3]$ such that
    \begin{equation}\label{Equation: monomial for Koszul homology}
    zx_0x_2^m= p_0x_0^2+ p_12x_0x_1+  p_2(x_1^2+2x_0x_2) +p_3( 2x_0x_3+2x_1x_2).
    \end{equation}
    if and only if $z\in 4 \mathbb Z$. Let $z\in \mathbb Z$ and suppose that $p_0,p_1,p_2,p_3$ satisfy Equation \ref{Equation: monomial for Koszul homology}. On the right hand side of Equation \ref{Equation: monomial for Koszul homology}, the only way to generate a $zx_0x_2^m$ monomial is to ensure that $p_2$ contains a $\frac{z}{2}x_2^{m-1}$ monomial. As a byproduct a $\frac{z}{2}x_1^2 x_2^{m-1}$ term is also generated and to offset it, $p_3$ needs to have a $-\frac{z}{4}x_1x_2^{m-2}$ monomial. Now $\frac{z}{4}$ has to make sense as an integer, so we obtain that $z\in 4\mathbb Z$.  On the other hand, given $z\in 4\mathbb Z$, a solution for Equation \ref{Equation: monomial for Koszul homology} is easy to generate:
    $$
    zx_0x_2^m= 0 \cdot x_0^2+ \left(\frac{z}{4}x_2^{m-2}x_3\right)2x_0x_1 +\frac{z}{2}x_2^{m-1}(x_1^2+2x_0x_2) -\frac{z}{4}x_1x_2^{m-2} ( 2x_0x_3+2x_1x_2).
    $$
    
    The monomial $x_0 x_2^{m}$ has bidegree $(4m,2+6m)$ so when dualizing the homological degree of the torsion groups increases by 1 yielding the result for $\operatorname{Kh}_{\operatorname{alg}} (T(n,\infty))$. 
    Observing Figures \ref{Figure: T4 middle} and \ref{Figure: T4 homology highest degrees}, mirroring both degrees $(i,j)\mapsto (-i,-j)$, dualizing and then renormalizing with the quantum degrees one can see that every $\operatorname{Kh}^{9+4k,14+6k}(T(4, \infty))$ also contains a $\mathbb Z / 4 \mathbb Z$ summand for $k\geq 0$. 
\end{proof}

\begin{proposition} \label{Proposition: GOR-comparison over F_2}
    Over $\mathbb F_2$, our stable Khovanov homology computations in Figure \ref{Figure: T4 middle} agree with $\operatorname{Kh}_{\operatorname{alg}} (T(4,\infty))$. More precisely, as $\mathbb F_2$ vector spaces
    $$
    \operatorname{Kh}^{i,j} (T(4,\infty);\mathbb F_2) \cong \operatorname{Kh}_{\operatorname{alg}}^{i,j} (T(4,\infty); \mathbb F_2)
    $$
    whenever $i\geq 42$ and $j\geq \frac{3}{2}i-1$.
\end{proposition}
\begin{proof}
    By Theorem 1.5 of \cite{Gorsky2013},  $\operatorname{Kh}_{\operatorname{alg}} (T(n,\infty); \mathbb F_2)$ has the following Poincaré series
    $$
    P_n(t,q ; \mathbb F_2)= \prod_{i=0}^{n-1} \frac{1+t^{2i+1}q^{2i+4}}{1-t^{2i}q^{2i+2}} \prod_{i=0}^{\lfloor \frac{n-1}{2} \rfloor } \frac{1-t^{4i}q^{4i+4}}{1+t^{4i+1}q^{4i+4}}  
    $$
    where $t$ represents the homological gradings and $q$ the quantum gradings. In the case of $n=4$ this can be simplified and expanded as a power series at $(0,0)$, so that 
    $$
    P_4(t,q; \mathbb F_2) = \frac{(1+q^2) (1+t^2q^4) (1+t^3q^6) (1+t^7q^{10}) }{(1-t^4q^6)(1-t^6q^8)}= \sum_{a,b\geq 0 } c_{a,b} t^aq^b
    $$
    for some $c_{a,b}\geq 0$. Consider another power series $P'(t,q)$, where the expansion of $(1-t^6q^8)^{-1}$ is truncated:
    $$
    P'(t,q) = \frac{(1+q^2) (1+t^2q^4) (1+t^3q^6) (1+t^7q^{10}) }{1-t^4q^6}\prod_{i=0}^6 (t^6q^8)^i = \sum_{a,b \geq 0} d_{a,b} t^aq^b.
    $$
    Straightforward calculations show that $b- \frac{3}{2}a \geq -1$ implies $c_{a,b}=d_{a,b}$ and that $c_{a,b}=d_{a,b}=0$ whenever $b-\frac{3}{2}a \geq 5 $. Any monomial contributing to $d_{a,b}$ with $a\geq 49$ must involve a $(t^4 q^6)^k$ term with $k \geq 1$ from the expansion of $(1-t^4q^6)^{-1}$. Hence $d_{a,b}=d_{a+4,b+6}$ whenever $a\geq 45$. 

    Denote $C(a,b)=\dim_{\mathbb F_2} \operatorname{Kh}^{a,b}(T(4,\infty);\mathbb F_2)$. Observing Figure \ref{Figure: T4 middle} and Figure $\ref{Figure: T4 homology highest degrees}$ for the homological degree $-41$, using the universal coefficient theorem and then mirroring the degrees we obtain that for $k\geq 11$:
    \begin{alignat*}{4}
        &C(4k-2,6k-4)&&=8, \qquad C(4k-2,6k-2)&&=7, \qquad C(4k-2,6k)&&=3,   \\
        &C(4k-1,6k-2)&&=8, \qquad C(4k-1,6k)&&=6, \qquad C(4k-1,6k+2)&&=2,   \\
        &C(4k,6k)&&=8, \qquad C(4k,6k+2)&&=5, \qquad C(4k,6k+4)&&=1,   \\
        &C(4k+1,6k+2)&&=7, \qquad C(4k+1,6k+4)&&=4, \qquad C(4k+1,6k+6)&&=1, 
    \end{alignat*}
    and $C(a,b)=0$ whenever $b\geq \frac{3}{2}a +5 $. We verify that $c_{a,b}=C(a,b)$ for all $42 \leq a \leq 49$ and $-1 \leq b- \frac{3}{2}a < 5$ by computing the coefficients $c_{a,b}$ from the rational function $P_4(t,q; \mathbb F_2)$ with \texttt{Mathematica}. The claim follows from induction.   
\end{proof}

\subsection{Rational replacements and \texorpdfstring{$\mathbb Z[G]$}{ZG}-complexes}

A unified way of incorporating several versions of Khovanov homology was introduced by Naot as $\mathbb Z[G]$-complexes \cite{zbMATH05118580}, where $G$ is a formal variable (not a group). Given a knot diagram $K$, one assigns to it $\operatorname{BN}(K)$ which is the reduced Khovanov complex associated to the Frobenius algebra $\mathbb Z[G][X]/(X^2+GX)$ over $\mathbb Z[G]$. For us, the key properties are:
\begin{enumerate}
    \item The chain homotopy type of $\operatorname{BN}(K)$ is a knot invariant. \label{ZG Property 1}
    \item Setting $G=0$ on $\operatorname{BN}(K)$ and taking the homology recovers the reduced Khovanov homology:
    $$H(\operatorname{BN}(K)\otimes_{\mathbb Z[G] } \mathbb Z [G]/G)\cong_{\mathbb Z \operatorname{Mod}} \overline{\operatorname{Kh}}(K).$$
    (In spectral sequence lingua, this is the fact that the reduced Lee spectral sequence starts at the the reduced Khovanov homology.)
    \item \label{ZG Property 3}Setting $G=1$ collapses the homology, except for one copy of $\mathbb Z$ in homological degree $0$: $$H(\operatorname{BN}(K)\otimes_{\mathbb Z[G]} \mathbb Z[G] /(G-1)) \cong_{\mathbb Z \operatorname{Mod}} \mathbb Z.$$
    (The reduced Lee spectral sequence abuts to the trivial Lee homology.)
    \item The $G$-torsion order on $H(\operatorname{BN}(K))$ gives an lower bound on the rational Gordian distance. \label{ZG Property 4}
\end{enumerate}
A diagrammatic introduction to $\mathbb Z[G]$-complexes and proofs to  \ref{ZG Property 1}-\ref{ZG Property 3} can be found in \cite{iltgen2021khovanovhomologyrationalunknotting} and \ref{ZG Property 4} will be expanded upon after the following algebraic result. Since in $\mathbb Z[G]$-modules the brackets $[\cdot ]$ are used for the formal variable, we will denote the homological and quantum gradings with superscripts of $i$ and $j$ respectively.    

\begin{proposition} \label{Proposition: wierd piece in ZG homology of T4}
    For any $n\geq2$ there is a split of $\mathbb Z[G]$-complexes $\operatorname{BN}(T(4,2n+1))\simeq \mathcal D_n\oplus \mathcal E_n$ where $\mathcal D_n$ is supported in homological gradings $[0,4n]$.  
     The four term complex $\mathcal{E}_n$ is
    $$
\begin{tikzcd}[row sep=0.4em, column sep=3em]
                                                  & {i^{4n+1}j^{12n+4}\mathbb Z[G]} \arrow[rd, "2"]  &                \\
{i^{4n}j^{12n}\mathbb Z[G]} \arrow[ru, "G^2"] \arrow[rd, "2G"'] & \oplus                          & {i^{4n+2}j^{12n+4}\mathbb Z[G]} \\
                                                  & {i^{4n+1}j^{12n+2}\mathbb Z[G]} \arrow[ru, "-G"'] &               
\end{tikzcd}
    $$ 
\end{proposition}
\begin{proof}
For any $n\geq 4$, the reduced Khovanov homology $\overline{\operatorname{Kh}}(T(4,2n+1))$  in homological degrees $[4n,4n+2]$ is given by:
    \begin{center}
        \begin{tabular}{|c|c|c|c|}
    \hline
    \diagbox[width=4.5em, height=1.4em]{\small$ j$}{\small$ i$} & $4n$ & $4n+1$ & $4n+2$ \\ \hline
    $12n+4$ &  &  & $\mathbb Z_2$ \\ \hline
    $12n+2$ & $\mathbb Z$ & $\mathbb Z$ &  \\ \hline
    $12n$ & $\mathbb Z \oplus \mathbb Z_2$ &  &  \\ \hline
    \end{tabular}
    \end{center}
and $\overline{\operatorname{Kh}}^{i,\ast}(T(4,2n+1))$ vanishes for $i\geq4n+3$. By duality of Khovanov homology, this can be obtained from $\overline{\operatorname{Kh}}^{i,\ast}(T(4,-2n-1))$ for $i\leq -4n+1$. In these degrees, the reduced Khovanov homology of negative torus knots can be obtained from Theorem \ref{Theorem: recursions for T4 homology}, Proposition \ref{Proposition: vanishing bounds for T4 homology} and computer data analogously to the proof of Theorem \ref{Theorem: T4 homology figures}.

By using Gaussian elimination and Smith normal form for the zero degree morphisms of $\operatorname{BN}(T(4,2n+1))$  one can obtain a $\mathbb Z[G]$ chain complex $(A,\partial)$  which is chain homotopic $\operatorname{BN}(T(4,2n+1))\simeq A$. Moreover chain spaces and zero degree morphisms of $A$ are of the form 
$$
\begin{tikzcd}[row sep=0em, column sep=1.5em]
       {} \arrow[r, "0"]           & {i^{4n}j^{12n+2}\mathbb Z [G]} \arrow[rrddd, "0"] &  &                                                   &  &                                  \\
                  & \oplus                                            &  & {i^{4n+1}j^{12n+4}\mathbb Z [G]} \arrow[rrd, "2"] &  &                                  \\
{} \arrow[r, "0"] & {i^{4n}j^{12n}\mathbb Z [G]}                      &  & \oplus                                            &  & {i^{4n+2}j^{12n+4}\mathbb Z [G].} \\
                  & \oplus                                            &  & {i^{4n+1}j^{12n+2}\mathbb Z [G]}                  &  &                                  \\
      {} \arrow[r, "2"]            & {i^{4n}j^{12n}\mathbb Z [G]}                      &  &                                                   &  &                                 
\end{tikzcd}
$$
Additionally the differentials $\partial^{4n}$ and $\partial^{4n+1} $ carry higher powers of $G$ which are determined by the differences in quantum degree, that is,
$$
\partial^{4n}= \begin{bmatrix}
a_1G & a_2G^2 & a_3G^2 \\
0 & a_4 G & a_5 G 
\end{bmatrix}, \qquad  \partial^{4n+1}= \begin{bmatrix}
2 & bG 
\end{bmatrix}
$$
for some $a_1,\dots, a_5,b\in \mathbb Z $. From $\partial^{4n+1} \partial^{4n}=0$ one can deduce that $a_1=0$. Subsequently $\partial^{4n} \partial^{4n-1}=0$  yields $a_3=a_5=0$. 

Since setting $G=1$ collapses the homology in these degrees, $b$ has to be odd and furthermore $a_2$ and $a_4$ have to be non-zero. Hence $A=\mathcal{D}_n \oplus \mathcal{E}_n'$ where $\mathcal{E}_n'$  is 
$$
\begin{tikzcd}[row sep=0.4em, column sep=3em]
                                                  & {i^{4n+1}j^{12n+4}\mathbb Z[G]} \arrow[rd, "2"]  &                \\
{i^{4n}j^{12n}\mathbb Z[G]} \arrow[ru, "a_2G^2"] \arrow[rd, "a_4G"'] & \oplus                          & {i^{4n+2}j^{12n+4}\mathbb Z[G]}. \\
                                                  & {i^{4n+1}j^{12n+2}\mathbb Z[G]} \arrow[ru, "bG"'] &               
\end{tikzcd}
$$
Since $2a_2=-a_4b$ and $b$ is odd, $a_4=2c$ for some $c\in \mathbb Z$. Now $(-b G,2 )$ is cycle in homological degree $4n+1$ and it is non-trivial if $c\neq \pm 1$. Moreover, if $c\neq \pm 1$ the element $(-bG,2)$ would continue to represent a nontrivial cycle after setting $G=1$. Hence $c$ has to be $\pm 1$ and $\mathcal{E}_n'$ is of the form  
$$
\begin{tikzcd}[row sep=0.4em, column sep=3em]
                                                  & {i^{4n+1}j^{12n+4}\mathbb Z[G]} \arrow[rd, "2"]  &                \\
{i^{4n}j^{12n}\mathbb Z[G]} \arrow[ru, "\mp b G^2"] \arrow[rd, "\pm 2 G"'] & \oplus                          & {i^{4n+2}j^{12n+4}\mathbb Z[G]}. \\
                                                  & {i^{4n+1}j^{12n+2}\mathbb Z[G]} \arrow[ru, "bG"'] &               
\end{tikzcd}
$$
which is isomorphic to the complex $\mathcal{E}_n$ we were aiming for.

The Khovanov homologies $\overline{\operatorname{Kh}}^{i_1,\ast} (T(4,5))$ and $\overline{\operatorname{Kh}}^{i_2,\ast} (T(4,7))$ for homological degrees $i_1\geq 8$ and $i_2\geq 12$ are given by 
    \begin{center}
        \begin{tabular}{|c|c|c|c|}
    \hline
    \diagbox[width=4.5em, height=1.4em]{\small$ j$}{\small$ i_1$} & $8$ & $9$ & $10$ \\ \hline
    $28$ &  &  & $\mathbb Z_2$ \\ \hline
    $26$ &  & $\mathbb Z$ &  \\ \hline
    $24$ & $\mathbb Z $ &  &  \\ \hline
    \end{tabular}
    $\hspace{1 cm}$
        \begin{tabular}{|c|c|c|c|}
    \hline
    \diagbox[width=4.5em, height=1.4em]{\small$ j$}{\small$ i_2$} & $12$ & $13$ & $14$ \\ \hline
    $40$ &  &  & $\mathbb Z_2$ \\ \hline
    $38$ & $\mathbb Z$ & $\mathbb Z$ &  \\ \hline
    $36$ & $\mathbb Z $ &  &  \\ \hline
    \end{tabular}
    \end{center}
respectively. Hence cases $n=2$ and $n=3$ can be proven similarly as $n\geq 4$ but with one or two fewer initial steps needed.
\end{proof}

Let $R$ be a ring and $A$ an $R[G]$ module. For $a\in A$ we denote $\operatorname{ord}_G(a)=\infty$ if $G^na\neq 0$ for all $n\geq 0$ and  $\operatorname{ord}_G(a)=\min \{ n\geq 0 \mid G^n a=0\}$ otherwise. Furthermore for a knot $K$ and $i\in \mathbb Z$ we denote
$$
\mathfrak u_i(K; R)=\max \{ \operatorname{ord}_G(a) \mid a\in H^{i,\ast} (\operatorname{BN}(K) \otimes_{\mathbb Z[G]} R[G]), \  \operatorname{ord}_G(a)< \infty \}.
$$

\begin{figure}
    \centering
    \input{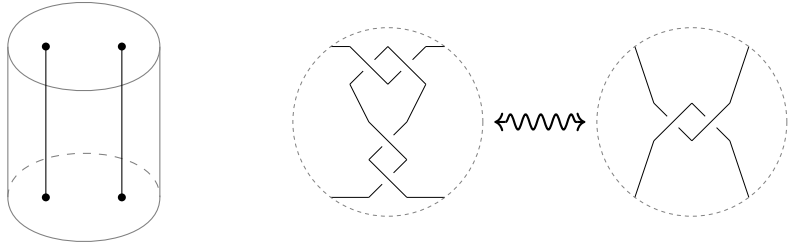}
    \caption{On the left, the pair $(D^2\times [0,1], \{(-\frac{1}{2},0)  \cup (\frac{1}{2},0) \}\times [0,1] )$ which is our model for rational tangles. In the center and on the right we see two rational tangles. Switching one with the other in some region of a knot is a proper rational tangle replacement, since both rational tangles share the same connectivity. }
    \label{Figure: Rational tangle and replacement}
\end{figure}

A 4-ended tangle $T$ is called \term{rational}, if $(B^3,T)$ is homeomorphic as a pair to $(D^2\times [0,1], \{(-\frac{1}{2},0)  \cup (\frac{1}{2},0) \}\times [0,1] )$, see Figure \ref{Figure: Rational tangle and replacement}. Two knots $K_1, K_2$ are said to be related by a \term{rational tangle replacement} if an isotopy representative of $K_1$ can be turned into an isotopy representative of $K_2$  by replacing one rational tangle with another inside a small $B^3$. The rational tangle replacement is said to be \term{proper}, if both rational tangles share the same connectivity. The \term{rational Gordian distance} between knots $K_1$ and $K_2$ is the minimal number of proper rational tangle replacements needed to transform $K_1$ into $K_2$. This defines a metric on knots and its values $u_q(K_1,K_2)$ are non-negative integers. 

\begin{theorem}[\cite{lewark2024khovanovhomologyrefinedbounds}, Theorem 1.2, Proposition 3.31] \label{Theorem: LMZ torsion bound for rational distance}
    For any knots $K_1,K_2$ and ring $R$
    $$
    \max_{i\in \mathbb Z} |\mathfrak u_i(K_1; R)-\mathfrak u_i(K_2; R)|\leq u_q(K_1,K_2).
    $$
\end{theorem}

Let us now restate  Proposition \ref{Proposition: rational distance 2 if not positive enough} from the introduction and prove it using Proposition \ref{Proposition: wierd piece in ZG homology of T4} and Theorem \ref{Theorem: LMZ torsion bound for rational distance}.

\ratDistProp*

\begin{proof}
    Let $K$ be a knot admitting a diagram with $c_+\leq 4n$ positive crossings. Thus chain homotopy class of $\operatorname{BN}(K)$ has a representative which vanishes at homological degree $4n+1$ and in particular $\mathfrak u_{4n+1}(K; \mathbb F_2)=0$. On the other hand Proposition \ref{Proposition: wierd piece in ZG homology of T4} yields $\mathfrak u_{4n+1}(T(4,2n+1); \mathbb F_2)=2$ which combined with Theorem \ref{Theorem: LMZ torsion bound for rational distance} proves $u_q(T(4,2n+1),K)\geq 2$. 
    
    The second claim similarly follows from Proposition \ref{Proposition: wierd piece in ZG homology of T4} and Theorem \ref{Theorem: LMZ torsion bound for rational distance}. Without loss of generality, one of the following cases applies:
    \begin{itemize}
        \item If $n>m>0$, then $\mathfrak u_{4n+1}(T(4,n); \mathbb F_2)=2$ and $\mathfrak u_{4n+1}(T(4,m); \mathbb F_2)=0$ which suffices for $u_q(T(4,n),T(4,m))\geq 2$.
        \item If $n>0>m$, then $T(4,m)$ has a diagram with $0$ positive crossings and the first claim can be applied.
        \item The set of rational tangles is closed under taking mirror images so by reversing the orientations throughout, the first case also yields $u_q(T(4,n),T(4,m))\geq 2$ whenever $0>m>n$.
    \end{itemize}
\end{proof}

\section{Proof of Proposition \ref{Proposition: abc bijections on cells}} \label{Section: proof of vertex correspondences}
In this section, we will prove Proposition \ref{Proposition: abc bijections on cells} which introduces the correspondences $\mathfrak{a}_n, \mathfrak{b}_n, \mathfrak{c}_n$. The first order of business is to introduce a purely combinatorial set of formal words 
$$
W\subset \mathcal{S}=\bigcup_{n=0}^\infty \{ 0,\bzero, \rzero,1,\bone,\rone, \bullet, \bbullet,\rbullet \}^n
$$ 
and prove that $W=U$ where $U= \bigcup_n U_n$ and $U_n$ is the set of unmatched cells of $G(\Psi \llbracket (\sigma_1\sigma_2 \sigma_3)^n\rrbracket, \Mgr )$. After this, the functions   $\mathfrak{a}_n, \mathfrak{b}_n, \mathfrak{c}_n$ are easy do define formally on $W$ and recursive properties of $W$ will allow us to inductively verify that these functions generate one-to-one correspondences in the designated gradings.

We establish functions $g_1,\dots,g_6, I_{\bA},I_{\bB},I_{\bC}$ whose domain and codomain are the power set $\mathcal{P}(\mathcal S)$ and which are defined by
\begin{align*}
    I_{\bA}(A)&=\{a .(111)^k \mid a\in A, \ k\geq 0\} \\
    I_{\bB}(A)&=\{a .(101011\bzero0011\bzero)^k \mid a\in A, \ k\geq 0\} \\
    I_{\bC}(A)&=\{a .(01\bzero \bzero 01)^k \mid a\in A, \ k\geq 0\} \\
    g_1(A)&=\{a.000,\ a.001,\ a.010,\ a.011 \mid a\in A\} \\
    g_2(A)&=\{a,\  a.100,\ a.110,\ a.00011\bzero \mid a\in A\} \\
    g_3(A)&=\{a.100001, \ a.110001,\ a.00011\bzero001\mid a\in A\} \\
    g_4(A)&=\{a,\  a.101,\ a.101010,\ a.101011.c\mid a\in A, \ c\in\{ e,\bzero 00,\bzero 01,\bzero 10,\bzero 11\}\} \\
    g_5(A)&=\{a.101010\bzero01,\  a.101011\bzero0011\bzero001\mid a\in A\} \\
    g_6(A)&=\{a,\  a.01\bzero\mid a\in A\}.
\end{align*}
Furthermore, we assign
\begin{multline}
    W= ( g_1 \circ I_{\bA}  ) (\{e\})  \cup  (g_4 \circ I_{\bB} \circ g_2 \circ I_{\bA}  ) (\{e\})  \\
    \cup  (g_6\circ I_{\bC} \circ g_5\circ  I_{\bB}\circ g_2\circ  I_{\bA}  )(\{e\}) 
    \cup  (g_6 \circ I_{\bC} \circ g_3 \circ I_{\bA}  )(\{e\}) \label{long equation in unmatched cell classification}
    \end{multline}
where $e$ denoted the length $0$ word of $\mathcal S$. The compositions of Equation \ref{long equation in unmatched cell classification} can be viewed as the maximal paths in the schematic:
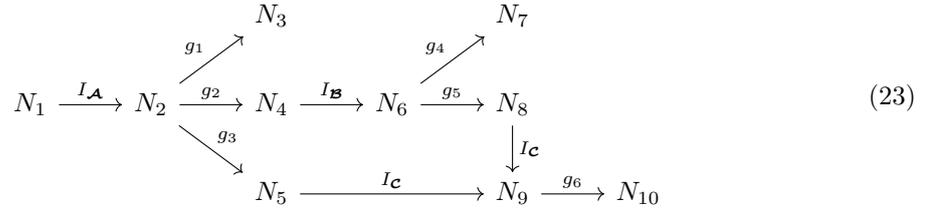
\begin{equation}    
\begin{tikzcd}
                                 &                                                          & N_3                               &                                        & N_7                            &        \\
N_1  \arrow[r, "I_{\bA}"] & N_2 \arrow[ru, "g_1"] \arrow[r, "g_2"] \arrow[rd, "g_3"] & N_4 \arrow[r, "I_{\bB}"] & N_6 \arrow[ru, "g_4"] \arrow[r, "g_5"] & N_8 \arrow[d, "I_{\bC}"] &        \\
                                 &                                                          & N_5 \arrow[rr, "I_{\bC}"]   &                                        & N_9 \arrow[r, "g_6"]           & N_{10}
\end{tikzcd} \label{composition schematic}
\end{equation}
where $N_1,\dots,N_{10}=\mathcal P(K)$ and $\{e \}$ is plugged in from $N_1$.

\begin{lemma}\label{unmatched cell classification lemma}
The set of unmatched cells $U$ coincides with the formal set of words $W$, that is, $W=U$.
\end{lemma}

\begin{proof}
    Let us begin by denoting $W_n=\{w\in W  : |w|=3n\}$ which turns our goal into showing that $W_n=U_n$ for each $n$. The base case of the induction is immediate: $W_0=\{e\}=U_0$. Now assume $W_n=U_n$. From Algorithm \ref{Algorithm: unmatched cells of Mgr} and Lemma \ref{lemma: red zeros matched} one can see  
$$
U_{n+1}=\left\{ac_1c_2c_3 \ \middle \vert \begin{array}{c}
     a\in W_n, \ c_1,c_2,c_3\in \{0,\bzero, 1\} \\
     ac_1\in \operatorname{cells}(\Psi\llbracket (\sigma_1\sigma_2\sigma_3)^{n}\sigma_1 \rrbracket), \ ac_1 \text{ is not matched in } M_{\ast,3n+1} \\
     ac_1c_2\in \operatorname{cells}(\Psi\llbracket (\sigma_1\sigma_2\sigma_3)^{n}\sigma_1 \sigma_2 \rrbracket), \ ac_1c_2 \text{ is not matched in } M_{\ast,3n+2} \\
     ac_1c_2c_3\in \operatorname{cells}(\Psi\llbracket (\sigma_1\sigma_2\sigma_3)^{n+1} \rrbracket), \ ac_1c_2c_3 \text{ is not matched in } M_{\ast,3n+3} \\
\end{array}\right\}.
$$
Hence our goal is to verify that for all $a\in W_n$ and $c_1,c_2,c_3\in \{0,\bzero, 1\}$:
\begin{equation}\label{Equivalence: induction for unmatched cells}
\begin{array}{l}
     ac_1c_2c_3\in \operatorname{cells}(\Psi\llbracket (\sigma_1\sigma_2\sigma_3)^{n+1} \rrbracket),  \text{ and }    ac_1, ac_1c_2 \text{ and }  ac_1c_2c_3 \\ \text{are not matched with Algorithm \ref{Algorithm: matching a single cell in Mgr} in the last 3 Morse layers}
    \end{array} \iff ac_1c_2c_3\in W_{n+1}.    
\end{equation}


We begin with the case $a\in (g_1 \circ I_{\bA}  ) (\{e\}) $ and consider the possible continuations $ac_1c_2c_3$ by first splitting to cases based on $c_1$, then on $c_2$ and finally on $c_3$. If $a=(111)^k010$, then both possible $c_1$ continuations $a.\bzero$ and $a.1$ will be matched
$$
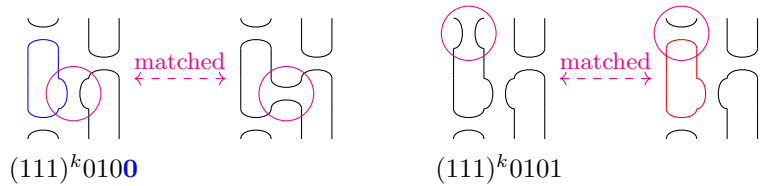
\begin{tikzpicture}[scale=0.4]
    \lineRightBendUpMore{0,0}
    \begin{scope}[blue]
    \lineRightBendDownMore{0,1}
    \lineUpBendRightMore{1,1}
    \lineUp{1,2}
    \lineRightBendUpMore{0,3}
    \lineUp{0,2}
    \lineUp{0,1}    
    \end{scope}
    
    \lineUp{2,0}
    \lineUp{3,0}    
    \lineUpBendLeftMore{2,1}
    \lineUp{3,1}
    \lineRightBendUpMore{2,2}
    \lineRightBendDownMore{2,3}
    \lineRightBendDownMore{0,4}
    \lineUp{2,3}
    \lineUp{3,3}

    \begin{scope}[shift={(7,0)}]
        
    \lineRightBendUpMore{0,0}
    \lineRightBendDownMore{0,1}
    \lineUp{2,0}
    \lineUp{3,0}
    \lineUp{0,1}
    \zeroSmoothing{1,1}
    \lineUp{3,1}
    \lineUp{0,2}
    \lineUp{1,2}
    \lineRightBendUpMore{2,2}
    \lineRightBendDownMore{2,3}
    \lineRightBendUpMore{0,3}
    \lineRightBendDownMore{0,4}
    \lineUp{2,3}
    \lineUp{3,3}
    \end{scope}

    \begin{scope}[shift={(14,0)}]
        
    \lineRightBendUpMore{0,0}
    \lineRightBendDownMore{0,1}
    \lineUp{2,0}
    \lineUp{3,0}
    \lineUp{0,1}
    \lineUpBendRightMore{1,1}
    \lineUpBendLeftMore{2,1}
    \lineUp{3,1}
    \lineUp{0,2}
    \lineUp{1,2}
    \lineRightBendUpMore{2,2}
    \lineRightBendDownMore{2,3}
    \oneSmoothing{0,3}
    \lineUp{2,3}
    \lineUp{3,3}
    \end{scope}

    \begin{scope}[shift={(21,0)}]
    \lineRightBendUpMore{0,0}
    \begin{scope}[red]
    \lineRightBendDownMore{0,1}
    \lineUpBendRightMore{1,1}
    \lineUp{1,2}
    \lineRightBendUpMore{0,3}
    \lineUp{0,2}
    \lineUp{0,1}    
    \end{scope}
    
    \lineUp{2,0}
    \lineUp{3,0}    
    \lineUpBendLeftMore{2,1}
    \lineUp{3,1}
    \lineRightBendUpMore{2,2}
    \lineRightBendDownMore{2,3}
    \lineRightBendDownMore{0,4}
    \lineUp{2,3}
    \lineUp{3,3}
    
    \end{scope}
    \draw[<->,dashed,magenta] (3.5,2) -- (6.5,2) node[midway, above, align=center] {\small matched};
    
    \draw[color=magenta] (1.5,1.5) circle [radius=0.9];
    \draw[color=magenta] (8.5,1.5) circle [radius=0.9];
    
    \draw[<->,dashed,magenta] (17.5,2) -- (20.5,2) node[midway, above, align=center] {\small matched};

    \draw[color=magenta] (14.5,3.5) circle [radius=0.9];
    
    \draw[color=magenta] (21.5,3.5) circle [radius=0.9];

    \node at (1.5,-1) {$(111)^k010\bzero$};

    \node at (15.5,-1) {$(111)^k0101$};
\end{tikzpicture}
$$
and the same argument shows that all continuations of $a=(111)^k 011$ will be matched. Similarly all possible continuations of $a=(111)^k 001$ are matched:
$$
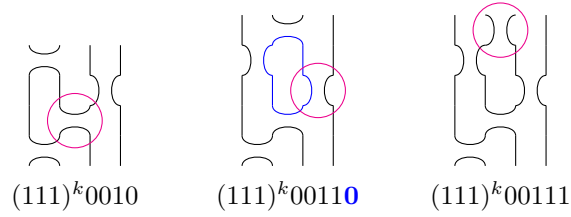
\begin{tikzpicture}[scale=0.4]
    
    \begin{scope}[shift={(0,0)}]
        
    \lineRightBendUpMore{0,0}
    \lineRightBendDownMore{0,1}
    \lineUp{2,0}
    \lineUp{3,0}
    \lineUp{0,1}
    \zeroSmoothing{1,1}
    \lineUp{3,1}
    \lineUp{0,2}
    \lineUp{1,2}
    \oneSmoothing{2,2}
    \lineRightBendUpMore{0,3}
    \lineRightBendDownMore{0,4}
    \lineUp{2,3}
    \lineUp{3,3}
    \end{scope}

\begin{scope}[shift={(7,0)}]
        
    \lineRightBendUpMore{0,0}
    \lineRightBendDownMore{0,1}
    \lineUp{2,0}
    \lineUp{3,0}
    \lineUp{0,1}
    \lineRightBendUpMore{1,1}
    
    \lineUp{3,1}
    \lineUp{0,2}
    \lineUpBendLeftMore{3,2}
    \lineUpBendRightMore{0,3}
    \lineUp{3,3}
    \begin{scope}[color=blue]
    \lineRightBendUpMore{1,4}
    \lineUp{2,3}
    \lineUp{1,2}
    \lineUpBendRightMore{2,2}
    \lineRightBendDownMore{1,2}
    \lineUpBendLeftMore{1,3}
    \end{scope}
    \lineRightBendDownMore{1,5}
    \lineUp{0,4}
    \lineUp{3,4}
    
\end{scope}

\begin{scope}[shift={(14,0)}]
        
    \lineRightBendUpMore{0,0}
    \lineRightBendDownMore{0,1}
    \lineUp{2,0}
    \lineUp{3,0}
    \lineUp{0,1}
    \lineRightBendUpMore{1,1}
    
    \lineUp{3,1}
    \lineUp{0,2}
    \lineUpBendLeftMore{3,2}
    \lineUpBendRightMore{0,3}
    \lineUp{3,3}
    \oneSmoothing{1,4}
    \lineUp{2,3}
    \lineUp{1,2}
    \lineUpBendRightMore{2,2}
    \lineRightBendDownMore{1,2}
    \lineUpBendLeftMore{1,3}
    
    \lineUp{0,4}
    \lineUp{3,4}
    
\end{scope}
    
    \draw[color=magenta] (1.5,1.5) circle [radius=0.9];
    
    \draw[color=magenta] (9.5,2.5) circle [radius=0.9];

    \draw[color=magenta] (15.5,4.5) circle [radius=0.9];

    \node at (1.5,-1) {$(111)^k0010$};
    
    \node at (8.5,-1) {$(111)^k0011\bzero$};
    
    \node at (15.5,-1) {$(111)^k00111$};

\end{tikzpicture}
$$
 where the magenta circles indicate the $i$ values of the matching arrows.  For $a=(111)^k 000$ the continuations $a.0$, $a.1\bzero$, $a.111$ are matched:
$$
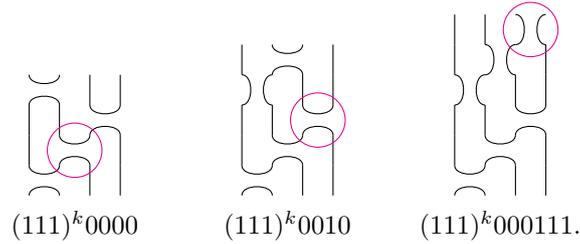
\begin{tikzpicture}[scale=0.4]

\begin{scope}[shift={(0,0)}]
                \lineUp{2,0}
        \lineUp{3,0}
        \lineUp{0,1}
        \lineUp{3,1}
        \lineUp{0,2}
        \lineUp{3,3}
        
        \zeroSmoothing{0,3}
        
        \lineRightBendDownMore{1,2}
        \lineUp{1,2}
        \lineUp{2,3}
        \zeroSmoothing{0,0}
        \lineRightBendUpMore{1,1}
        \zeroSmoothing{2,2}
        
    \end{scope}

        \begin{scope}[shift={(7,0)}]
                \lineUp{2,0}
        \lineUp{3,0}
        \lineUp{0,1}
        \lineUp{3,1}
        \lineUp{0,2}
        \lineUp{3,3}
        
        \lineUp{3,4}

        \lineRightBendUpMore{1,4}
        \lineUpBendLeftMore{1,3}
        \lineRightBendDownMore{1,2}
        \lineUp{1,2}
        \zeroSmoothing{2,2}
        \lineUp{2,3}
        \lineUp{0,4}
        
        \zeroSmoothing{0,0}
        \lineRightBendUpMore{1,1}
        \lineUpBendRightMore{0,3}
        \lineRightBendDownMore{1,5}
        
    \end{scope}

    \begin{scope}[shift={(14,0)}]
                \lineUp{2,0}
        \lineUp{3,0}
        \lineUp{0,1}
        \lineUp{3,1}
        \lineUp{0,2}
        
        \lineUp{1,2}
        \lineRightBendDownMore{1,2}
        \lineRightBendUpMore{2,2}
        \lineUpBendLeftMore{1,3}
        \lineUpBendRightMore{1,4}
        
        \lineUpBendLeftMore{2,4}
        
        \lineUp{2,3}
        \lineUp{3,3}
        \lineUp{3,4}
        \lineRightBendDownMore{2,3}
        \lineUp{0,4}
        \lineUp{0,5}
        \lineUp{1,5}

        \oneSmoothing{2,5}
        \zeroSmoothing{0,0}
        \lineRightBendUpMore{1,1}
        \lineUpBendRightMore{0,3}        
    \end{scope}
    
    \draw[color=magenta] (1.5,1.5) circle [radius=0.9];
    
    \draw[color=magenta] (9.5,2.5) circle [radius=0.9];
    
    \draw[color=magenta] (16.5,5.5) circle [radius=0.9];
    \node at (1.5,-1) {$(111)^k0000$};
    \node at (8.5,-1) {$(111)^k0010$};
    \node at (15.5,-1) {$(111)^k000111.$};

\end{tikzpicture}
$$

By running Algorithm \ref{Algorithm: matching a single cell in Mgr} on the cell $(111)^k00011\bzero$ we can see that the continuation $11\bzero$ ends up being unmatched:
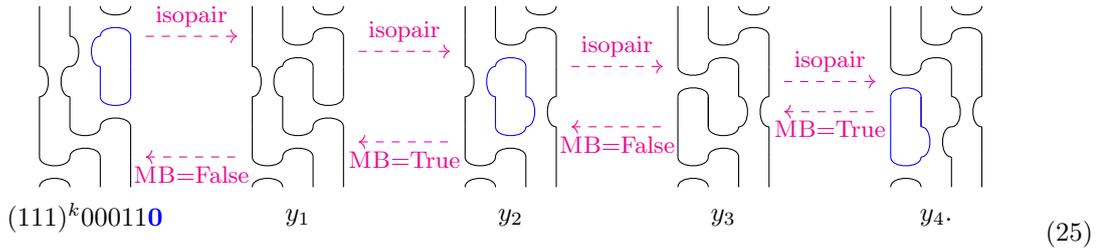
\begin{equation} \label{Chain of pictures: unmatching of 000110}
\begin{tikzpicture}[scale=0.4]
        \lineUp{2,0}
        \lineUp{3,0}
        \lineUp{0,1}
        \lineUp{3,1}
        \lineUp{0,2}
        \lineUp{1,2}
        \begin{scope}[color=blue]
        \lineRightBendUpMore{2,5}
        \lineUpBendLeftMore{2,4}
        \lineRightBendDownMore{2,3}
        \lineUp{2,3}
        \lineUp{3,3}
        
        \lineUp{3,4}
        \end{scope}
        \lineUp{0,4}
        \lineUp{0,5}
        \lineUp{1,5}

        \zeroSmoothing{0,0}
        \zeroSmoothing{1,1}
        \lineRightBendUpMore{2,2}
        \oneSmoothing{0,3}
        \lineUpBendRightMore{1,4}
        \lineRightBendDownMore{2,6}

        \begin{scope}[shift={(7,0)}]
            
        \lineUp{2,0}
        \lineUp{3,0}
        \lineUp{0,1}
        \lineUp{3,1}
        \lineUp{0,2}
        \lineUp{1,2}
        \lineUp{2,3}
        \lineUp{3,3}
        \lineUp{0,4}
        \lineUp{3,4}
        \lineUp{0,5}
        \lineUp{1,5}

        \zeroSmoothing{0,0}
        \zeroSmoothing{1,1}
        \zeroSmoothing{2,2}
        \oneSmoothing{0,3}
        \zeroSmoothing{1,4}
        \zeroSmoothing{2,5}
        \end{scope}

    \begin{scope}[shift={(14,0)}]
                \lineUp{2,0}
        \lineUp{3,0}
        \lineUp{0,1}
        \lineUp{3,1}
        \lineUp{0,2}
        \lineUp{3,3}
        
        \lineUp{3,4}

        \begin{scope}[color=blue]
        \lineRightBendUpMore{1,4}
        \lineUpBendLeftMore{1,3}
        \lineRightBendDownMore{1,2}
        \lineUp{1,2}
        \lineUpBendRightMore{2,2}
        \lineUp{2,3}
        \end{scope}
        \lineUp{0,4}
        \lineUp{0,5}
        \lineUp{1,5}

        \zeroSmoothing{0,0}
        \lineRightBendUpMore{1,1}
        \lineUpBendLeftMore{3,2}
        \lineUpBendRightMore{0,3}
        \lineRightBendDownMore{1,5}
        \zeroSmoothing{2,5}

    \end{scope}
\begin{scope}[shift={(21,0)}]
                \lineUp{2,0}
        \lineUp{3,0}
        \lineUp{0,1}
        \lineUp{3,1}
        \lineUp{0,2}
        \lineUp{1,2}
        \lineUp{3,3}
        
        \lineUp{3,4}

        \lineRightBendUpMore{1,4}
        \lineRightBendDownMore{1,2}
        \lineUp{1,2}
        \lineUpBendRightMore{2,2}
        \lineUp{2,3}
        \lineUp{0,4}
        \lineUp{0,5}
        \lineUp{1,5}

        \zeroSmoothing{0,3}
        \zeroSmoothing{0,0}
        \lineUpBendLeftMore{3,2}
        \lineRightBendDownMore{1,5}
        \lineRightBendUpMore{1,1}
        \zeroSmoothing{2,5}

    \end{scope}
    \begin{scope}[shift={(28,0)}]
                \lineUp{2,0}
        \lineUp{3,0}
        \lineUp{3,1}
        \lineUp{3,3}
        
        \lineUp{3,4}

        \lineRightBendUpMore{1,4}
        \lineUpBendLeftMore{2,1}
        \lineUp{1,2}
        \lineUpBendRightMore{2,2}
        \lineUp{2,3}
        \lineUp{0,4}
        \lineUp{0,5}
        \lineUp{1,5}

        \lineRightBendDownMore{0,4}
        \lineRightBendUpMore{0,0}
        \lineUpBendLeftMore{3,2}
        \lineRightBendDownMore{1,5}
        \zeroSmoothing{2,5}

        \begin{scope}[color=blue]
            \lineRightBendDownMore{0,1}
        \lineUpBendRightMore{1,1}
        \lineRightBendUpMore{0,3}
        \lineUp{1,2}
        \lineUp{0,2}
        \lineUp{0,1}
        \end{scope}

    \end{scope}
    \draw[->,dashed,magenta] (3.5,5) -- (6.5,5) node[midway, above, align=center] {\small isopair};
    \draw[<-,dashed,magenta] (3.5,1) -- (6.5,1) node[midway, below, align=center] {\small \textsc{MB}=False};
    \draw[->,dashed,magenta] (10.5,4.5) -- (13.5,4.5) node[midway, above, align=center] {\small isopair};
    
    \draw[<-,dashed,magenta] (10.5,1.5) -- (13.5,1.5) node[midway, below, align=center] {\small \textsc{MB}$=$True};
    
    \draw[->,dashed,magenta] (17.5,4) -- (20.5,4) node[midway, above, align=center] {\small isopair};
    \draw[<-,dashed,magenta] (17.5,2) -- (20.5,2) node[midway, below, align=center] {\small \textsc{MB}$=$False};
    \draw[->,dashed,magenta] (24.5,3.5) -- (27.5,3.5) node[midway, above, align=center] {\small isopair};
    \draw[<-,dashed,magenta] (24.5,2.5) -- (27.5,2.5) node[midway, below, align=center] {\small \textsc{MB}=True};

    \node at (1.5,-1) {$(111)^k00011\bzero$};
    \node at (8.5,-1) {$y_1$};
    \node at (15.5,-1) {$y_2$};
    \node at (22.5,-1) {$y_3$};
    \node at (29.5,-1) {$y_4.$};
        
    \end{tikzpicture}
\end{equation}
Let us elaborate on (\ref{Chain of pictures: unmatching of 000110}): The isopair arrows are executing lines $\ref{algoline: isopair 1}$ and \ref{algoline: isopair 2} whereas the \textsc{MB}=True/False  arrows execute lines \ref{algoline: MatchesBack 1} and \ref{algoline: MatchesBack 2} and state whether \textsc{MatchesBack} returns True or False. From this we can see that $(y_3\to y_4), (y_1\to y_2)\in \Mgr$ and thus  $(111)^k00011\bzero\in U_{n+1}$. Clearly, $(111)^k00011\bzero\in W_{n+1}$ as well and $W_{n+1}$ contains those  precisely those continuations which are in  $U_{n+1}$. In other words we have proven Equivalence \ref{Equivalence: induction for unmatched cells} in the case of $a\in  (g_1 \circ I_{\bA}  ) (\{e\}) $. 

Next assume that $a\in (g_6\circ I_{\bC} \circ g_5\circ  I_{\bB}\circ g_2\circ  I_{\bA}  )(\{e\})     \cup  (g_6 \circ I_{\bC} \circ g_3 \circ I_{\bA}  )(\{e\})$, so $a=\dots\gzero \gzero 01$ or $a=\dots \gzero \gzero 0101\bzero$ where each $\gzero$ is individually a placeholder for either $0$ or $\bzero$. Let us first investigate $a=\dots\gzero \gzero 01$ first: The continuations $11$ and $011$ get matched via merge-type isomorphisms whereas $00$ and $1\bzero$ get matched with each other:
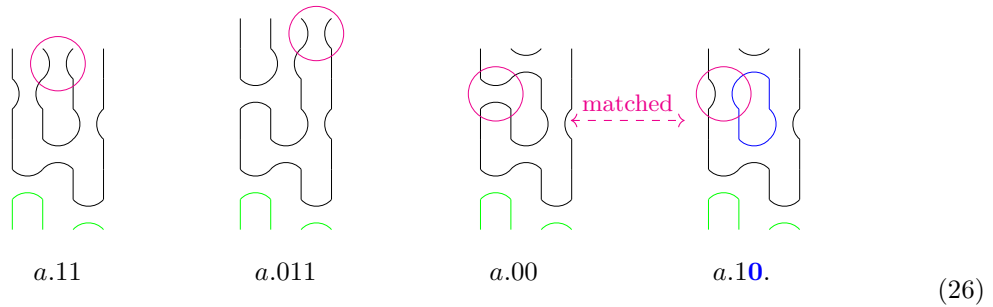
\begin{equation}\label{Picture: classification of arrows with high tC part 1}
    \begin{tikzpicture}
        \node at (0,0) {

\begin{tikzpicture}
\node at (0, 0) [anchor=south west, scale=0.4] { 
            \input{WUclassification/000111}
        };  
\node at (3, 0) [anchor=south west, scale=0.4] { 
            \input{WUclassification/0001011} 
        };  
\node at (6, 0) [anchor=south west, scale=0.4] { 
            \input{WUclassification/000100}
        };  
\node at (9, 0) [anchor=south west, scale=0.4] { 
            \input{WUclassification/00011x} 
        };

    \draw[<->,dashed,magenta] (7.4,1.5) -- (8.9,1.5) node[midway, above, align=center] {\small matched};

\node at (0.65,-0.5) {$a.11$};

\node at (3.65,-0.5) {$a.011$};

\node at (6.65,-0.5) {$a.00$};

\node at (9.65,-0.5) {$a.1\bzero.$};

\end{tikzpicture}

        };
    \end{tikzpicture}
\end{equation}
Again, running Algorithm \ref{Algorithm: matching a single cell in Mgr} we can see that the cell $a.01\bzero$ remains unmatched:
    $$
    
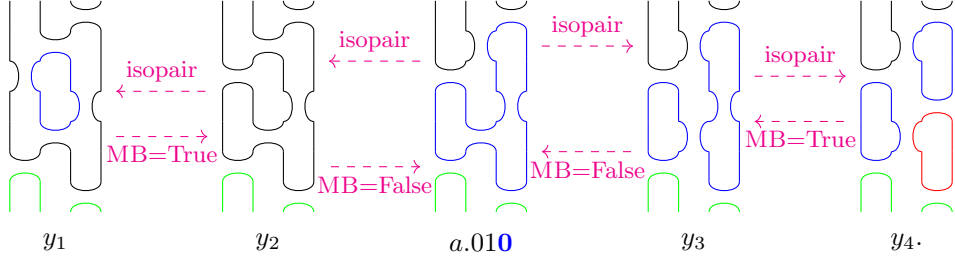
\begin{tikzpicture}[scale=0.4]
        
        \begin{scope}[green]
            
        \lineUp{0,0}
        \lineUp{1,0}
        \lineRightBendUpMore{2,0}
        
        \lineRightBendUpMore{0,1}

        \end{scope}

        \begin{scope}[color=blue]
            \lineRightBendDownMore{2,1}
            \lineRightBendDownMore{0,2}
            \lineUp{2,1}
            \lineUp{3,1}
            \lineUp{0,2}
            \lineRightBendUpMore{1,2}
            \lineRightBendDownMore{1,3}
            \lineUp{3,2}
            \lineUp{0,3}
            \lineUp{1,3}
            \lineUpBendRightMore{2,3}
            \lineUpBendLeftMore{3,3}
            \lineRightBendUpMore{0,4}
            \lineUp{2,4}
            \lineUp{3,4}
            \lineUp{3,5}
            \lineUpBendLeftMore{2,5}
            
            \lineRightBendUpMore{2,6}

        \end{scope}
        
            \lineRightBendDownMore{0,5}
            \lineUp{0,5}
            \lineUpBendRightMore{1,5}
            \lineUp{0,6}
            \lineUp{1,6}
        \lineRightBendDownMore{2,7}

        \begin{scope}[shift={(-7,0)}]
                \begin{scope}[green]
            
        \lineUp{0,0}
        \lineUp{1,0}
        \lineRightBendUpMore{2,0}
        
        \lineRightBendUpMore{0,1}

        \end{scope}

            \lineRightBendDownMore{2,1}
            \lineRightBendDownMore{0,2}
            \lineUp{2,1}
            \lineUp{3,1}
            \lineUp{0,2}
            \lineRightBendUpMore{1,2}
            \lineRightBendDownMore{1,3}
            \lineUp{3,2}
            \lineUp{0,3}
            \lineUp{1,3}
            \lineUpBendRightMore{2,3}
            \lineUpBendLeftMore{3,3}
            \lineRightBendUpMore{0,4}
            \lineUp{2,4}
            \lineUp{3,4}
            \lineUp{3,5}

            \lineRightBendUpMore{1,5}
            
            \lineRightBendUpMore{2,6}

            \lineRightBendDownMore{0,5}
            \lineUp{0,5}
            \lineRightBendDownMore{1,6}
            \lineUp{0,6}
            \lineUp{1,6}
        \lineRightBendDownMore{2,7}

        \end{scope}

        \begin{scope}[shift={(-14,0)}]
               \begin{scope}[green]
            
        \lineUp{0,0}
        \lineUp{1,0}
        \lineRightBendUpMore{2,0}
        
        \lineRightBendUpMore{0,1}

        \end{scope}

        \begin{scope}[color=blue]

            \lineUp{1,3}
            \lineUpBendRightMore{2,3}
            
            \lineRightBendDownMore{1,3}
            \lineUp{2,4}

            \lineRightBendUpMore{1,5}
            \lineUpBendLeftMore{1,4}
        \end{scope}
            \lineRightBendDownMore{2,1}
            \lineRightBendDownMore{0,2}
            \lineUp{2,1}
            \lineUp{3,1}
            \lineUp{0,2}
            
            \lineRightBendUpMore{1,2}

            \lineUp{3,4}
        
          \lineUp{3,5}
          
            \lineUpBendLeftMore{3,3}
                        \lineUp{3,2}
            \lineUp{0,3}
            \lineUpBendRightMore{0,4}

            \lineRightBendUpMore{2,6}

            \lineUp{0,5}
            \lineRightBendDownMore{1,6}
            \lineUp{0,6}
            \lineUp{1,6}
        \lineRightBendDownMore{2,7}

        \end{scope}

        \begin{scope}[shift={(7,0)}]
                 \begin{scope}[green]
            
        \lineUp{0,0}
        \lineUp{1,0}
        \lineRightBendUpMore{2,0}
        
        \lineRightBendUpMore{0,1}

        \end{scope}

        \begin{scope}[color=blue]
            \lineRightBendDownMore{2,1}
            \lineRightBendDownMore{0,2}
            \lineUp{2,1}
            \lineUp{3,1}
            \lineUp{0,2}
            \lineUpBendRightMore{1,2}
            \lineUpBendLeftMore{2,2}
            \lineUp{3,2}
            \lineUp{0,3}
            \lineUp{1,3}
            \lineUpBendRightMore{2,3}
            \lineUpBendLeftMore{3,3}
            \lineRightBendUpMore{0,4}
            \lineUp{2,4}
            \lineUp{3,4}
            \lineUp{3,5}
            \lineUpBendLeftMore{2,5}
            
            \lineRightBendUpMore{2,6}

        \end{scope}
        
            \lineRightBendDownMore{0,5}
            \lineUp{0,5}
            \lineUpBendRightMore{1,5}
            \lineUp{0,6}
            \lineUp{1,6}
        \lineRightBendDownMore{2,7}

        \end{scope}

    \begin{scope}[shift={(14,0)}]
             \begin{scope}[green]
            
        \lineUp{0,0}
        \lineUp{1,0}
        \lineRightBendUpMore{2,0}
        
        \lineRightBendUpMore{0,1}

        \end{scope}

        \begin{scope}[color=blue]
            
            \lineRightBendDownMore{0,2}
            \lineUp{0,2}
            \lineUp{0,3}
            \lineUp{1,3}
            \lineRightBendUpMore{0,4}
            \lineUp{2,4}
            \lineUp{3,4}
            \lineRightBendDownMore{2,4}
            
            \lineUp{3,5}
            \lineUpBendLeftMore{2,5}
            
            \lineUpBendRightMore{1,2}
            
            \lineRightBendUpMore{2,6}

        \end{scope}
        \begin{scope}[red]
            
            \lineRightBendDownMore{2,1}
            \lineUp{2,1}
            \lineUp{3,1}
            \lineRightBendUpMore{2,3}
            
            \lineUpBendLeftMore{2,2}
            \lineUp{3,2}
        \end{scope}

            \lineRightBendDownMore{0,5}
            \lineUp{0,5}
            \lineUpBendRightMore{1,5}
            \lineUp{0,6}
            \lineUp{1,6}
        \lineRightBendDownMore{2,7}
        \end{scope}

    \draw[dashed,magenta,->] (-0.5,5) -- (-3.5,5) node[midway, above, align=center, color=magenta] {\small isopair}; 
    \draw[dashed,magenta,<-] (-10.5,4) -- (-7.5,4) node[midway, above, align=center, color=magenta] {\small isopair};
    \draw[dashed,magenta,->] (-10.5,2.5) -- (-7.5,2.5) node[midway,below, align=center, color=magenta] {\small \textsc{MB}$=$True};

    \draw[dashed,magenta,<-] (-0.5,1.5) -- (-3.5,1.5) node[midway,below, align=center, color=magenta] {\small \textsc{MB}$=$False};
    
    \draw[dashed,magenta,->] (3.5,5.5) -- (6.5,5.5) node[midway, above, align=center, color=magenta] {\small isopair};
    
    \draw[dashed,magenta,<-] (3.5,2) -- (6.5,2) node[midway,below, align=center, color=magenta] {\small \textsc{MB}$=$False};
    
    \draw[dashed,magenta,->] (10.5,4.5) -- (13.5,4.5) node[midway, above, align=center, color=magenta] {\small isopair};
    \draw[dashed,magenta,<-] (10.5,3) -- (13.5,3) node[midway,below, align=center, color=magenta] {\small \textsc{MB}$=$True};

    \node at (1.5,-1) {$a.01\bzero$};
    \node at (8.5,-1) {$y_3$};
    \node at (15.5,-1) {$y_4.$};
    \node at (-5.5,-1) {$y_2$};
    \node at (-12.5,-1) {$y_1$};

    \end{tikzpicture}
    $$
as the arrows $(y_2\to y_1),(y_4\to y_3)$ are contained in $\Mgr$.

Next, case $a=\dots \gzero \gzero 0101\bzero$: Continuations $1$ and $\bzero11$ get similarly matched with merge-type arrows whereas $\bzero1\bzero$ get pairs up with $\bzero00$:
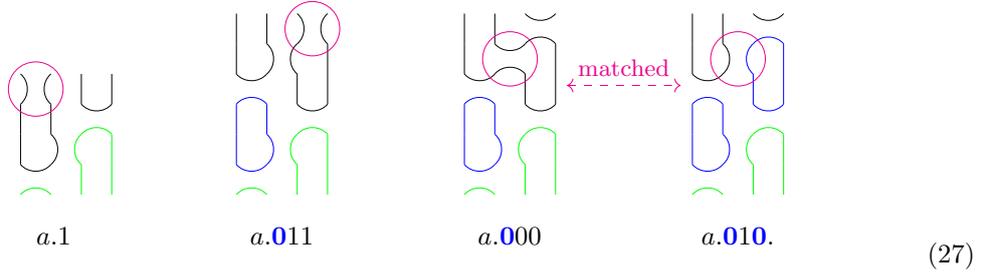
\begin{equation} \label{Picture: classification of arrows with high tC part 2}
    \begin{tikzpicture}
        \node at (0,0) {

\begin{tikzpicture}
\node at (0, 0) [anchor=south west, scale=0.4] { 
            \input{WUclassification/abca0101}
        };  
\node at (3, 0) [anchor=south west, scale=0.4] { 
            \input{WUclassification/abcabc010x11} 
        };  
\node at (6, 0) [anchor=south west, scale=0.4] { 
            \input{WUclassification/abcabc010x00}
        };  
\node at (9, 0) [anchor=south west, scale=0.4] { 
            \input{WUclassification/abcabc010x1x} 
        };

    \draw[<->,dashed,magenta] (7.4,1.5) -- (8.9,1.5) node[midway, above, align=center] {\small matched};

\node at (0.65,-0.5) {$a.1$};

\node at (3.65,-0.5) {$a.\bzero11$};

\node at (6.65,-0.5) {$a.\bzero00$};

\node at (9.65,-0.5) {$a.\bzero 1 \bzero$.};

\end{tikzpicture}

        };
    \end{tikzpicture}
\end{equation}
The continuation $\bzero01$ remains unmatched
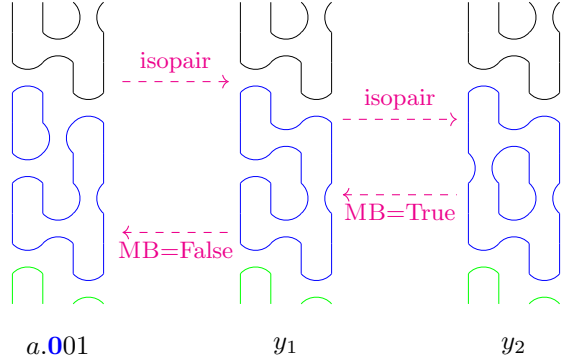
\begin{center}
\input{WUclassification/cabcabcabc}
\end{center}
since $(y_1 \to y_2) \in \Mgr$. By comparing these unmatched continuations with the definition of $W$, we can see that Equivalence \ref{Equivalence: induction for unmatched cells} holds also in the case  $a\in (g_6\circ I_{\bC} \circ g_5\circ  I_{\bB}\circ g_2\circ  I_{\bA}  )(\{e\})     \cup  (g_6 \circ I_{\bC} \circ g_3 \circ I_{\bA}  )(\{e\})$. We omit the details of the last case: showing that Equivalence \ref{Equivalence: induction for unmatched cells} holds for $a\in (g_4 \circ I_{\bB} \circ g_2 \circ I_{\bA}  ) (\{e\})$.   
\end{proof}
To additionally ensure that no mistakes were made in the vast amount of cases of the previous proof, $U_n=W_n$ is verified with \texttt{braidalgo.py} for $n=0,\dots , 50$.

In order to manage the recursions of $W$, we define the formal words of length 12
$$
\bA=1^{12},\qquad \bB=101011\bzero0011\bzero, \qquad \bC= (01\bzero \bzero 01)^2
$$
and for every $\ast\in \{ \bA,\bB, \bC\}$ and $r\in \mathbb Z_{\geq 0}$ set 
$$
\mathcal W_{\ast,r}=\{a\in W \mid a \text{ contains at least } r \text{ disjoint instances of } \ast \text{ as a subword}\}.
$$
We also define bijections $\gamma_{\ast}\colon \mathcal W_{\ast,1}\to \mathcal W_{\ast,2} $, where the value $\gamma_{\ast}(w)$ is obtained by replacing the first occurrence of subword $\ast$ in the word $w$ by the subword $\ast^2$. Let us establish two lemmas which allow us to connect the formal recursions $\gamma_\ast$ on $W$ with quantum and homological gradings in $C_n$.
\begin{lemma}\label{Lemma: recursions of W from 24 on}
    For $n\geq 24$ it holds that 
    $$U_n\subset \bigcup_{\ast \in\{\bA,\bB, \bC\}} \mathcal W_{\ast,2}=\bigcup_{\ast \in\{\bA,\bB, \bC\}} \operatorname{im}\gamma_{\ast}.$$
\end{lemma}
\begin{proof}
    By staring at the definition of $W$ one can confirm that the longest word in $W$ which does not contain $\bA^2,\bB^2$ or $\bC^2$ has length $3\cdot 23$ and longer words must contain one of these patterns.
\end{proof}

\begin{lemma}\label{Lemma: t value implies W}
    Suppose $w\in U$ and $\ast, \ast' \in\{\bA,\bB, \bC\}$. If $w\in \mathcal W_{\ast,n}$ we have
    \begin{equation}\label{Equation: t gamma delta}
    t_{\ast'}(w) + \delta_{\ast,\ast'}=t_{\ast'}(\gamma_{\ast}(w)) 
    \end{equation}
    where $\delta$ is the Kronecker delta ($\delta_{\ast,\ast'}=1$, if $\ast=\ast'$ and $\delta_{\ast,\ast'}=0$ otherwise).  The inequality $t_{\ast}(w)\geq m$ implies $w\in \mathcal W_{\ast,m}$. 
\end{lemma}
\begin{proof}
    The first claim amounts to a simple degree computation and  for the second claim 
    proceed with an induction on $n$. For words $w\in U_n$ with $n\leq 23$ the implication is verified with \texttt{basecases.py}. Suppose now that it holds for all $w\in U_{n-4}$ with $n\geq 24$ and let $w\in U_n$ with $t_{\ast}(w)\geq m$. By Lemma \ref{Lemma: recursions of W from 24 on} we have $w= \gamma_{\ast'}(v)$ for some $v$ and $\ast'$. Applying Equation \ref{Equation: t gamma delta} to this yields  
    $t_{\ast}(v)\geq m-\delta_{\ast,\ast'}$. Hence by the induction assumption $v\in \mathcal W_{\ast,k}$ where $k=m-\delta_{\ast, \ast'}$. Thus $w=\gamma_{\ast'}(v)\in \mathcal W_{\ast,m}$.
\end{proof}
We can  now prove Lemma \ref{Lemma: minimal tA,tC} with an analogous induction.
\begin{proof}[Proof of Lemma \ref{Lemma: minimal tA,tC}] 
With \texttt{basecases.py} we verify that $t_{\bA}(w)$, $t_{\bC}(w)\geq -\frac{3}{2}$ for all $w\in U_n$ with $n\leq 23$. Now, let $w\in U_n$ and assume that $t_{\bA}(v)$, $t_{\bC}(v) \geq -\frac{3}{2}$ for all $v\in U_{n-4}$. By Lemma \ref{Lemma: recursions of W from 24 on} and Equation \ref{Equation: t gamma delta} we have
$$
t_{\bA}(w)=t_{\bA}(\gamma_{\ast}(v))\geq t_{\bA}(v) + \delta_{\ast, \bA}\geq -\frac{3}{2}
$$
using some $v\in U_{n-4}$ and $\ast\in \{\bA,\bB, \bC\}$. Similarly $t_{\bC}(w)\geq -\frac{3}{2}$.
\end{proof}

We will now define the functions $\mathfrak a_n$, $\mathfrak b_n$ and $\mathfrak c_n$ as restrictions and corestrictions:  
    \begin{alignat*}{3}
        \mathfrak{a}_n\colon& \{ w\in U_n \mid t_{\bA}(w)\geq 1 \} &&\to \{ w\in U_{n+4} \mid t_{\bA}(w)\geq 2 \}, \  &&\mathfrak a_n (w)=\gamma_{\bA}(w) \\ 
        \mathfrak{b}_n\colon& \{ w\in U_n \mid t_{\bB}(w)\geq 1 \} &&\to \{ w\in U_{n+4} \mid t_{\bB}(w)\geq 2 \}, \ &&\mathfrak b_n (w)=\gamma_{\bB}(w) \\ 
        \mathfrak{c}_n\colon& \{ w\in U_n \mid t_{\bC}(w)\geq 1 \} &&\to \{ w\in U_{n+4} \mid t_{\bC}(w)\geq 2 \}, \ &&\mathfrak c_n (w)=\gamma_{\bC}(w).  
    \end{alignat*}
Thus proving Proposition \ref{Proposition: abc bijections on cells} amounts to verifying that $\mathfrak a_n$, $\mathfrak b_n$ and $\mathfrak c_n$ are well-defined bijections which respect the  connectivities and shift the gradings in the assigned way.
\begin{proof}[Proof of Proposition \ref{Proposition: abc bijections on cells}]
    Let $\ast \in \{ \bA, \bB, \bC \}$ and $w\in U_n$ with $t_\ast (w)\geq 1$. Then by Lemma \ref{Lemma: t value implies W} $w\in \mathcal W_{\ast,1}$ and we can apply $\gamma_\ast$ to $w$. The Equation \ref{Equation: t gamma delta} yields $t_\ast(\gamma_\ast(w))\geq 2$ and hence $\mathfrak a_n$, $\mathfrak b_n$ and $\mathfrak c_n$ are well defined. The proof of surjectivity works similarly through Lemma \ref{Lemma: t value implies W} and Equation \ref{Equation: t gamma delta}. The functions $\mathfrak a_n$, $\mathfrak b_n$ and $\mathfrak c_n$ are also injections, since $\gamma_{\bA}$, $\gamma_{\bB}$ and $\gamma_{\bC}$ are.

    A simple grading computation shows that for all $a\in \mathcal W_{a,1}$, $b\in \mathcal W_{b,1}$, $c\in \mathcal W_{c,1}$
    \begin{align*}
        (h(a),q(a))&=(h(\gamma_{\bA}(a)),q(\gamma_{\bA}(a)))+(0,12) \\
        (h(b),q(b))&=(h(\gamma_{\bB}(b)),q(\gamma_{\bB}(b)))+(6,20) \\
        (h(c),q(c))&=(h(\gamma_{\bC}(c)),q(\gamma_{\bC}(c)))+(8,24). 
    \end{align*}
    It is also easy to see from Figure \ref{Figure: periodic building blocks} that connectivities of the underlying Temperley-Lieb diagrams remain unchanged.
\end{proof}


\section{Proof of Proposition \ref{Proposition: ac matrix element equations}}\label{Section: proof of morphism correspondences}

In this section we prove Proposition \ref{Proposition: ac matrix element equations} separately for cases $\mathfrak a_n$ and $\mathfrak c_n$. 

\subsection{Case \texorpdfstring{$\mathfrak a_n$}{an}}

The proof for $\mathfrak a_n$ is more or less straightforward, as the equalities of matrix elements are induced by a correspondences of paths which in turn are given by a graph embedding. The proof comes as a consequence of the following slightly more general lemma.  
\begin{lemma}\label{Lemma: adding 1 to the front}
    Let $B$ be a negative braid word on $n$ strands, $1 \leq i \leq n-1$ and assume $\Mgr$ is a Morse matching on  $G(\Psi\llbracket \sigma_i B \rrbracket)$.  Then $\Mgr$ is a Morse matching on $G(\llbracket B \rrbracket)$ and the injective map 
    $$
    \varphi \colon  \operatorname{cells}(\Psi \llbracket B \rrbracket) \to \operatorname{cells}(\Psi \llbracket \sigma_i B \rrbracket\{1\}), \ \varphi(w)=(1.w)\{1\}
    $$
    sends any matrix element of $\Mgr \Psi \llbracket B \rrbracket $ to its additive inverse: $d_{b,a}=-\partial_{\varphi(b),\varphi(a)}$.
\end{lemma}
\begin{proof}
    It is easy to see that $\varphi$ preserves the degrees and the connectivity of diagrams. Additionally, $\varphi$ defines a graph homomorphism 
    $$
    \varphi \colon G(\Psi\llbracket B \rrbracket) \to G(\Psi\llbracket \sigma_i B \rrbracket \{1\}) 
    $$
    which is a graph isomorphism into the induced subgraph of its image. We want to show that $\varphi$ also gives a graph isomorphism with the arrows of $\Mgr$ reversed on both graphs. In other words, we want to show that 
    $$
    \rho \colon G(\Psi\llbracket B \rrbracket, \Mgr) \to G(\Psi\llbracket \sigma_i B \rrbracket\{1\}, \Mgr) , \  \rho(x) =\varphi(x) 
    $$
    is a graph isomorphism into the induced subgraph of its image. This amounts to showing that an edge $a\to b$ is contained in the greedy matching of $\Psi \llbracket B \rrbracket$ if and only if the edge $\varphi (a) \to \varphi (b)$ is contained in the greedy matching of $\Psi \llbracket \sigma_i B \rrbracket$. Morally, this follows from the fact that $1$-smoothing of the first crossing does not change the geometry of the pictures and from Lemma \ref{lemma: red zeros matched}. 

    A rigorous argument proving this falls back on Algorithm \ref{Algorithm: matching a single cell in Mgr}. We write $\operatorname{isopair}_B$, $\operatorname{isopair}_{\sigma_i B}$, $\textsc{MatchesBack}_B$ and $\textsc{MatchesBack}_{\sigma_i B}$ to distinguish the implicit underlying tangle diagram in Algorithm $\ref{Algorithm: matching a single cell in Mgr}$.  By additionally writing $\varphi(\star)=\star$, we get that for all $j,k\geq 0$, $w\in \operatorname{cells}(\Psi \llbracket B \rrbracket)$ 
    $$
    \varphi(\operatorname{isopair}_B(w,j,k))=\operatorname{isopair}_{\sigma_i B}(\varphi(w),j+1,k+1).
    $$
    By doing an induction on the sequence
    $$
    (j,k)=(1,1),(2,2),(1,2),(3,3),(2,3),(1,3),(4,4),\dots
    $$
    and utilizing Lemma \ref{lemma: red zeros matched} one can show that for all $a,b,j,k$ it holds that
    $$
    \varphi( \textsc{MatchesBack}_B(\operatorname{isopair}_B(a,j,k),b))=
    \textsc{MatchesBack}_{\sigma_iB}(\operatorname{isopair}_{\sigma_i B}(\varphi(a),j+1,k+1),\varphi(b)).
    $$
    Thus $(a\to b)\in \Mgr \iff \varphi(a) \to \varphi(b)$ and $\rho$ is a graph isomorphism to the induced subgraph of its image. It follows that, the acyclicity of $G(\Psi\llbracket \sigma_i B \rrbracket)$ can be pulled back to $G(\Psi\llbracket B \rrbracket)$ making $\Mgr$ a Morse matching on the latter graph.
    
    The graph embedding $\rho$ gives rise to injective maps 
    $$
    \rho_{y,x}\colon \{\text{paths from } x \text{ to } y \} \to \{\text{paths from } \rho(x) \text{ to } \rho(y) \} 
    $$
    by
    $$
    \rho_{y,x} (x\to a_1\to \dots \to a_n \to y)= (\rho(x)\to \rho(a_1)\to \dots \to \rho(a_n) \to \rho(y)).
    $$
    Moreover, the maps $\rho_{y,x}$ are surjective, since Lemma \ref{lemma: red zeros matched} obstructs the existence of an arrow $(a\to b)\in G(\Psi \llbracket\sigma_i B \rrbracket , \Mgr)$ with $a\in \operatorname{im}\rho$ and $b\not \in \operatorname{im} \rho $.

    Recall that the differentials $d$ of $\Mgr \Psi \llbracket B \rrbracket $ and $\partial$ of $\Mgr \Psi \llbracket \sigma_i B \rrbracket$ are defined as sums:
    $$
    d_{b,a}=\sum_{p\in \{ \text{paths}\colon a\to b\}} R_B (p), \qquad \qquad \partial_{\varphi(b),\varphi(a)}=\sum_{p'\in \{ \text{paths}\colon \rho(a)\to \rho(b)\}} R_{\sigma_iB} (p) 
    $$
    where the paths are evaluated by the remembering functors:
    \begin{align*}
    R_B\colon G(\Psi \llbracket B \rrbracket, \Mgr ) &\to \operatorname{Mat} (\Cob(8)) \\
    R_{\sigma_i B}\colon G(\Psi \llbracket \sigma_i B \rrbracket, \Mgr ) &\to \operatorname{Mat} (\Cob(8)).        
    \end{align*}
    Under the two remembering functors, a single edge $x\to y$ gets mapped to equal cobordisms but with opposite signs: 
    $$
    R_B(x\to y)=-R_{\sigma_i B} (\varphi (x) \to \varphi (y)).
    $$
    Every path contributing to $d_{b,a}$ and to $\partial_{\varphi(b),\varphi(a)}$ zig-zags between two homological layers and thus has an odd number of edges. Hence the sign swaps for every path and we can use bijectivity of $\rho_{y,x}$ to get $d_{b,a}=-\partial_{\varphi(b), \varphi(a)}$.
\end{proof}

\begin{proof}[Proof of Proposition \ref{Proposition: ac matrix element equations}, case $\mathfrak{a}_n$.]
Proposition \ref{Proposition: Mgr Morse on T4} shows that $\Mgr$ is a Morse matching on the graph $G(\Psi \llbracket (\sigma_1 \sigma_2 \sigma_3)^n \rrbracket)$, for all $n$. From Lemma \ref{Lemma: adding 1 to the front}, it follows that $\Mgr$ is a Morse matching on graphs $G(\Psi \llbracket \sigma_3(\sigma_1 \sigma_2 \sigma_3)^n \rrbracket)$ and $G(\Psi \llbracket \sigma_2\sigma_3(\sigma_1 \sigma_2 \sigma_3)^n \rrbracket)$ as well. Now, let $w,u \in U_n$ with $t_{\bA}(w), t_{\bA}(u) \geq 1$ and assume that $\partial_{u,w}$ is a matrix element of $C_n=\Mgr (\Psi\llbracket (\sigma_1\sigma_2\sigma_3)^n \rrbracket) $. Applying $\varphi$ from Lemma \ref{Lemma: adding 1 to the front} twelve times coincides with the function $\mathfrak{a}_n$ from Proposition \ref{Proposition: abc bijections on cells}, that is,
$$
(\varphi \circ \varphi \circ \dots \circ \varphi)(x)=\mathfrak{a}_n(x)
$$
for any unmatched cell $x$. Thus
$$
\partial_{u,w}=(-1)^{12} \partial_{(\varphi \circ \dots \circ \varphi)(u),(\varphi \circ \dots \circ \varphi)(w)}=\partial_{\mathfrak a_n(u),\mathfrak a_n(w)}.
$$    
\end{proof}

\subsection{Case \texorpdfstring{$\mathfrak c_n$}{cn}}
Unfortunately, in the proof of Proposition \ref{Proposition: ac matrix element equations}, case $\mathfrak c_n$, we do not manage to generate a simple graph isomorphism which would automatically yield a correspondence of paths. Nevertheless, the basic idea remains the same: a correspondence of paths $\Phi_n$ will give out an equality of matrix elements. This time, constructing the map $\Phi_n$ will require some amount of manual labor and showing that $\Phi_n$ works in the case of $\mathfrak c_n$ is also more demanding.

With directed graphs we have respected to the usual convention that paths follow the orientation of edges. To contrast these usual \term{orientation preserving paths} we say that a sequence of vertices $a_1,\dots,a_k$ is \term{an orientation neglecting path} if for every $i$ there exists either an edge $a_i \to a_{i+1}$ or an edge $a_{i+1}\to a_i$ in the underlying directed graph. Given $a,b\in U_n$ we define $\operatorname{ZZ}(a,b)$ to be the set of orientation preserving zig-zag paths from $a$ to $b$ and $Z_n=\{p\in \operatorname{ZZ}(a,b) \mid a,b\in U_n, \ t_{\bC}(a)\geq 6, \ t_{\bC}(b)\geq 1\} $.

The proof of Proposition \ref{Proposition: ac matrix element equations}, case $\mathfrak{c}_n$ is broken down into the following three tasks:
\begin{enumerate}
    \item Define a set $A_n\subset Z_n$. For each orientation preserving path $(p\colon a\to b)\in A_n$ define an orientation neglecting path $\Phi_n(p)\colon \mathfrak c_n(a)\to \mathfrak c_n(b)$.
    \item Show that $\Phi_n(p)$ is in fact always orientation preserving and hence $\Phi_n \colon A_n \to Z_{n+4}$ is a well-defined and injective map. Prove that $R(p)=R(\Phi_n (p))$ as morphisms in $\operatorname{Mat}(\Cob (8))$ for all $p\in A_n$.
    \item Prove that for all $a,b\in U_n$ with $t_{\bC}(a)\geq 6$, $t_{\bC}(b)\geq 1$ we have
    \begin{align}
    \sum_{p\in \operatorname{ZZ}(a,b)} R(p)&= \sum_{p\in \operatorname{ZZ}(a,b)\cap A_n} R(p) \\
    \sum_{p'\in \operatorname{ZZ}(\mathfrak{c}_n(a), \mathfrak c_n (b))}R(p') &=\sum_{p'\in \operatorname{ZZ}(\mathfrak{c}_n(a), \mathfrak c_n (b))\cap  \operatorname{im} \Phi_n}R(p'). 
    \end{align}
\end{enumerate}
Once the three tasks are completed, one can calculate 
\begin{multline*}
\partial_{b,a}=\sum_{p\in \operatorname{ZZ}(a,b) } R(p) = \sum_{p\in \operatorname{ZZ}(a,b) \cap A_n} R(p) = \\ 
=\sum_{p'\in \operatorname{ZZ}(\mathfrak{c}_n(a), \mathfrak c_n (b))\cap  \operatorname{im} \Phi_n}R(p')=\sum_{p'\in \operatorname{ZZ}(\mathfrak{c}_n(a), \mathfrak c_n (b))}R(p') =\partial_{\mathfrak{c}_n(b), \mathfrak c_n (a)}
\end{multline*}
for any $a,b\in U_n$ with $t_{\bC}(a)\geq 6$ and $t_{\bC}(b)\geq 1$ concluding the proof.

\subsubsection{Task 1: setting up  \texorpdfstring{$\Phi_n$}{Phin}}
Informally, the definition of $\Phi_n$ will consist of \term{local subpaths} $B_m(p)$ which are transformed into different \term{global environments} $s,r\in \mathcal{R}$ with \term{surgery maps} $\Lambda_{s,t}$ and then appended on both sides with \term{local moves} such as $\lozenge_m$ and $\triangledown_m$. Let us dive into the specifics.

Let $p=(p_1\to \dots \to p_{2k})\in Z_n$ be a zig-zag path. For $m\leq n$ and $j< k$ we define the index 
$$
\beta_m(p,2j+1)=\min \{ l>j \mid l=k \text{ or } i(p_{2l}\to p_{2l+1}) >3m \}.
$$
We further define $B_m(p)$ to be the first subpath of $p$ satisfying the following conditions:
\begin{enumerate}
    \item $B_m(p)=(p_{2a+1}\to \dots \to p_{2b})$
    \item $a=0$ or $i(p_{2a} \to p_{2a+1})>3m$
    \item $b=k$ or $i(p_{2b} \to p_{2b+1})>3m$
    \item $i(p_{2j} \to p_{2j+1})\leq 3m$ for all $a<j<b$
    \item $i(p_{2a+1}\to p_{2a+2})\leq 3m$.
\end{enumerate}
Such subpath need not to exist which is why one needs to be slightly careful when using the notation. 

The set of global environments $\mathcal R$ consists of upper-halves dissected cells of $\Psi \llbracket (\sigma_1\sigma_2 \sigma_3)^n \rrbracket$ for varying $n$, that is,
$$
\mathcal{R} = \{ r\in \mathcal S \mid vr\in \operatorname{cells} \Psi \llbracket (\sigma_1\sigma_2 \sigma_3)^n \rrbracket , \ n\in \mathbb Z_{\geq 0}, \ |r| \in 3\mathbb Z \}.  
$$
Recall that a cell of $\Psi \llbracket (\sigma_1\sigma_2 \sigma_3)^n \rrbracket$ contains a choice of smoothing, 0 or 1, for each crossing as well as a choice of color, red or blue, for each loop. This data is also present in every $r\in \mathcal R$ but additionally there is extra information on the colors (red, blue or black) of the arcs whose boundary is at the bottom of $r$. This extra data determines whether the arcs in $r$ were boundary arcs in $vr$ (color black) or whether the arcs in $r$ were part of complete loops in $vr$ with either color red or blue.  

\begin{figure}
    \centering
    \input{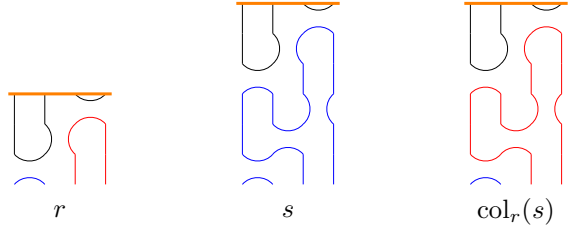}
    \caption{Two global environments $r,s\in \mathcal R$. Since $r\sim s$ we can apply $\operatorname{col}_r$ to $s$, and the stronger relation applies: $s \approx \operatorname{col}_r(s)$.}
    \label{Figure: Equivalence relations on global enviroments}
\end{figure}

For words $r,s\in \mathcal R$ we denote $r\sim s$, if the connectivity of the underlying Temperley-Lieb diagrams of $r$ and $s$ are the same. This defines an equivalence relation on $\mathcal R$ and we denote the equivalence classes by $[r]_{\sim}$. Given $r\in  \mathcal R$ we define a function $\operatorname{col}_r\colon [r]_{\sim } \to  [r]_{\sim }$  where $\operatorname{col}_r(s)$ is the element in $\mathcal{R}$ which is obtained from $s$ by changing the colors of the bottom boundary arcs of $s$ to match those of $r$, see Figure \ref{Figure: Equivalence relations on global enviroments}. A stronger equivalence relation $\approx$ on $\mathcal{R}$ is also needed: $r\approx s$, if $r\sim s$ and $\operatorname{col}_r(s)=r$, that is, $r$ and $s$ almost exactly the same; they may differ only by the coloring of their bottom components. 

Given a global environment $s\in \mathcal{R}$ and $n\in \mathbb Z_{\geq 0} $ we define
$$
V_{s,n}=\{vr \in \operatorname{cells} \Psi \llbracket (\sigma_1\sigma_2 \sigma_3)^n \rrbracket \mid r\approx s, \ t_{\bC}(v)\geq 1  \}
$$
and for any $t\in [s]_{\sim}$ we have the surgery map 
$$
\Lambda_{s,t,n} \colon V_{s,n} \to V_{t,m}, \ \Lambda_{s,t,n}(vr)=v \operatorname{col}_r(t).
$$
where $m=n+\frac{1}{3}(|t|-|w|)$. Suppose that $a,b\in V_{s,n}$ and there is an edge  $(a\to b)\in G(\Psi \llbracket (\sigma_1\sigma_2 \sigma_3)^n \rrbracket, \Mgr )$. It follows that there will also exist an edge $\Lambda_{s,t,n} (a) \to \Lambda_{s,t,n} (b)$ or an edge $\Lambda_{s,t,n} (b) \to \Lambda_{s,t,n} (a)$ in the codomain $G(\Psi \llbracket (\sigma_1\sigma_2 \sigma_3)^m \rrbracket, \Mgr )$. The orientation of $a\to b$ might change under $\Lambda_{s,t,m}$, as a priori it cannot be guaranteed that both or neither of the edges are contained in their respective greedy matchings. Hence $\Lambda_{s,t,m}$ will map an orientation preserving path $p$ to an orientation neglecting path $\Lambda_{s,t,m}(p)$, and  it is not immediately clear whether $\Lambda_{s,t,m}(p)$ is an honest orientation preserving path in the directed graph $G(\Psi \llbracket (\sigma_1\sigma_2 \sigma_3)^m \rrbracket, \Mgr )$.

\begin{figure}
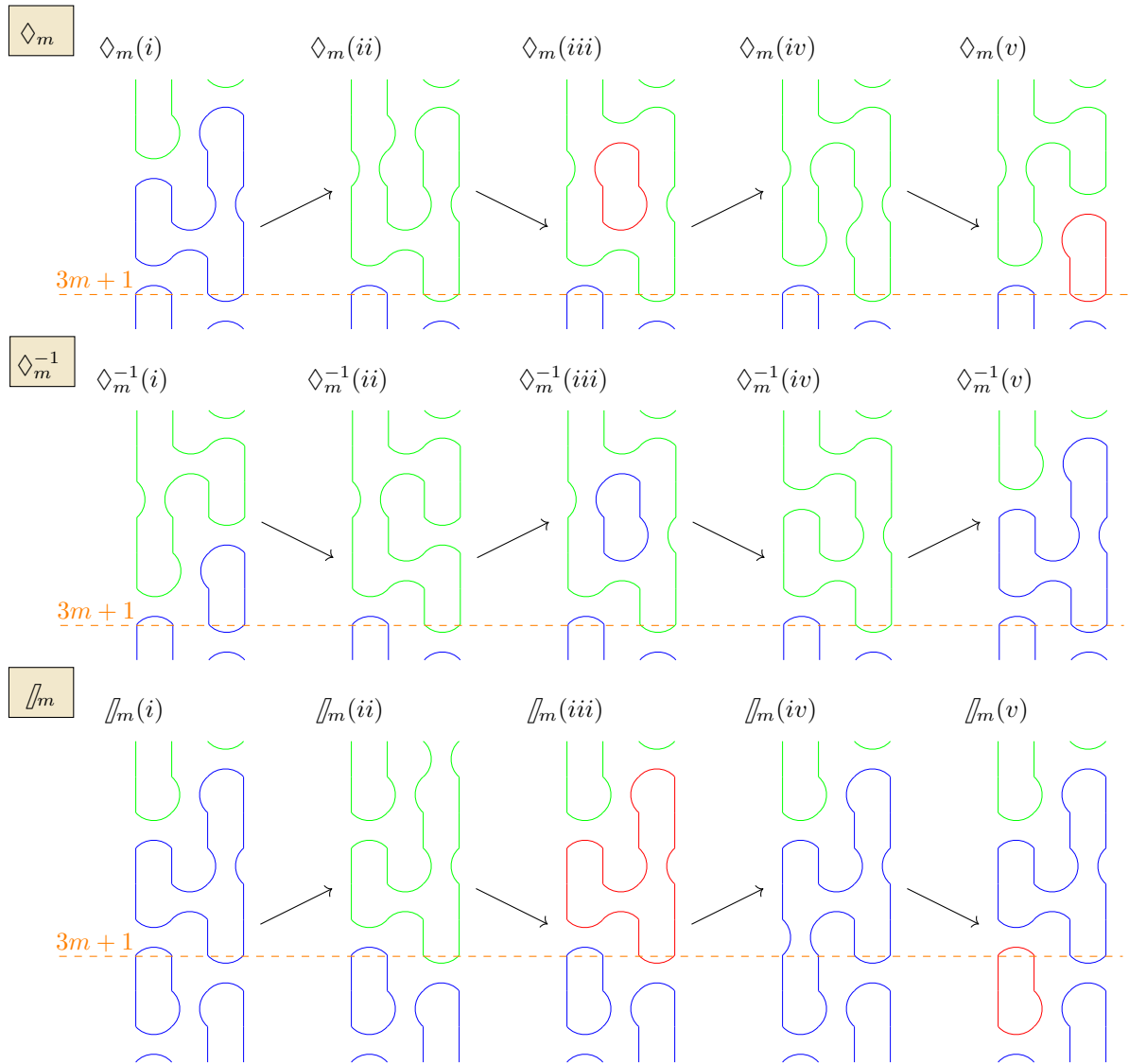

    \centering
    \input{diamond/diamondtext}
    \input{invdiamond/invdiamondtext}
    
    \input{jump/jumptext}
    \caption{Local moves part I. The typographic naming of the moves tries to mimic their patterns visually: the move $\lozenge_m$ consists of 4 edges which change the smoothings at  4 different crossings forming a diamond. In the move $\fatslash_m$ there are two changes of smoothings in the top right and two in the bottom left and the move $\triangledown_m$ is a half of $\lozenge_m$. The superscripts in $\lozenge_m^{-1}, \ \triangledown_m^{-1}, \ \overline{\triangledown}^{-1}$  suggest that these inverses of the moves $\lozenge_m, \ \triangledown_m, \ \overline{\triangledown}$ in an informal sense. The subscript $m$ determines the Morse layer, where the local move is applied. The slope on the arrows indicate whether they are going upwards or downwards in the homological degree. In the Morse complexes, all paths that contribute towards the differentials have to zig-zag between homological layer and hence all the local moves are also always zig-zagging. 
}
    \label{Figure: lego moves 1}
\end{figure}

\begin{figure}
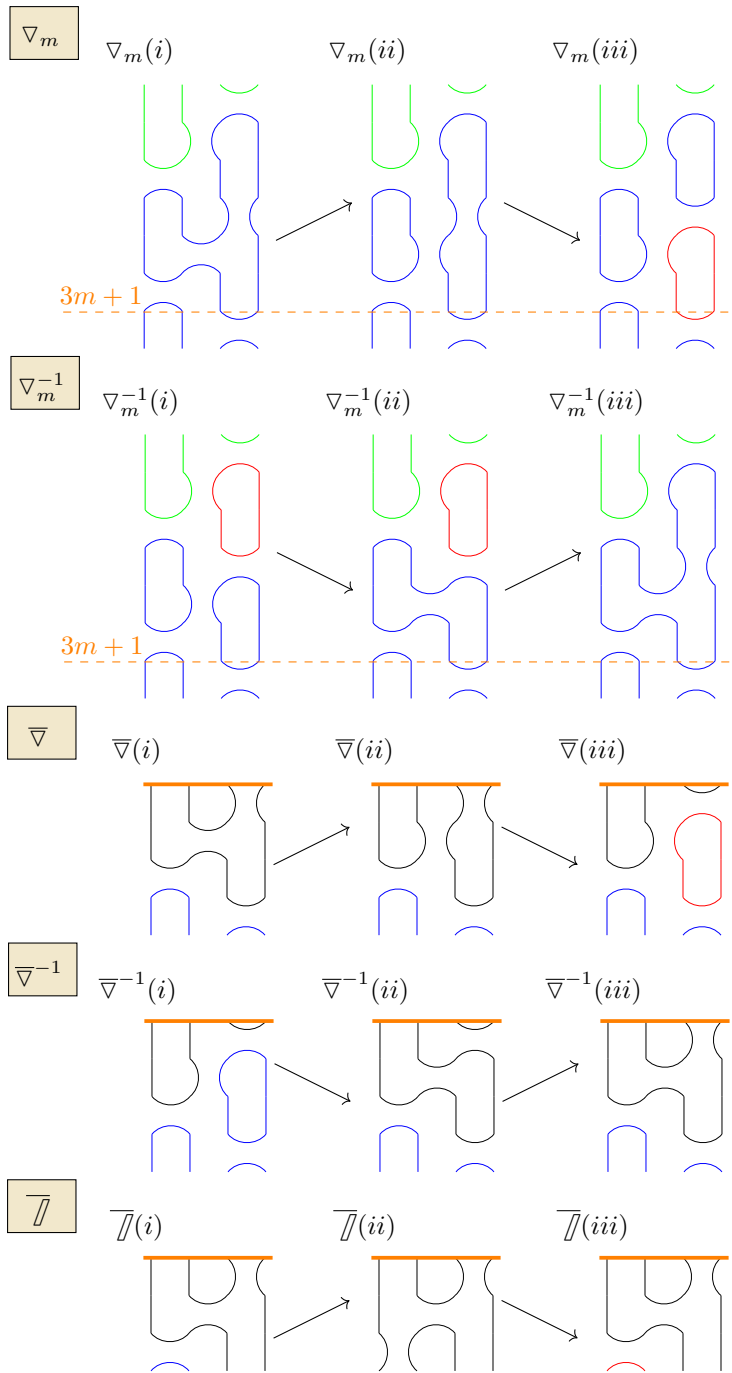

    \centering
    
    \input{tridown/tridowntext}
    \input{invtri/invtritext}
    
    \input{bartridown/bartridowntext}
    \input{invbartri/invbartritext}
    
    \input{bjump/bjumptext}
    \caption{Local moves part II. The moves $\overline{\triangledown}, \  \overline{\triangledown}^{-1}, \ \overline{\fatslash}$ are always located at the very top of the braid which is drawn with an orange line in the pictures. 
    \label{Figure: lego moves 2} }
\end{figure}

As previously stated, $\Phi_n$ cannot be built only by using graph isomorphisms, e.g. suitable maps $\Lambda$. To patch this lack of bijectivity, Figures \ref{Figure: lego moves 1} and \ref{Figure: lego moves 2} present additional local moves $\lozenge_m \triangledown_m \fatslash_m \lozenge_m^{-1}, \ \triangledown_m^{-1}, \ \overline{\triangledown}, \  \overline{\triangledown}^{-1}, \ \overline{\fatslash} $  which are employed to describe $\Phi_n$, whenever it fails to map vertices of paths in a one-to-one correspondence. The local moves are chained together to stretch the paths and make them longer. For example by stating that a path $p=(p_1\to\dots \to p_{2k})\in Z_n$ satisfies equation 
$$
p=\overline{\triangledown} \fatslash_{n-3} B_{n-3}(p)  \overline{\triangledown}^{-1}
$$
we mean that $p_1\to p_2\to p_3$ follows $\overline{\triangledown}$, then $p_3\to p_7$ follows $\fatslash_{n-3}$, then $p_7\to p_{2k-2}$ follows $B_{n-3}(p)$ and finally $p_{2k-2}\to p_{2k}$ follows $\overline{\triangledown}^{-1}$. More specifically, Morse layers $3n-11$ to $3n-3$ of $p_3$ are as in picture $\overline{\triangledown}(\romannumeral 3)$ and Morse layers $3(n-3)$ to $3(n-1)$ of the same vertex $p_3$ are as in picture $\fatslash_{n-3}(\romannumeral 1)$.

Finally we can define $\Phi_n$ and $A_n$. Assuming that $p\in Z_n$ is a path for which $B_{n-2}(p)$ exists we assign
\begin{alignat}{3}
&p=B_{n-2}(p)    &&\quad \implies  \quad   &&\Phi_n(p)=\Lambda(B_{n-2}(p)) \label{PhiDef 1} \\ 
&p=\fatslash_{n-2} B_{n-2}(p)    &&\quad \implies  \quad   &&\Phi_n(p)=\fatslash_{n+2} \lozenge_n  \fatslash_{n-2} \Lambda(B_{n-2}(p))  \lozenge_n^{-1} \label{PhiDef 2}\\ 
&p=\lozenge_{n-2} B_{n-2}(p)\lozenge_{n-2}^{-1}   &&\quad \implies  \quad &&     \Phi_n(p)=\lozenge_{n+2}\fatslash_{n} \lozenge_{n-2} \Lambda(B_{n-2}(p))  \lozenge_{n-2}^{-1}\lozenge_{n+2}^{-1} \label{PhiDef 3}\\ 
&p=\triangledown_{n-2} B_{n-2}(p)\overline{ \triangledown}^{-1}    &&\quad \implies  \quad   &&\Phi_n(p)=\triangledown_{n+2}\triangledown_{n}\triangledown_{n-2} \Lambda(B_{n-2}(p)) \triangledown_{n-1}^{-1}\triangledown_{n+1}^{-1}\overline{\triangledown}^{-1}  \label{PhiDef 4}
\end{alignat}
and assuming that $B_{n-3}(p)$ exists we assign
\begin{alignat}{3}
&p=\overline{\triangledown} \triangledown_{n-3} B_{n-3}(p) \triangledown_{n-2}^{-1}    &&\quad \implies  \quad   &&\Phi_n(p)=\overline{\triangledown}\triangledown_{n+1}\triangledown_{n-1} \triangledown_{n-3} \Lambda(B_{n-3}(p)) \triangledown_{n-2}^{-1}\triangledown_{n}^{-1} \triangledown_{n+2}^{-1} \label{PhiDef 5}\\ 
&p=\overline{\triangledown} \fatslash_{n-3} B_{n-3}(p)  \overline{\triangledown}^{-1}   &&\quad \implies  \quad   &&\Phi_n(p)=\overline{\triangledown}\fatslash_{n+1}\lozenge_{n-1} \fatslash_{n-3} \Lambda(B_{n-3}(p)) \lozenge_{n-1}^{-1} \overline{\triangledown}^{-1} \label{PhiDef 6}\\ 
&p=\overline{\fatslash} \, \lozenge_{n-3} B_{n-3}(p) \lozenge_{n-3}^{-1}    &&\quad \implies  \quad   &&\Phi_n(p)=\overline{\fatslash}\, \lozenge_{n+1}\fatslash_{n-1} \lozenge_{n-3} \Lambda(B_{n-3}(p)) \lozenge_{n-3}^{-1}\lozenge_{n+1}^{-1}  \label{PhiDef 7}\\ 
&p=\fatslash_{n-3} B_{n-3}(p)    &&\quad \implies  \quad   &&\Phi_n(p)=\fatslash_{n+1} \lozenge_{n-1}  \fatslash_{n-3} \Lambda(B_{n-3}(p))  \lozenge_{n-1}^{-1}.  \label{PhiDef 8}
\end{alignat}
Denote $A_n$ as the subset of $Z_n$ for which we have defined $\Phi_n$ in Implications \ref{PhiDef 1}-\ref{PhiDef 8}. The subscripts of $\Lambda$ are omitted and they are to be understood from the context: for example in Implication \ref{PhiDef 2} we mean $\Lambda=\Lambda_{s,t,n}$ where $vs$ is the 5th vertex of $p$,  $wt$ is $13$th  vertex of $\Phi_n(p)$ and $|v|=|w|=3(n-2)$.
This concludes the definitions of $\Phi_n$ and $A_n$ and we have therefore completed Task 1.

\subsubsection{Task 2:  \texorpdfstring{$\Phi_n$}{Phin} is well defined and invariant under \texorpdfstring{$R$}{R}} 
We will now begin to make use of the assumption $t_{\bC}(a)\geq 6$. Recall that by Lemma \ref{Lemma: t value implies W} and the classification of unmatched cells the inequality $t_{\bC}(a)\geq 6$ implies that 
$$
a=v(01\bzero \bzero01)^{12} \quad \text{ or } \quad a=v(01\bzero \bzero01)^{12} 01\bzero
$$
for some $v\in \mathcal{S}$. Inside this repeating region, it will be easier for us to control how edges and paths behave. The next combinatorial lemma confines the quantum degree of the relevant zig-zag paths.
\begin{lemma}\label{Lemma: q goes up at most 3 in a path}
    Let $a,b$ be unmatched cells of $U_n$ with $\operatorname{hdeg}(a)+1=\operatorname{hdeg}(b)$. Assume that there exists a path from $a$ to $b$ in $G(\Psi \llbracket (\sigma_1\sigma_2 \sigma_3)^n\rrbracket, \Mgr)$. Then $\operatorname{qdeg}(b)\leq \operatorname{qdeg}(a)+3$. 
\end{lemma}
\begin{proof}
    Denote $o(w)$ as the number of characters 1 at the start of the word $w$ for any vertex $w$ in  $G(\Psi \llbracket (\sigma_1\sigma_2 \sigma_3)^n\rrbracket, \Mgr)$. It follows from Lemma \ref{lemma: red zeros matched} that $o$ is increasing along edges and paths. Hence it suffices to show that for all $n\geq 0$ and unmatched cells $a,b$ of $U_n$ the following implication holds:
    \begin{equation}
    h(a)+1=h(b), \  o(a)\leq o(b) \implies q(b)\leq q(a)+3. \label{Equation: h O implies q}    
    \end{equation}
    The differences of $h,o,q$ do not vary under $\gamma_{\ast}$ recursions, that is, for all $\ast \in \{\bA,\bB, \bC\}$ and unmatched cells $a,b$ we have:
    \begin{align*}
        h(\gamma_{\ast}(a))-h(\gamma_{\ast}(b))&=h(a)-h(b) \\
        o(\gamma_{\ast}(a))-o(\gamma_{\ast}(b))&=o(a)-o(b) \\
        q(\gamma_{\ast}(a))-q(\gamma_{\ast}(b))&=q(a)-q(b).
    \end{align*}
    Thus if Implication \ref{Equation: h O implies q} holds for a pair $a,b$, then for all  $\ast \in \{\bA,\bB, \bC\}$, it will hold for the pair $\gamma_{\ast}(a),\gamma_{\ast}(b)$ too.

    We start out by verifying Implication \ref{Equation: h O implies q} for all $n=0,\dots, 82$ with \texttt{basecases.py}. Next we assume that Implication \ref{Equation: h O implies q} holds for $n-4$ where $n\geq 83$ and take $a,b$ to be unmatched cells of $U_n$ with $h(a)+1=h(b)$ and $o(a)\leq o(b)$. 

    \textbf{Case:} $t_{\bA}(a)\geq 2$. From Lemma \ref{Lemma: t value implies W} it follows that $a=\gamma_{\bA}(x)$ for some unmatched cell $x$ of $C_{n-4}$. Since $o(a)\leq o(b)$, it is easy to see that $b=\gamma_{\bA}(y)$ for some $y$ also. Since Implication \ref{Equation: h O implies q} holds for $x,y$ by the induction assumption, we get $q(b)\leq q(a)+3$.

    \textbf{Case:} $t_{\bA}(a)< 2$, $t_{\bC}(b)< 2$. From our assumptions and Lemma \ref{Lemma: minimal tA,tC} we get the following set of equalities and inequalities:
    \begin{multline*}
        -\frac{3}{2}\leq t_{\bA}(a) < 2, \quad
        -\frac{3}{2}\leq t_{\bC}(a),\quad
        -\frac{3}{2}\leq t_{\bA}(b),\quad         
        -\frac{3}{2}\leq t_{\bC}(b) < 2, \quad 
        h(a)+1=h(b),\quad n\geq 83.       
    \end{multline*}
    By using \texttt{Lean}'s linarith tactic, these can be seen to imply $t_{\bB}(a)\geq 2$ and $t_{\bB}(b)\geq 2$. It follows that there exists $x,y$ with $b=\gamma_{\bB}(x)$, $a=\gamma_{\bB}(y)$ and hence  $q(b)\leq q(a)+3$ holds by the induction assumption.

    \textbf{Case:} $t_{\bC}(b)\geq 2$, $t_{\bC}(a)\geq 2$. As in the previous case, this immediately suffices.

    \textbf{Case:} $t_{\bC}(b)\geq 2$, $t_{\bC}(a)< 2$. Opening up the definition of $t_{\bC}$ gives 
    $$
    -\frac{9}{4}n +h(a) -\frac{3}{4}q(a) -\frac{1}{4}< 2 \leq -\frac{9}{4}n +h(a)+1 -\frac{3}{4}q(b) -\frac{1}{4}
    $$
    which directly implies $q(a)+\frac{4}{3} > q(b)$.
\end{proof}

In order to prove that the maps $\Lambda$ induce graph isomorphisms between graphs $G$, we need to partially classify the edges of $\Mgr$ in $G(\Psi \llbracket \sigma_1 \sigma_2 \sigma_3)^n \rrbracket)$.

\begin{figure}
    \centering
    \input{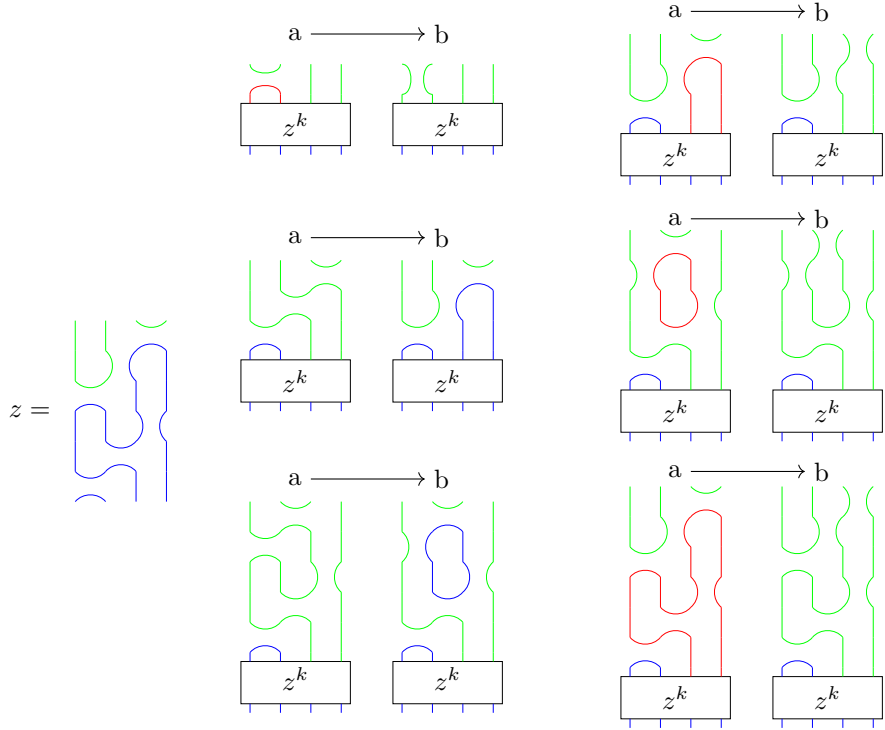}
    \caption{Classification of reversed arrow with $u$ evaluated in the repeating region. The first Morse layer of the pictured subword $z^k=(\bzero 0101\bzero)^k$ is $3m-5$ or $3m-8$ in all of the diagrams.}
    \label{Figure: reversed arrows with high t_C}
\end{figure}

\begin{lemma}\label{Lemma: Classifying reversed arrows with high t_C}
    Let $(a\to b)\in \Mgr$ for the graph $G(\Psi \llbracket (\sigma_1 \sigma_2 \sigma_3)^n \rrbracket)$. Suppose that $u(a\to b)> 3m$ and write $a=c_1\dots c_{3n}$, $b=d_1\dots d_{3n}$ as strings of characters. Assume that $t_{\bC} (c_1\dots c_{3m})\geq 1$. Then $(c_1\dots c_{3m}) =(d_1 \dots d_{3m})\in U$ and $a\to b$ is of one of the following the types in Figure \ref{Figure: reversed arrows with high t_C} for some $k\geq 1$.
In particular $i(a\to b)> 3m$.
\end{lemma}
\begin{proof}
    Since $t_{\bC}(c_1\dots c_{3m})\geq 1$, it follows from Lemma \ref{Lemma: t value implies W} that $c_1\dots c_{3m}$ contains at least one instance of the subword $\bC= (01\bzero \bzero 01)^2$. By comparing this with the definition of $W$ one can see that $c_{3m-6}\dots c_{3m}=0\bzero 0101\bzero$ or $c_{3m-9} \dots c_{3m}=0\bzero 0101\bzero\bzero01$. In the proof of Lemma \ref{unmatched cell classification lemma}, we ran Algorithm \ref{Algorithm: matching a single cell in Mgr} for such cells to classify the unmatched cells $U$. As a by-product we obtained all of the relevant matched arrows in Pictures  \ref{Picture: classification of arrows with high tC part 1} and \ref{Picture: classification of arrows with high tC part 2}. Those six arrows are exactly the ones depicted in Figure \ref{Figure: reversed arrows with high t_C}.
\end{proof}

Notice that all of the arrow types of Figure \ref {Figure: reversed arrows with high t_C} are represented in the local moves, Figures \ref{Figure: lego moves 1} and \ref{Figure: lego moves 2}. However, in the graph $G(\Psi \llbracket (\sigma_1\sigma_2\sigma_3)^n \rrbracket, \Mgr)$ these arrows are reversed so the downward arrows from Figures \ref{Figure: lego moves 1} and \ref{Figure: lego moves 2} have opposite direction compared to the ones pictured in Figure \ref {Figure: reversed arrows with high t_C}.

We define a function $\tau$ on $\mathcal R$ so that computing $t_{\bC}$ can be done in parts:
\begin{equation}\label{Equation: connection of tC and remainder}
t_{\bC}(v)+\tau(u)=t_{\bC}(vu).    
\end{equation}
Concretely, this is obtained by assigning
$$
\tau\colon \mathcal R \to \mathbb Q, \ \tau(r)=-\frac{1}{4}|r|+\frac{1}{4} (\#1 \text{ in } r) +\frac{3}{4} (\#\bzero \text{ in } r) -\frac{3}{4} (\#\rzero\text{ in } r).
$$ 
In order to control the $t_{\bC}(v)$ while having a handle on the a whole word $w=vr$ and the global environment $r \in \mathcal{R}$, we set 
$$
\mathcal{T}\colon X\times \mathcal Q \to \mathbb R , \ \mathcal T(w,r)= t_{\bC}(w) -\max \{\tau(s) \mid s\approx r \} -\frac{9}{4}
$$
where $X=\bigcup_{n=0}^\infty \operatorname{cells}(\Psi \llbracket (\sigma_1 \sigma_2 \sigma_3)^n \rrbracket)$. Recall that given $s\in \mathcal R $ and $n\in \mathbb Z_{\geq 0}$ we wrote 
$$
V_{s,n}=\{vr \in \operatorname{cells} \Psi \llbracket (\sigma_1\sigma_2 \sigma_3)^n \rrbracket \mid r\approx s, \ t_{\bC}(v)\geq 1  \}.
$$
This set of vertices accompanied with a restricted set of edges
$$
    E_{s,n}= \{ (a\to b) \in G(\Psi \llbracket (\sigma_1\sigma_2 \sigma_3)^n\rrbracket, \Mgr )\mid a,b\in V_{s,n},\  i(a\to b)\leq 3n- |s|   \}. 
$$
together make up a useful subgraph $G_{s,n}=(V_{s,n},E_{s,n})$ of the graph $G(\Psi \llbracket (\sigma_1\sigma_2 \sigma_3)^n\rrbracket, \Mgr)$. From these subgraph, the paths can be copied with the graph isomorphisms $\Lambda$. The next lemma gives us a computable condition which ensures that certain subpaths are contained in $G_{s,n}$.

\begin{lemma}\label{Lemma: beta subpath in G and other properties}
    Let $p=(p_1\to \dots \to p_{2k})\in Z_n$ be a zig-zag path and for each $j$ write $p_j=v_jr_j$ for some  $v_j\in \operatorname{cells}(\Psi \llbracket (\sigma_1 \sigma_2 \sigma_3 )^m \rrbracket)$ and  $r_j\in \mathcal{R}$. Suppose $i(p_{2l+1} \to p_{2l+2}) \leq 3m$ and $\mathcal T(p_1, r_{2l+1}) \geq 1$ and denote $b= \beta_{m} (p,2j+1)$. Then
    \begin{enumerate}
        \item $p_{2l+1}\to \dots \to p_{2b}$ is a path in $G_{r,n}$ where $r=r_{2l+1}$ \label{Claim 1 in combined T path after beta}
        \item $r_{2l+1} \approx r_{2b}$ \label{Claim 2 in combined T path after beta}
        \item $v_{2b}$ is either of the form $c_1\dots c_{3m-6}  \bzero 01 01\bzero$ or $c_1\dots c_{3m-9} \bzero 01 01\bzero \bzero10 $ for some characters $c_i\in \{ 0,\bzero, 1\}$. \label{Claim 3 in combined T path after beta}
        \item If $p_{2b} \notin U$, then $i(p_{2b} \to p_{2b+1} ) >3m$ and $p_{2b} \to p_{2b+1}$ is one of the arrows in Figure \ref {Figure: reversed arrows with high t_C}. \label{Claim 4 in combined T path after beta}
    \end{enumerate}
\end{lemma}
\begin{proof}
    Claim \ref{Claim 1 in combined T path after beta} is proven with an induction on the subpath $p_{2l+1}\to\dots  \to p_{2b}$ and it is not difficult to see that Claims \ref{Claim 2 in combined T path after beta}-\ref{Claim 4 in combined T path after beta} follow from it. From Lemma \ref{Lemma: q goes up at most 3 in a path} we get that for all $j$, 
    $$
    (h(p_j),q(p_j))\in (h(p_1),q(p_1))+ ([0,1]\times [0,3])
    $$
    and hence $t_{\bC}(p_j)\in t_{\bC}(p_1)+[-\frac{9}{4},1]$. Using it, we can compute
    \begin{equation}
    t_{\bC}(v_{2l+1})= t_{\bC}(p_{2l+1}) - \tau({r_{2l+1}})\geq t_{\bC}(p_1) -\frac{9}{4}  -\max \{\tau(s) \mid s\approx r \} \geq \mathcal{T}(p_1,r_{2l+1}) \geq 1 \label{Equation: tc tau T calc}    
    \end{equation}
    which means $p_{2l+1}\in G_{r,n}$. Since $i(p_{2l+1}\to p_{2l+2})\leq 3m$ we get $r_{2l+1}\approx r_{2l+2}$ and thus the lower bound \ref{Equation: tc tau T calc} can be applied for $v_{2l+2}$ as well, yielding $(p_{2l+1} \to p_{2l+2}) \in G_{r,n}$.
    
    Assume next that $(p_{2l+1}\to \dots \to p_{2t}) \in G_{r,n}$ for some $t>l$ and that $i(p_{2t} \to p_{2t+1}) \leq 3m$. Then by Lemma \ref{Lemma: Classifying reversed arrows with high t_C} $u(p_{2t} \to p_{2t+1}) \leq 3m$  and hence we can get from Lemma \ref{Lemma: u/i} that $i(p_{2t+1} \to p_{2t+2}) \leq 3m$. It follows that $r_{2t} \cong r_{2t+1} \cong r_{2t+2}$ and so  $(p_{2l+1}\to \dots \to p_{2t+2}) \in G_{r,n}$ which proves Claim \ref{Claim 1 in combined T path after beta} by induction.
 \end{proof}

We can now proceed in showing that in certain cases the maps $\Lambda$ are graph isomorphisms.
\begin{lemma}\label{Lemma: bijectivity of cores}
    Let $n\in \mathbb Z_{\geq 0}$ and $s,t\in \mathcal{R}$ with $s\sim t$. Then $\Lambda_{s,t,n}\colon  V_{s,n} \to V_{t,m}$ induces a graph isomorphism from $G_{s,n} \to G_{t,m}$ where $m=n+\frac{1}{3}(|t|-|s|)$. Thus any path $p$ in $G_{s,n}$ will be mapped to an orientation preserving path $\Lambda_{s,t,n}(p)$ and in addition $R(p)=R(\Lambda_{s,t,n}(p))$.
\end{lemma}
\begin{proof}
    It is immediate that $\Lambda=\Lambda_{s,t,n}$ is a bijection and that, both edge sets are the same if one disregards the direction of the arrows. Therefore to show that $\Lambda$ is a graph isomorphism, it suffices to show that
    \begin{equation} \label{Equation: Lambda Mgr iff}
    (a\to b) \in \Mgr \quad \iff \quad  (\Lambda(a)\to \Lambda(b))\in \Mgr.    
    \end{equation}
    for all edges $(a\to b)\in \Mgr$. In case $a\to b$ is not an isomorphism, it is immediate that Equivalence \ref{Equation: Lambda Mgr iff} holds since neither $a\to b$ nor $\Lambda(a) \to \Lambda(b)$ are contained in their respective greedy matchings. 

    Thus we can assume that $a\to b$ is an isomorphism and writing $u(a\to b)$ makes sense. If $u(a\to b)\leq 3n- |s|$ then Equivalence \ref{Equation: Lambda Mgr iff} holds by Lemma \ref{Lemma: equivalence of matching based on first characters}. On the other hand  $u(a\to b)> 3n- |s|\geq i(a\to b)$ is impossible by Lemma \ref{Lemma: Classifying reversed arrows with high t_C}.
\end{proof}

\begin{proof}[Proof for Task 2]
A straightforward evaluation of $\mathcal T$ (albeit tedious due to the number of cases)  shows that for all paths $p\in A_n$ and $p\in \operatorname{im} \Phi_n$ we have 
$$
\mathcal T(p_1,r_{2l+1})\geq 1 
$$
where $p_1$ is the first vertex of $p$ and $v_{2l+1}r_{2l+1}$ is the first vertex of the either subpath $B_{n-2}(p)$ or $B_{n-3}(p)$ and $\operatorname{length(r_{2l+1})}$ is $6$ or $9$ respectively. Thus, it follows from Lemma  \ref{Lemma: beta subpath in G and other properties} that the subpaths are contained in their respective graphs $G_{s,n}$. Hence, we can use Lemma \ref{Lemma: bijectivity of cores} to see that $\Lambda(B_{n-2}(p))$ and $\Lambda (B_{n-3}(p))$ are orientation preserving subpaths and that $R(B_{n-2}(p))=R(\Lambda(B_{n-2}(p)))$ and $R(B_{n-3}(p))=R(\Lambda(B_{n-3}(p)))$ as morphisms in $\operatorname{Mat}(\Cob (8))$. Since the local moves are always orientation preserving subpaths, it follows that $\Phi_n(p)$ is an orientation preseving path for all $p\in A_n$.

In order to wrap up the proof, we still need to show that the first and last steps of $p$ and $\Phi_n(p)$ agree under $R$. This in turn is straightforward evaluation of the functor $R$. More concretely, say in the case of Implication \ref{PhiDef 3}, given 
$$
(p_1\to \dots \to p_{2k})=p =\lozenge_{n-2} B_{n-2}(p) \lozenge_{n-2}^{-1}, \quad \text{ and } \quad (q_1\to \dots \to q_{2m} ) =\Phi_n(p)
$$
one can see that $R(p_1\to p_5)$ and $R(q_1\to q_{13})$ are equal as cobordisms (including signs!) and that similarly $R(p_{2k-4} \to p_{2k})=R(q_{2m-8} \to q_{2m})$.  Analogous statements can be verified for Implications \ref{PhiDef 2} and \ref{PhiDef 4}-\ref{PhiDef 8}. 
    \end{proof}

\subsubsection{Task 3: (co)restrictions to \texorpdfstring{$\Phi_n$}{Phin} coincide with the differentials} 
The remaining task in order to prove Proposition \ref{Proposition: ac matrix element equations}, case $\mathfrak c_n$, is to show that under the functor $R$ the set of all zig-zag paths $Z_n$ and $Z_{n+4}$ behave similarly to the sets of paths $A_n$ and $\operatorname{im} \Phi_n$. This is the most laborious part, since it consists of classifying all of the paths of $Z_n$ and $Z_{n+4}$ to an extent. At first, we will work at the minuscule edge-level scale; then the edges are patched up to medium-sized blocks and finally the large scale classification of paths is completed by sewing together the medium-sized blocks.



In any path zig-zag $(p_1\to \dots \to p_{2k})\in Z_n$ a downward edge $p_{2j} \to p_{2j+1}$ is uniquely determined by the vertex $p_{2j}$ since the reversed edges of $\Mgr$ form a matching on the graph. In most cases, the upward edges $p_{2j+1} \to p_{2j+2}$ are determined by the value $i(p_{2j+1}\to p_{2j+2})$, only in the case of splitting a red circle, isomorphism iv) from Figure \ref{Figure: 2d local dictionary}, an additional choice of color needs to be made. This is useful: in Lemma \ref{Lemma: u/i} we obtained an upper bound for $i$  and the next lemma grants us lower bounds for it.

\begin{lemma}\label{Lemma: i lower bounds}
    Let $p=(p_1\to \dots \to p_{2k})\in Z_n$ be a zig-zag path, write $p_{2j+1}=c_1\dots c_{3n}$ and assume that $\mathcal T (p_1,c_{3m+1} \dots c_{3n}) \geq 1$. If $c_{3m+1}c_{3m+2}c_{3m+3}= \gzero 00$ or $c_{3m+1}\dots c_{3m+5}=\gzero 011\rzero$, where $\gzero \in \{ 0, \rzero, \bzero \}$,
    \begin{center}
        \input{ilowerbound/il13combined}
    \end{center}
    then $i(p_{2j+1}\to p_{2j+2})\geq 3m+1$.
    Also if in either case $i(p_{2j+1}\to p_{2j+2})= 3m+1$ then it cannot hold that $i(p_{2j+2}\to p_{2j+3})=u(p_{2j+2}\to p_{2j+3})= 3m+1$. Alternatively, if $i(p_{2j}\to p_{2j+1})>3m+6$ and $c_{3m+1}\dots c_{3m+5}=\bzero 0101\bzero$
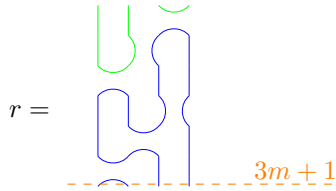
    \begin{center}
        \input{ilowerbound/il2withtext}
    \end{center}
    then $i(p_{2j+1}\to p_{2j+2})\geq 3m+1$. 
\end{lemma}
\begin{proof}
    Let us consider the first case $c_{3m+1}c_{3m+2}c_{3m+3}= \gzero 00$ and assume towards contradiction that $i(p_{2j+1}\to p_{2j+2})\leq  3m$. Let $b=\beta_{m}(p,2j+1) $ and write $p_{2b}=d_1\dots d_{3n}$. Lemma \ref{Lemma: beta subpath in G and other properties} yields that $d_{3m+1}\dots d_{3n}\approx c_{3m+1} \dots c_{3n}$ and that $p_{2b}\to p_{2b+1}$ is one of the arrows from Figure \ref {Figure: reversed arrows with high t_C}. The only arrows compatible with this data match $p_{2b}$ to a higher homological degree which contradicts the fact that $p$ is a zig-zag path. The case $c_{3m+1}\dots c_{3m+5}=\gzero 011\rzero$ is proven identically. 

    If in either of the cases above, we have  $i(p_{2j+1}\to p_{2j+2})= 3m+1$ and $i(p_{2j+2}\to p_{2j+3})=u(p_{2j+2}\to p_{2j+3})= 3m+1$, then we can use the claim again to get $3m+1 \leq i(p_{2j+3}\to p_{2j+4})$ and Lemma \ref{Lemma: u/i} to show that $i(p_{2j+3}\to p_{2j+4}) \leq 3m+1$. This is impossible; the arrow $p_{2j+3}\to p_{2j+4}$ has the wrong direction as we just traversed it from $p_{2j+2}$ to  $p_{2j+3}$.

    The last case $u(p_{2j}\to p_{2j+1})>3m+6$ and $c_{3m+1}\dots c_{3m+5}=\gzero 0101\bzero$ is also similar. If $i(p_{2j+1}\to p_{2j+2})< 3m+1$, assign $b=\beta_{m}(p,2j+1) $ and then $p_{2b}=d_1\dots d_{3n}$ will be matched up the homological degree due to the classification in Lemma \ref{Lemma: Classifying reversed arrows with high t_C}.
\end{proof}
Now, we can more concretely pin down implications regarding the next edges after subpaths.
\begin{lemma}\label{Lemma: path after beta}
    Let $p=(p_1\to \dots \to p_{2k})\in Z_n$ be a zig-zag path and  write $p_{2j+1}=c_1\dots c_{3n}$. Assume that $\mathcal T (p_1,c_{3m+1} \dots c_{3n}) \geq 1$ and  $i(p_{2j+1}\to p_{2j+2})\leq 3m$ and denote $b=\beta(p,2j+1)$. The following implications determine the Morse layers of $p_{2b}$ based on the Morse layers of $p_{2j+1}$:
    \begin{center}
        \input{betaclass2/implications1}
        \hspace{0.5 cm}
        \input{betaclass2/implications3}
        
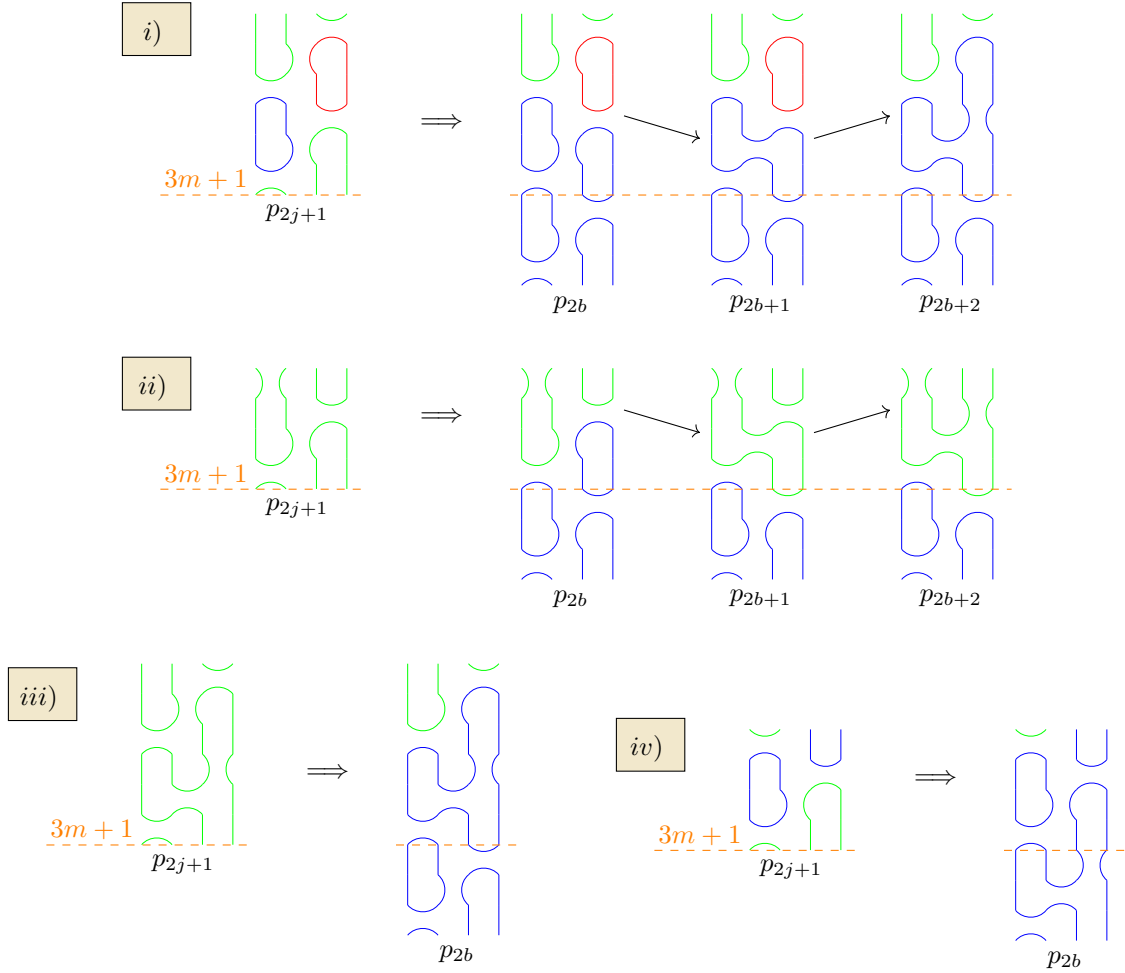
\hspace{0.5 cm}

\begin{tikzpicture}
\node at (0, 0) { 
            \input{betaclass2/implications4}
        };  
\node at (8, -0.42) {

\begin{tikzpicture}
\node at (0, 0) [anchor=south west, scale=0.4] { 
            \input{betaclassification/beta2a}
        };  
\node at (3.5, -1.2) [anchor=south west, scale=0.4] { 
            \input{betaclassification/beta2b} 
        };

\node at (0.6,-0.2) {$p_{2j+1}$};
\node at (4.2,-1.4) {$p_{2b}$};

\node at (2.5,1) {$\implies$};

\draw[orange,dashed] (-1.2,0.05) -- (1.5,0.05);
\node[orange] at (-.6,0.25) {$3m+1$};
\draw[orange,dashed] (3.4,0.05) -- (5,0.05);

\begin{scope}[shift={(0,-3.8)}]
\definecolor{lighttan}{rgb}{0.95, 0.91, 0.80}
\fill[lighttan] (-1.7,4.9) rectangle (-0.8,5.6);
\draw (-1.7,4.9) rectangle (-0.8,5.6);
\node at (-1.3,5.2) {$iv)$};
\end{scope}

\end{tikzpicture}

        };

\end{tikzpicture}

    \end{center}
    If additionally $\mathcal T (p_1,c_{3m+1} \dots c_{3n}) \geq \frac{3}{2}$, then
    \begin{center}
        \input{betaclass2/implications2}
    \end{center}
\end{lemma}
\begin{proof}
    All of the implications follow from Lemma \ref{Lemma: beta subpath in G and other properties}: In $i)$, $ii)$ and $v)$, $p_{2b}\to p_{2b+1}$ is the only arrow pictured in Figure \ref {Figure: reversed arrows with high t_C} which is downward pointing and compatible with $c_{3m+1}\dots c_{3m+4}$ or $c_{3m+1}\dots c_{3m+6}$. In those cases, the second arrow $p_{2b+1}\to p_{2b+2}$ can be pinned down by excluding other possibilities with Lemmas \ref{Lemma: u/i} and \ref{Lemma: i lower bounds}. The distinct bound in Implication $v)$ follows from the fact that $\mathcal{T}(p_1,d_{3m+1}\dots d_{3n})+\frac{1}{2}= \mathcal{T}(p_1,c_{3m+1}\dots c_{3n})$ where $d_{3m+1}\dots d_{3n}$ are the last characters of $p_{2j+3}$. Hence one needs extra cushion in order to guarantee that Lemma \ref{Lemma: i lower bounds} can be applied. In Implications $iii)$ and $iv)$, the depicted $p_{2b}$ is the only configuration which does not get matched up in homological degree.
\end{proof}

\begin{lemma} \label{Lemma: middle exclusion i}
    Let $p=(p_1\to \dots \to p_{2k} ) \in Z_n$ be a zig-zag path and write $p_1=c_1\dots c_{3n}$. If $c_{3m+1} \dots c_{3m+9}=\bzero0101\bzero \bzero 01$ 
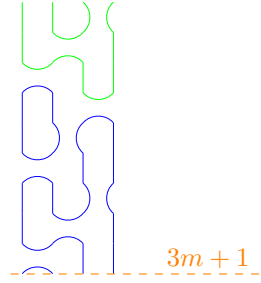
    \begin{center}

\begin{tikzpicture}
\node at (0, 0) [anchor=south west, scale=0.4] {  
            \input{betaclassification/betaonehalf} 
        };

\draw[orange,dashed] (-0.1,0.05) -- (3.2,0.05);
\node[orange] at (2.5,0.25) {$3m+1$};

\end{tikzpicture}

    \end{center}
    and $\mathcal{T}(p_1, c_{3m+1} \dots c_{3n}) \geq \tfrac{5}{2}$, then $i(p_1\to p_2)\neq 3m+2$.
\end{lemma}
\begin{proof}
    Assuming $i(p_1 \to p_2)=3m+2$  leads to 
    \begin{center}
        \input{Proof_of_lego_implications/pic3}
    \end{center}
    where the $b=\beta_m(p,3)$. At $p_3$, Lemma \ref{Lemma: u/i} was used to deduce that $i(p_3\to p_4)\leq 3m+3$. There are no arrows with values $i(p_3\to p_4)=3m+1,3m+2,3m+3$: at $3m+1$ two blue circles are merged, $3m+2$ is already $1$-smoothed and at $3m+3$ the arrow has the wrong direction. By writing $p_3=d_1 \dots d_{3n}$, we can compute that 
    $$
    \frac{5}{2} \leq \mathcal{T}(p_1, c_{3m+1} \dots c_{3n}) =\mathcal{T}(p_1, d_{3m+1} \dots d_{3n})+\frac{3}{2}
    $$
    and thus Lemma \ref{Lemma: path after beta}  iv) can also be applied at $p_3$ to obtain $p_{2b}$ and $p_{2b+1}$ where $b=\beta_m(p,3)$. Then at $p_{2b+1}$ we can use Lemmas \ref{Lemma: u/i} and \ref{Lemma: i lower bounds} to pin down the $i$ value and acquire the edge $p_{2b+1}\to p_{2b+2}$. However, the vertex $p_{2b+2}$ is matched upwards and hence $p$ cannot continue as a zig-zag path between the homological levels. This renders the first edge with $i(p_1 \to p_2)=3m+2$ impossible.
\end{proof}

\begin{lemma}\label{Lemma: Lego moves}
    Let  $p=(p_1\to \dots \to p_{2k} ) \in Z_n$ be a zig-zag path. The following implications hold:
    \begin{alignat}{3}
& (p_{2j-1} \to p_{2j+1})= \barfatslash(\romannumeral 1 \to \romannumeral 3)  && \quad \implies \quad &&   (p_{2j+1}\to p_{2j+3})=\lozenge_{n-3}(\romannumeral 1\to\romannumeral 3) \label{Implication: barjumpmove -> diamond} \\  
& (p_{2j-1} \to p_{2j+1})= \overline{\triangledown}(\romannumeral 1 \to \romannumeral 3)  && \quad \implies \quad &&   (p_{2j+1}\to p_{2j+3}) \in \{ \triangledown_{n-3}(\romannumeral 1 \to \romannumeral 3), \, \fatslash_{n-3}(\romannumeral 1 \to \romannumeral 3)\} \label{Implication: barhalfdiamond -> halfdiamond or jump}\\
& (p_{2j} \to p_{2j+1})= \overline{\triangledown}^{-1}(\romannumeral 1 \to \romannumeral 2 ) && \quad \implies \quad &&   (p_{2j}\to p_{2j+2})=\overline{\triangledown}^{-1}(\romannumeral1 \to \romannumeral 3) \label{Implication: barinverse halfdiamond -> continues}
\end{alignat}
    If we further write $p_{2j+1}=c_1\dots c_{3n}$ 
    and  $T=\mathcal{T}(p_1,c_{3m+1}\dots c_{3n})$, then 
    \begin{alignat}{3}
T\geq 1, \quad& (p_{2j-1} \to p_{2j+1})= \lozenge_m(\romannumeral 1 \to \romannumeral 3)  && \quad \implies \quad &&   (p_{2j-1}\to p_{2j+3})=\lozenge_m(\romannumeral 1  \to \romannumeral 5) \label{Implication: diamond -> continues}\\      
T\geq 2, \quad& (p_{2j-1} \to p_{2j+1})= \fatslash_{m+2}(\romannumeral 3 \to \romannumeral 5)  && \quad \implies \quad &&   (p_{2j+1}\to p_{2j+3})=\lozenge_m(\romannumeral 1  \to \romannumeral 3)\label{Implication: jump -> diamond}\\ 
T\geq \tfrac{3}{2}, \quad& (p_{2j-3} \to p_{2j+1})= \lozenge_{m+2}(\romannumeral 1 \to \romannumeral 5)  && \quad \implies \quad &&   (p_{2j+1}\to p_{2j+3})=\fatslash_m(\romannumeral 1  \to \romannumeral 3)\label{Implication: diamond -> jump}\\ 
T\geq 3, \quad& (p_{2j-1} \to p_{2j+1})= \triangledown_{m+2}(\romannumeral 1 \to \romannumeral 3)  && \quad \implies \quad &&   (p_{2j+1}\to p_{2j+3})=\triangledown_m(\romannumeral 1  \to \romannumeral 3) \label{Implication: halfdiamond -> halfdiamond} \\
T\geq \tfrac{3}{2}, \quad& (p_{2j} \to p_{2j+1})= \lozenge_m^{-1}(\romannumeral 1 \to \romannumeral 2)  && \quad \implies \quad &&   (p_{2j}\to p_{2j+4})=\lozenge_m^{-1}(\romannumeral 1\to\romannumeral 5)\label{Implication: invdiamond -> continues}\\ 
T\geq 1, \quad& (p_{2j} \to p_{2j+1})= \triangledown_m^{-1}(\romannumeral 1 \to \romannumeral 2)  && \quad \implies \quad &&   (p_{2j}\to p_{2j+2})=\triangledown_m^{-1}(\romannumeral 1  \to \romannumeral 3).\label{Implication: invhalfdiamond -> continues}
\end{alignat}
\end{lemma}
\begin{proof}
    \textbf{Implication \ref{Implication: barjumpmove -> diamond}:} There are 4 meaningful possible values for $i(p_{2j+1}\to p_{2j+2})$: it is either $3n-5$, $3n-3$, $3n-1$ or  $i(p_{2j+1}\to p_{2j+2})\leq 3n-6$. The claim follows from the first option, so we can prove it by excluding the three others. If $i(p_{2j+1}\to p_{2j+2})=3n-3$, then one can see that the quantum degree increases too much: 
    $$
    q(p_1)\leq q(p_{2j+1})=q(p_{2j})-2 =q(p_{2j+2})-4 \leq q(p_{2k}) -4
    $$
    while Lemma \ref{Lemma: q goes up at most 3 in a path} imposes $q(p_1) \geq q(p_{2k})-3$. The value $i(p_{2j+1}\to p_{2j+2})=3n-1$ is excluded by Lemma \ref{Lemma: u/i}. In the case $i(p_{2j+1}\to p_{2j+2})\leq 3n-6$ we can use  Implication ii) of Lemma \ref{Lemma: path after beta} which makes $p_{2b+2}$ be matched upwards, where $b=\beta_{n-2}(p,2j+1)$.

    \textbf{Implications \ref{Implication: barinverse halfdiamond -> continues}, \ref{Implication: diamond -> continues}, \ref{Implication: invdiamond -> continues} and \ref{Implication: invhalfdiamond -> continues}} follow by excluding the other possible values of $i(p_{2j+1}\to p_{2j+2})$ with Lemmas \ref{Lemma: u/i} and \ref{Lemma: i lower bounds}. For the moment, we postpone the proofs of Implications \ref{Implication: barhalfdiamond -> halfdiamond or jump} and \ref{Implication: halfdiamond -> halfdiamond}.




    \textbf{Implication \ref{Implication: jump -> diamond}:} Using Implication $v)$ of Lemma  \ref{Lemma: path after beta} at $3m+3$ excludes $i(p_{2j+1})\to i(p_{2j+2})\leq 3m+3$ since otherwise $p_{2b+2}$ would be matched upwards, where $b=\beta_{m+1}(p,2j+1)$.

    \textbf{Implication \ref{Implication: diamond -> jump}:} There are 4 meaningful possible values for $i(p_{2j+1}\to p_{2j+2})$: it is either  $3m+6$, $3m+7$,  $i(p_{2j+1}\to p_{2j+2})\leq 3m+3$ or $i(p_{2j+1}\to p_{2j+2})\geq 3m+10$. We deem the case $3m+7$ impossible by considering the color of the top left green component in the $\lozenge_{m+2}(i)$ diagram. It cannot be blue, since then $p_{2j-3}\to p_{2j-2}$ would be a merge of two blue circles and thus not a arrow. It cannot be red, since then that component would be blue in the $\lozenge_{m+2}(v)$ which would make $p_{2j+1}\to p_{2j+2}$ with $i(p_{2j+1}\to p_{2j+2})=3m+7$ to be a merge of two blue circles. It also cannot be black, as that would increase the quantum degree too much as in the case of Implication \ref{Implication: barjumpmove -> diamond}. 

    The case $i(p_{2j+1}\to p_{2j+2})\geq 3m+10$ can be ruled out with Lemma \ref{Lemma: u/i}. If $i(p_{2j+1}\to p_{2j+2})\leq 3m+3$, then Implication $i)$ of Lemma \ref{Lemma: path after beta} at $3m+3$ forces $p_{2b+2}$ to be:
    
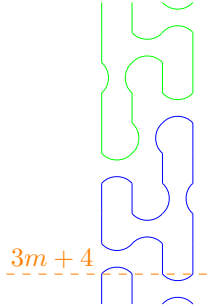
\begin{center}

\begin{tikzpicture}

\node at (0, 0) [anchor=south west, scale=0.4] { 
            \input{Proof_of_lego_implications/cabcabcabc000101x100}
        };  
\draw[orange,dashed] (-1.2,0.45) -- (1.5,0.45);
\node[orange] at (-.6,0.65) {$3m+4$};

\end{tikzpicture}

    \end{center}
    Then $i(p_{2b+2}\to p_{2b+3})=3m+10$ which causes a contradiction: By Lemma \ref{Lemma: i lower bounds} at $3m+10$, $i(p_{2b+3}\to p_{2b+4})\geq 3m+10$ and Lemma \ref{Lemma: u/i} forces $i(p_{2b+3}\to p_{2b+4})\leq 3m+10$. However, the value $i(p_{2b+3}\to p_{2b+4})= 3m+10$ is impossible, as that arrow is reversed.
\end{proof}

Given two unmatched cells $a,b\in U_n$, we say that two sets of paths $K,L\subset \operatorname{ZZ}(a,b)$  are $R$-equivalent and denote $K \overset{R}{\sim} L$, if 
    $$
    \sum_{p\in K} R(p)=\sum_{p\in L}R(p). 
    $$
This enables us to rewrite Claim 3 as 
$$
\operatorname{ZZ}(a,b) \overset{R}{\sim} \operatorname{ZZ}(a,b)\cap A_n \quad \text{and} \quad \operatorname{ZZ}(\mathfrak c_n(a),\mathfrak c_n(b)) \overset{R}{\sim} \operatorname{ZZ}(\mathfrak c_n(a),\mathfrak c_n(b))\cap \operatorname{im} \Phi_n.
$$
The following lemma addresses the case where $\Phi_n$ does not map all paths of $Z_n$ bijectively to all paths of $Z_{n+4}$. Instead, the mapping is done up to $R$-equivalence which suffices for our purposes.

\begin{lemma}\label{Lemma: jump move R-equivalence}
    Let $p_1\to\dots \to  p_{2j+1}$ be a zig-zag path in $G(\Psi \llbracket (\sigma_1 \sigma_2 \sigma_3)^n\rrbracket, \Mgr)$ such that $i(p_l\to p_{l+1})>3m$ for all $l$. Assume  $p_{2j+1}=v.r$ where $r\in \mathcal S$, $|v|=3m$ and $\mathcal{T}(p_1,r)\geq \tfrac{3}{2}$. Let $y\in U_n$ be an unmatched cell and denote 
    \begin{align*}
        K&=\{(q_1\to \dots \to q_{2k}) \in \operatorname{ZZ}(x,y) \mid q_l=p_l \text{ for } l=1,\dots,2j+1, \ (q_{2j-1}\to q_{2j+1})=\fatslash_m(\romannumeral 1 \to \romannumeral 3)\} \\ 
        L&=\{(q_1\to \dots \to q_{2k}) \in \operatorname{ZZ}(x,y) \mid q_l=p_l \text{ for } l=1,\dots,2j+1, \ (q_{2j-1}\to q_{2j+3})=\fatslash_m(\romannumeral 1 \to \romannumeral 5)\}.
    \end{align*}
    Then $K$ is $R$-equivalent to $L$.
\end{lemma}

\begin{proof}
Throughout this proof, given $q\in K$  we will denote $2b+1$ as the last indices of the subpath $B_{3m}(q)$. In addition to $L$, we shall define $4$ other subsets $N_1,\dots,N_4$ of $K$. Let the subset $N_1 $ consist of those zig-zag paths $q$ whose subpath $q_{2j+1}\to q_{2b+1}$ is of the form
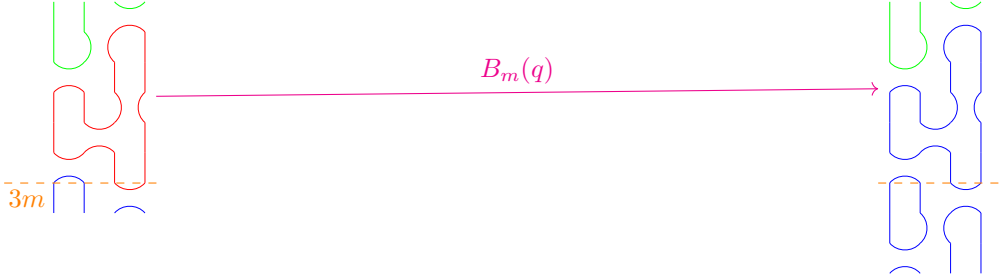
\begin{center}
    \input{JMC/N1}
\end{center}
the subset $N_2$ consist of those zig-zag paths $q$ whose subpath $q_{2j+1}\to q_{2b+3}$ is of the form
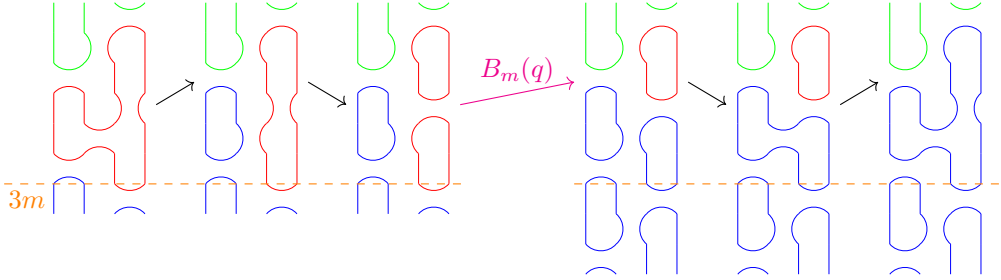
\begin{center}
\input{JMC/N2}
\end{center}
the subset $N_3$ consist of those zig-zag paths $q$ whose subpath $q_{2j+1}\to q_{2b+3}$ is of the form
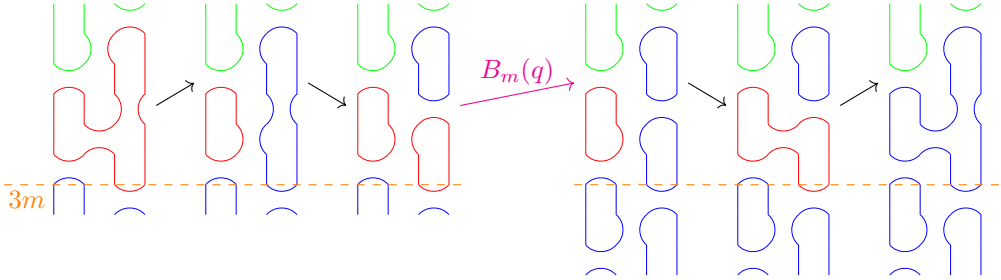
\begin{center}
    \input{JMC/N3}
\end{center}
the subset $N_4$ consist of those zig-zag paths $q$ whose subpath $q_{2j+1}\to q_{2b+5}$ is of the form
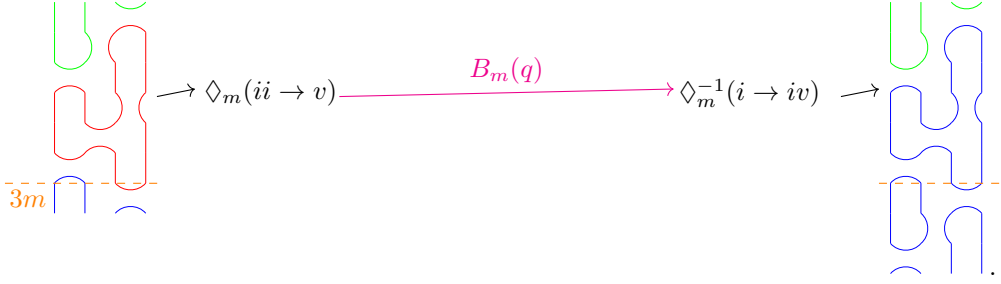
\begin{center}
    \input{JMC/N4}
\end{center}
Splitting into cases based on potential values of $i(q_{2j+1}\to q_{2j+2})$, and further using Lemmas \ref{Lemma: u/i}, \ref{Lemma: i lower bounds}, \ref{Lemma: path after beta} and Implications \ref{Implication: diamond -> continues} and \ref{Implication: invdiamond -> continues} from Lemma \ref{Lemma: Lego moves} it can be deduced that 
$$
K=L \sqcup N_1 \sqcup N_2  \sqcup N_3 \sqcup N_4.
$$

We can now define a bijection $\Omega_2\colon N_1\to N_2$ which sends $B_m(q)$ to $\Lambda(B_m(q))$, pads it with small subpaths and acts as identity elsewhere. More concretely, if $p=cB_{m}(q)d\in N_1$, then
$$
\Omega_2(q)=c \gamma_2\Lambda(B_m(q)) \delta_2 d
$$
where $\gamma_2$ and $\delta_2$ are the 2 edge subpaths from the definition of $N_2$. Similarly we define bijections $\Omega_3\colon N_1\to N_3$ and $\Omega_4 \colon N_1\to N_4$. It can be verified that $\Omega_i$ are well-defined by observing that the subpaths are glued along common vertices. Additionally, the maps $\Omega_i$ are bijections, since the maps $\Lambda$ are graph isomorphisms.

It follows from Lemma \ref{Lemma: bijectivity of cores} that $R(B_m(q))=R(\Lambda(B_m(q)))$ so by computing the value of $R$ on the padding subpaths one can see that 
$$
R(q)=-R(\Omega_2(q))=-R(\Omega_3(q))=R(\Omega_4(q))
$$
for all $q\in N_1$. Thus 
$$
\sum_{q\in K} R(q)=\sum_{q\in L}  R(q)+ \sum_{q\in N_1} \Big ( R(q)+R(\Omega_2(q))+R(\Omega_3(q))+R(\Omega_4(q)) \Big )=\sum_{q\in L} R(q)
$$
which proves the claim.
\end{proof}

\begin{proof}[Proof of Implications \ref{Implication: barhalfdiamond -> halfdiamond or jump} and \ref{Implication: halfdiamond -> halfdiamond} from Lemma \ref{Lemma: Lego moves}]
    In Implication \ref{Implication: halfdiamond -> halfdiamond}, there are 4 meaningful possible values for $i(p_{2j+1}\to p_{2j+2})$: it is either  $3m+2$, $3m+6$, or  $i(p_{2j+1}\to p_{2j+2})\leq 3m+3$ or $i(p_{2j+1}\to p_{2j+2})\geq 3m+10$. The two inequalities can be excluded with Lemma \ref{Lemma: i lower bounds} and \ref{Lemma: u/i} respectively, so we may assume towards contradiction that $i(p_{2j+1}\to p_{2j+2})=3m+6$. Then $i(p_{2j+2}\to p_{2j+3})=3m+6$, so $p_{2j+3}$ locally coincides with $q_{2j+1}$ from Lemma \ref{Lemma: jump move R-equivalence}. Hence we can use the same classification to deduce that all paths from $p_{2j+3}$ must arrive to a vertex $p_{2b}$ of the form:
    
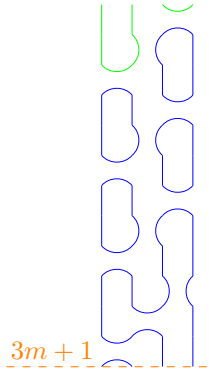
\begin{center}

\begin{tikzpicture}

\node at (0, 0) [anchor=south west, scale=0.4] { 
            \input{Proof_of_lego_implications/abcabcabcabc001010x1xx1x}
        };  
\draw[orange,dashed] (-1.2,0.05) -- (1.5,0.05);
\node[orange] at (-.6,0.25) {$3m+1$};

\end{tikzpicture}

    \end{center}
    and $i(p_{2b}\to p_{2b+1})=3m+8$. Then a contradiction is reached, since Lemmas \ref{Lemma: u/i} and \ref{Lemma: i lower bounds} enforce that $3m+7\leq i(p_{2b+1}\to p_{2b+2}) \leq 3m+9$, but there are no such arrows in the graph.  The candidate arrows are either merges of two blue loops, or a reversed which both make them non-viable.  Implication \ref{Implication: barhalfdiamond -> halfdiamond or jump} is proven similarly by making use of the classification in the proof of Lemma \ref{Lemma: jump move R-equivalence}.
\end{proof}

\begin{proof}[Proof for Task 3]
Let $a,b\in U_n$ be unmatched cells with $t_{\bC}(a)\geq 6$ and $t_{\bC}(b)\geq 1$. Since $t_{\bC}(a) \geq 1$, it follows that $a=v.\bzero0101\bzero$ or $a=v.\bzero 0101\bzero\bzero 01$ for some $v\in U$. Assume first that  $a=v.\bzero0101\bzero$, that is, 
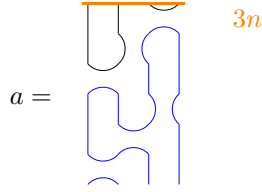
\begin{center}

\begin{tikzpicture}
\node at (0, 0) [anchor=south west, scale=0.4] {

    \begin{tikzpicture}
    \begin{scope}[color=blue]
    
        \draw (0, 0) to[bend left=50] (1, 0);
    \draw (2, 0) to (2, 1);
    \draw (1, 1) to[bend left=50] (2, 1);
    \draw (0, 1) to[bend right=50] (1, 1);
    \draw (0, 1) to (0, 2);
    \draw (0, 2) to (0, 3);
    \draw (0, 3) to[bend left=50] (1, 3);
    \draw (1, 2) to (1, 3);
    \draw (1, 2) to[bend right=50] (2, 2);
    \draw (2, 2) to[bend right=50] (2, 3);
    \draw (2, 3) to (2, 4);
    \draw (2, 4) to[bend left=50] (2, 5);
    \draw (2, 5) to[bend left=50] (3, 5);
    \draw (3, 4) to (3, 5);
    \draw (3, 3) to (3, 4);
    \draw (3, 2) to[bend left=50] (3, 3);
    \draw (3, 1) to (3, 2);
    \draw (3, 0) to (3, 1);
    \end{scope}

    \draw (0, 5) to (0, 6);
    \draw (0, 4) to (0, 5);
    \draw (0, 4) to[bend right=50] (1, 4);
    \draw (1, 4) to[bend right=50] (1, 5);
    \draw (1, 5) to (1, 6);
    \draw (2, 6) to[bend right=50] (3, 6);

\draw[orange,line width=1mm] (-0.2,6) -- (3.2,6);

    \end{tikzpicture}

        };  

\node[orange] at (2.25,2.25) {$3n$};
\node at (-0.6,1.2) {$a=$};

\end{tikzpicture}

\end{center}
The strategy of proof is to divide $\operatorname{ZZ}(a,b)$ into smaller subsets by the $i$ value of the first edge. To that end, we define 
$$
I(l)=\{(p_1\to \dots \to p_k)\in \operatorname{ZZ}(x,y) \mid x,y\in C_n, \ n\in \mathbb Z_{\geq 0}, \ i(p_1\to p_2)=l\}.
$$
and our goal is to show the subsets $\operatorname{ZZ}(a,b)\cap I(l)$ are $R$-equivalent with subsets that define $A_n$ resulting in $\operatorname{ZZ}(a,b)\overset{R}{\sim}\operatorname{ZZ}(a,b)\cap A_n$.

It is immediate that 
\begin{equation}
\operatorname{ZZ}(a,b) \cap \big(I(3n-5)\cup I(3n-3)\cup I(3n-1)\big)=\emptyset  \label{Equation: Claim 3 emptyset 1}  
\end{equation}
as  characters at those indices are either  $1$:s or correspond to a merge of two blue circles. Next, let $(p_1\to \dots \to  p_{2k} )\in \operatorname{ZZ}(a,b)\cap I(3n)$, that is, $i(p_1\to p_2)=3n$. It follows $p_2$ is matched downwards to $p_3$ with $i(p_2\to p_3)=3n$ and so $(p_1\to p_3)=\fatslash_{n-2}(\romannumeral 1 \to \romannumeral 3)$. Hence 
\begin{align*}
    \operatorname{ZZ}(a,b)\cap I(3n)&= \{(p_1\to \dots \to p_{2k}) \in \operatorname{ZZ}(a,b) \mid (p_1\to p_3)=\fatslash_{n-2}(\romannumeral 1 \to \romannumeral 3) \} \\
    &  \overset{R}{\sim} \{(p_1\to \dots \to p_{2k}) \in \operatorname{ZZ}(a,b) \mid (p_1\to p_5)=\fatslash_{n-2}(\romannumeral 1 \to \romannumeral 5) \}
\end{align*}
where Lemma \ref{Lemma: jump move R-equivalence} is used for the $R$-equivalence. Supposing now that $p\in \operatorname{ZZ}(a,b)$  and  $(p_1\to p_5)=\fatslash_{n-2}(\romannumeral 1 \to \romannumeral 5)$ we can use Lemma \ref{Lemma: u/i} to obtain that $i(p_5 \to p_6)\geq 3n-5$. Since the unique arrow $f$ from $p_5$ with $i(f)=3n-5$ has the wrong direction, we can further deduce that $i(p_5 \to p_6)\geq 3n-6$. Thus $B_{n-2}(p)$ is well-defined and by Lemma \ref{Lemma: path after beta} iii) we get that:
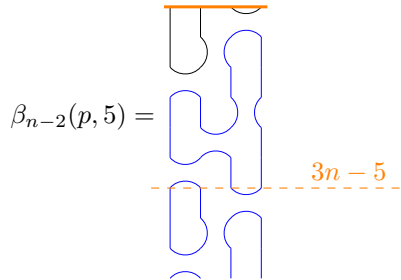
\begin{center}

\begin{tikzpicture}
\node at (0, -1.2) [anchor=south west, scale=0.4] {

    \begin{tikzpicture}
    
    \draw (0, 8) to (0, 9);
    \draw (0, 7) to (0, 8);
    \draw (0, 7) to[bend right=50] (1, 7);
    \draw (1, 7) to[bend right=50] (1, 8);
    \draw (1, 8) to (1, 9);
    \draw (2, 9) to[bend right=50] (3, 9);

    \begin{scope}[color=blue]

        \draw (0, 0) to[bend left=50] (1, 0);
    \draw (2, 0) to (2, 1);
    \draw (2, 1) to[bend left=50] (2, 2);
    \draw (2, 2) to[bend left=50] (3, 2);
    \draw (3, 1) to (3, 2);
    \draw (3, 0) to (3, 1);
    \draw (0, 3) to[bend left=50] (1, 3);
    \draw (1, 2) to (1, 3);
    \draw (1, 1) to[bend right=50] (1, 2);
    \draw (0, 1) to[bend right=50] (1, 1);
    \draw (0, 1) to (0, 2);
    \draw (0, 2) to (0, 3);
    \draw (2, 8) to[bend left=50] (3, 8);
    \draw (3, 7) to (3, 8);
    \draw (3, 6) to (3, 7);
    \draw (3, 5) to[bend left=50] (3, 6);
    \draw (3, 4) to (3, 5);
    \draw (3, 3) to (3, 4);
    \draw (2, 3) to[bend right=50] (3, 3);
    \draw (2, 3) to (2, 4);
    \draw (1, 4) to[bend left=50] (2, 4);
    \draw (0, 4) to[bend right=50] (1, 4);
    \draw (0, 4) to (0, 5);
    \draw (0, 5) to (0, 6);
    \draw (0, 6) to[bend left=50] (1, 6);
    \draw (1, 5) to (1, 6);
    \draw (1, 5) to[bend right=50] (2, 5);
    \draw (2, 5) to[bend right=50] (2, 6);
    \draw (2, 6) to (2, 7);
    \draw (2, 7) to[bend left=50] (2, 8);
\end{scope}
\begin{scope}[color=red]
\end{scope}

\draw[orange,line width=1mm] (-0.2,9) -- (3.2,9);

    \end{tikzpicture}

        };

\node at (-1,1) {$ \beta_{n-2}(p,5)=$};

\draw[orange,dashed] (-0.1,0.05) -- (3.2,0.05);
\node[orange] at (2.5,0.25) {$3n-5$};

\end{tikzpicture}

\end{center}
where $\beta_{n-2}(p,5)$ is not matched with any arrow $f$ having $u(f)\leq 3n-6$. It follows quickly that $\beta_{n-2}(p,5)$ is an unmatched cell and thus we have shown
\begin{equation}\label{Equation: Claim 3 relation 1}
    \operatorname{ZZ}(a,b)\cap I(3n)  \overset{R}{\sim} \{ p\in \operatorname{ZZ}(a,b)  \mid p=\fatslash_{n-2} B_{n-2}(p) \}.    
\end{equation}

Similarly, repeatedly using Lemmas \ref{Lemma: Lego moves} and \ref{Lemma: jump move R-equivalence} one proves that
\begin{align}
        \operatorname{ZZ}(a,b)\cap \bigcup_{l=1}^{3n-6} I(l) & \overset{R}{\sim} \{ p\in \operatorname{ZZ}(a,b)  \mid p=B_{n-2}(p) \} \label{Equation: Claim 3 relation 2}\\
    \operatorname{ZZ}(a,b)\cap I(3n-2) & \overset{R}{\sim} \{ p\in \operatorname{ZZ}(a,b)  \mid p=\lozenge_{n-2} B_{n-2}(p)\lozenge_{n-2}^{-1} \} \label{Equation: Claim 3 relation 3}\\
    \operatorname{ZZ}(a,b)\cap I(3n-4) & \overset{R}{\sim} \{ p\in \operatorname{ZZ}(a,b)  \mid p=\triangledown_{n-2} B_{n-2}(p)\overline{ \triangledown}^{-1} \}. \label{Equation: Claim 3 relation 4}
\end{align}
Now combining Relations \ref{Equation: Claim 3 emptyset 1}, \ref{Equation: Claim 3 relation 1}, \ref{Equation: Claim 3 relation 2}, \ref{Equation: Claim 3 relation 3} and \ref{Equation: Claim 3 relation 4} yields
\begin{align*}
    \operatorname{ZZ}(a,b) = & \left(\operatorname{ZZ}(a,b)\cap \bigcup_{l=1}^{3n-6} I(l) \right) \cup \big(  \operatorname{ZZ}(a,b)\cap  I(3n-4) \big) \notag \\ 
    & \quad \cup \big(  \operatorname{ZZ}(a,b)\cap  I(3n-2) \big)\cup \big(  \operatorname{ZZ}(a,b)\cap  I(3n) \big) \\
    \overset{R}{\sim}& \{ p\in \operatorname{ZZ}(a,b) \mid p\in \{ B_{n-2}(p),\  \fatslash_{n-2} B_{n-2}(p), \  \lozenge_{n-2} B_{n-2}(p)\lozenge_{n-2}^{-1}, \ \triangledown_{n-2} B_{n-2}(p)\overline{ \triangledown}^{-1} \}\} \\
    =&\operatorname{ZZ}(a,b) \cap A_n.
\end{align*}

For showing that $\operatorname{ZZ}(\mathfrak c_n(a),\mathfrak c_n(b)) \overset{R}{\sim} \operatorname{ZZ}(\mathfrak c_n(a),\mathfrak c_n(b)) \cap \operatorname{im}\Phi_n$ one analogously verifies  
\begin{align}
        \operatorname{ZZ}(\mathfrak c_n(a),\mathfrak c_n(b))\cap \bigcup_{l=1}^{n-6} I(l) & \overset{R}{\sim} \{ p\in \operatorname{ZZ}(\mathfrak c_n(a),\mathfrak c_n(b))  \mid p=B_{n-2}(p) \}  \label{Equation: Claim 3 relation 5} \\
    \operatorname{ZZ}(\mathfrak c_n(a),\mathfrak c_n(b))\cap I(3n+12) & \overset{R}{\sim} \{ p\in \operatorname{ZZ}(\mathfrak c_n(a),\mathfrak c_n(b))  \mid p=\fatslash_{n+2} \lozenge_n  \fatslash_{n-2} B_{n-2}(p)  \lozenge_n^{-1} \}   \label{Equation: Claim 3 relation 6} \\
    \operatorname{ZZ}(\mathfrak c_n(a),\mathfrak c_n(b))\cap I(3n+10) & \overset{R}{\sim} \{ p\in \operatorname{ZZ}(\mathfrak c_n(a),\mathfrak c_n(b))  \mid p=\lozenge_{n+2}\fatslash_{n} \lozenge_{n-2} B_{n-2}(p)  \lozenge_{n-2}^{-1}\lozenge_{n+2}^{-1} \}    \label{Equation: Claim 3 relation 7}  \\
    \operatorname{ZZ}(\mathfrak c_n(a),\mathfrak c_n(b))\cap I(3n+8) & \overset{R}{\sim} \{ p\in \operatorname{ZZ}(\mathfrak c_n(a),\mathfrak c_n(b))  \mid p=\triangledown_{n+2}\triangledown_{n}\triangledown_{n-2} B_{n-2}(p) \triangledown_{n-1}^{-1}\triangledown_{n+1}^{-1}\overline{\triangledown}^{-1} \}.    \label{Equation: Claim 3 relation 8}
\end{align}
The subpaths $B_{n-2}$ in the right-hand side of Relations \ref{Equation: Claim 3 relation 5}-\ref{Equation: Claim 3 relation 8}  belong in their respective graphs $G_{r,n+4}$. These graphs $G_{r,n+4}$ are isomorphic to graphs $G_{s,n}$ of Relations \ref{Equation: Claim 3 relation 1}-\ref{Equation: Claim 3 relation 4} via graph isomorphisms $\Lambda$. It follows that the paths in the right-hand side of   Relations \ref{Equation: Claim 3 relation 5}-\ref{Equation: Claim 3 relation 8} are precisely those of $\operatorname{im} \Phi_n$.

One can also see that for all $l\in 3n-5,\dots, 3n+6, 3n+7,3n+9,3n+11$ 
$$
\operatorname{ZZ}(\mathfrak c_n(a),\mathfrak c_n(b))\cap I(l)=\emptyset.
$$
by using Lemma \ref{Lemma: middle exclusion i}. Combining all of this together yields:
$$
\operatorname{ZZ}(\mathfrak c_n(a),\mathfrak c_n(b)) \overset{R}{\sim} \operatorname{ZZ}(\mathfrak c_n(a),\mathfrak c_n(b))\cap \operatorname{im} \Phi_n.
$$

The secondary case $a=v.\bzero 0101\bzero\bzero01$ is proven similarly: Relations 
\begin{align}
    \operatorname{ZZ}(a,b)\cap \bigcup_{l=1}^{3n-6} I(l) &\overset{R}{\sim} \{ p\in \operatorname{ZZ}(a,b)  \mid p=B_{n-2}(p) \} \\
    \operatorname{ZZ}(a,b)\cap I(3n-1) & \overset{R}{\sim} \{ p\in \operatorname{ZZ}(a,b)  \mid p\in  \{\overline{\triangledown} \triangledown_{n-3} B_{n-3}(p) \triangledown_{n-2}^{-1} ,\ \overline{\triangledown} \hspace{-1.5mm} \fatslash_{n-3} B_{n-3}(p)  \overline{\triangledown}^{-1} \}   \} \\
    \operatorname{ZZ}(a,b)\cap I(3n-2) & \overset{R}{\sim} \{ p\in \operatorname{ZZ}(a,b)  \mid p=\overline{\fatslash} \, \lozenge_{n-3} B_{n-3}(p) \lozenge_{n-3}^{-1}  \} \\
    \operatorname{ZZ}(a,b)\cap I(3n-3) & \overset{R}{\sim} \{ p\in \operatorname{ZZ}(a,b)  \mid p=\fatslash_{n-3} B_{n-3}(p) \} 
\end{align}
yield $\operatorname{ZZ}(a,b) \overset{R}{\sim} \operatorname{ZZ}(a,b)\cap A_n $  and $\operatorname{ZZ}(a,b) \cap I(3n-8)=\emptyset $ by Lemma \ref{Lemma: middle exclusion i}. For the codomain, relations 
\begin{align}
    \operatorname{ZZ}(\mathfrak c_n(a),\mathfrak c_n(b))\cap \bigcup_{l=1}^{n-6} I(l) & \overset{R}{\sim} \{ p\in \operatorname{ZZ}(\mathfrak c_n(a),\mathfrak c_n(b))  \mid p=B_{n-2}(p) \} \\
    \operatorname{ZZ}(\mathfrak c_n(a),\mathfrak c_n(b))\cap I(3n+11) & \overset{R}{\sim} \{ p\in \operatorname{ZZ}(\mathfrak c_n(a),\mathfrak c_n(b))  \mid p\in \notag \\ 
    \{\overline{\triangledown}\triangledown_{n+1}\triangledown_{n-1} \triangledown_{n-3} &B_{n-3}(p) \triangledown_{n-2}^{-1}\triangledown_{n}^{-1} \triangledown_{n+2}^{-1}, 
     \overline{\triangledown}\fatslash_{n+1}\lozenge_{n-1} \fatslash_{n-3} B_{n-3}(p) \lozenge_{n-1}^{-1} \overline{\triangledown}^{-1} \} \} \\
    \operatorname{ZZ}(\mathfrak c_n(a),\mathfrak c_n(b))\cap I(3n+10) & \overset{R}{\sim} \{ p\in \operatorname{ZZ}(\mathfrak c_n(a),\mathfrak c_n(b))  \mid p= \overline{\fatslash}\, \lozenge_{n+1}\fatslash_{n-1} \lozenge_{n-3} B_{n-3}(p) \lozenge_{n-3}^{-1}\lozenge_{n+1}^{-1} \} \\
    \operatorname{ZZ}(\mathfrak c_n(a),\mathfrak c_n(b))\cap I(3n+9) & \overset{R}{\sim} \{ p\in \operatorname{ZZ}(\mathfrak c_n(a),\mathfrak c_n(b))  \mid p=\fatslash_{n+1} \lozenge_{n-1}  \fatslash_{n-3} B_{n-3}(p)  \lozenge_{n-1}^{-1} \} 
\end{align}
yield 
$
\operatorname{ZZ}(\mathfrak c_n(a),\mathfrak c_n(b)) \overset{R}{\sim} \operatorname{ZZ}(\mathfrak c_n(a),\mathfrak c_n(b))\cap \operatorname{im} \Phi_n 
$ whereas $\operatorname{ZZ}(\mathfrak c_n(a),\mathfrak c_n(b))\cap \bigcup_{l=3n-5}^{3n+8} I(l)= \emptyset$. 

\end{proof}

\section{Matchings on small braids: a numerical investigation}\label{Section: numerical evidence}

As a test set for our numerical examination, we consider knots up to 12 crossings (12 included). For this set of knots we obtain their (highly non-unique) braid diagrams from KnotInfo \cite{knotinfo} and call this set of 2977 diagrams $\mathcal{B}$. By running an exhaustive search we obtain the following.
\begin{numericalresult}
   For 2976 out of 2977 braid diagrams $B\in \mathcal{B}$, the graphs $G(\Psi \llbracket B \rrbracket, \Mgr) $ are acyclic. For exactly one braid diagram $\operatorname{BD}(12_{n784})\in \mathcal{B}$ the graph  $G(\Psi \llbracket \operatorname{BD}(12_{n784}) \rrbracket, \Mgr)$, contains directed cycles. Those cycles are of length 12 or longer and the graph $G(\Psi \llbracket \operatorname{BD}(12_{n784}) \rrbracket, \Mgr) $ contains no cycles of length strictly less than 12.
\end{numericalresult}
The exceptional braid diagram $\operatorname{BD}(12_{n784})\in \mathcal{B}$ has 5 strands and 14 crossings and is defined as 
$$
\operatorname{BD}(12_{n784})=\sigma_1^{-1} \sigma_4^{-1} \sigma_2^{-1}\sigma_3\sigma_2 \sigma_4^{-1} \sigma_3\sigma_3 \sigma_2\sigma_2 \sigma_1^{-1} \sigma_4^{-1} \sigma_3 \sigma_2.
$$ 
By removing crossings from $\Gamma$ through trial and error, we find a smaller braid diagram $\Gamma$ with 10 crossings such that $G(\Psi \llbracket \Gamma \rrbracket,\Mgr )$ also contains an analogous directed cycle and
$$
\Gamma = \sigma_1^{-1} \sigma_4^{-1} \sigma_3 \sigma_2 \sigma_3 \sigma_3\sigma_2 \sigma_2 \sigma_1^{-1} \sigma_4^{-1}.
$$
Additionally the cycle is unique in $G(\Psi \llbracket \Gamma \rrbracket,\Mgr )$ and due to its more reasonable size, we can depict it along side $\Gamma$ in Figure \ref{Figure: exceptional braid and its cycle}:
\begin{figure}[ht]
    \centering
    \input{cycle_as_separate_files/cycle_total}
    \caption{The braid diagram $\Gamma$ and the unique cycle $x_1\to \dots \to x_{12} \to x_1$ of length 12 in the graph $G(\Psi \llbracket \Gamma_2 \rrbracket,\Mgr )$ }
    \label{Figure: exceptional braid and its cycle}
\end{figure}

To numerically measure the effectiveness of discrete Morse theory, we compare the sizes of delooped complexes $\Psi \llbracket B \rrbracket$ to Morse complexes $\Mgr \Psi \llbracket B \rrbracket$, $\Mlex \Psi \llbracket B \rrbracket$  and to the minimal complexes. We limit our test set to  $\mathcal{B}\setminus \{\operatorname{BD}(12_{n784}) \}$ since the greedy matching does not generate a sensible complex on the exceptional braid diagram. 

Given a braid $B$, we write
$
\min \llbracket B \rrbracket= \min \{ \dimkom(C)  \mid C\simeq \llbracket B \rrbracket \} \in \mathbb Z_{\geq 0}.
$
There are a number of programs which compute integral Khovanov homology of link diagrams using Bar-Natan's algorithm and therefore internally compute $\min\llbracket B \rrbracket$. However (at least without minor tweaking), most of them do not offer inputs and outputs for braid/tangle diagrams as opposed to complete links. As a subsidiary for $\min \llbracket B \rrbracket$, we compute $\min (\llbracket B \rrbracket; {\mathbb F_2})$ which is the minimal size of the complex when taken with $\mathbb F_2$ coefficients and which can be readily computed with \texttt{kht++} \cite{khtpp}. 

It is easy to see that $\min (\llbracket B \rrbracket; {\mathbb F_2})\leq \min \llbracket B \rrbracket$ and since odd torsion seems to be quite rare in Khovanov homology of complete links it seems reasonable to assume that $\min (\llbracket B \rrbracket; {\mathbb F_2})$ is a good estimate for $\min \llbracket B \rrbracket$. A summary of these results can be seen in Table \ref{Figure: numerical complex sizes} and a more detailed list of the results and the software generated them can be found in \cite{computerCodeAndDataFor4BraidPaper}.

\begingroup

\begin{table}[t]
\centering
\setlength{\tabcolsep}{10pt} 
\renewcommand{\arraystretch}{1.5} 

\begin{tabular}{|c||c|c|c|c|c}
    \hline
     $\mathcal{B}\setminus \{\operatorname{BD}(12_{n784}) \}$& $\min( \llbracket B \rrbracket ; \mathbb F_2 )$ & $\Mgr \Psi \llbracket B \rrbracket$ & $\Mlex \Psi \llbracket B \rrbracket$ & $\Psi \llbracket B \rrbracket$ \\
     \hline
     \hline
     total cell count & $2720982$ & $3217752$ & $7632284$& $226444836$ \\
     \hline
      \rule{0pt}{5ex}
      \parbox[c]{3.5cm}{ratio of total cell count  \\ to total cell count of \\  minimal complexes}
    & 1.00 & $ 1.18$ &$2.80$&$83.22$ \\
     \hline
     minimality obtained & 2976/2976 &871/2976 &396/2976& 0/2976\\
     \hline 
\end{tabular}
\caption{
A summary of sizes of various complexes for the test set   $\mathcal{B}\setminus \{\operatorname{BD}(12_{n784}) \}$ of 2976 braid diagrams. On the first row containing numbers, the sum of minimizers $\sum_B \min( \llbracket B \rrbracket ; \mathbb F_2 ) $ as well as the sums $\sum_B \dimkom (-)$ for different complexes of $B$ are displayed. The second numerical row displays ratios in the row above and the last row displays the number of cases where the Morse complexes end up being minimal. It is noteworthy that the complexes $\Mlex \Psi \llbracket B \rrbracket$ obtain minimality exactly when $B$ is an alternating braid diagram. The greedy complexes $\Mgr \Psi \llbracket B \rrbracket$ obtain minimality in 475 non-alternating braid diagrams in addition to the 396 alternating ones. 
}
\label{Figure: numerical complex sizes}
\end{table}

\endgroup

\FloatBarrier
\newpage

\section*{Appendix A: \texorpdfstring{$\operatorname{Kh}^{i,j} (T(4,-n))$}{Hij (T(4,-n))} in lowest and highest degrees}
\addcontentsline{toc}{section}{Appendix A: \texorpdfstring{$\operatorname{Kh}^{i,j} (T(4,-n))$}{Hij (T(4,-n))} in lowest and highest degrees}

This appendix is devoted to five tables presenting the even, unreduced integral Khovanov homology of negative 4-strand torus links in the lowest and highest non-trivial homological degrees. In combination, Figures \ref{Figure: T4 middle}, \ref{Figure: T4 homology lowest degrees} and \ref{Figure: T4 homology highest degrees} describe non-vanishing homology groups $\operatorname{Kh}^{i,j}(T(4,-n))$ for every $i\in \mathbb Z$ and $n\geq 28$  whenever such nontrivial groups exist. 
To save space in the large tables, the symbols $\oplus$ are omitted, e.g., $\mathbb Z \ \mathbb Z^4_2 \ \mathbb Z_4$ denotes $\mathbb Z \oplus (\mathbb Z/ 2\mathbb Z)^4 \oplus \mathbb Z/ 4\mathbb Z$.

\begin{figure}[ht]
    \begin{subfigure}{\textwidth}
        \centering
          \resizebox{\textwidth}{!}{

\footnotesize\begin{tabular}{|c||c|c|c|c|c|c|c|c|c|c|c|}
\hline
  &$-8n$ &&&&&&&&&& $-8n+10$\\
\hline
\hline
  $-24n+24$&&  &  &  &  &  &  &  &  &  & $\mathbb{{Z}}^{  }_{  }$\\
\hline
  &&  &  &  &  &  &  &  &  &  & $\mathbb{{Z}}^{ 4 }_{  }  \ \mathbb{{Z}}^{  }_{ 2 }$\\
\hline
  &&  &  &  &  &  &  &  & $\mathbb{{Z}}^{  }_{  }$ & $\mathbb{{Z}}^{  }_{  }$ & $\mathbb{{Z}}^{  }_{  }  \ \mathbb{{Z}}^{ 4 }_{ 2 }  \ \mathbb{{Z}}^{  }_{ 4 }$\\
\hline
  &&  &  &  &  &  &  &  & $\mathbb{{Z}}^{ 4 }_{  }  \ \mathbb{{Z}}^{  }_{ 2 }$ & $\mathbb{{Z}}^{ 5 }_{  }  \ \mathbb{{Z}}^{  }_{ 2 }$ & $\mathbb{{Z}}^{  }_{  }  \ \mathbb{{Z}}^{ 2 }_{ 2 }  \ \mathbb{{Z}}^{  }_{ 4 }$\\
\hline
  &&  &  &  &  &  & $\mathbb{{Z}}^{  }_{  }$ & $\mathbb{{Z}}^{  }_{  }$ & $\mathbb{{Z}}^{  }_{  }  \ \mathbb{{Z}}^{ 4 }_{ 2 }  \ \mathbb{{Z}}^{  }_{ 4 }$ & $\mathbb{{Z}}^{ 2 }_{  }  \ \mathbb{{Z}}^{ 2 }_{ 2 }  \ \mathbb{{Z}}^{  }_{ 4 }$ & $\mathbb{{Z}}^{  }_{  }  \ \mathbb{{Z}}^{  }_{ 2 }$\\
\hline
  &&  &  &  &  &  & $\mathbb{{Z}}^{ 4 }_{  }  \ \mathbb{{Z}}^{  }_{ 2 }$ & $\mathbb{{Z}}^{ 5 }_{  }  \ \mathbb{{Z}}^{  }_{ 2 }$ & $\mathbb{{Z}}^{  }_{  }  \ \mathbb{{Z}}^{  }_{ 2 }  \ \mathbb{{Z}}^{  }_{ 4 }$ & $\mathbb{{Z}}^{  }_{ 2 }  \ \mathbb{{Z}}^{  }_{ 4 }$ & \\
\hline
  &&  &  &  & $\mathbb{{Z}}^{  }_{  }$ & $\mathbb{{Z}}^{  }_{  }$ & $\mathbb{{Z}}^{ 4 }_{ 2 }  \ \mathbb{{Z}}^{  }_{ 4 }$ & $\mathbb{{Z}}^{  }_{  }  \ \mathbb{{Z}}^{ 2 }_{ 2 }  \ \mathbb{{Z}}^{  }_{ 4 }$ & $\mathbb{{Z}}^{  }_{  }$ &  & \\
\hline
  &&  &  &  & $\mathbb{{Z}}^{ 4 }_{  }  \ \mathbb{{Z}}^{  }_{ 2 }$ & $\mathbb{{Z}}^{ 5 }_{  }  \ \mathbb{{Z}}^{  }_{ 2 }$ & $\mathbb{{Z}}^{  }_{  }  \ \mathbb{{Z}}^{  }_{ 2 }$ & $\mathbb{{Z}}^{  }_{ 2 }$ &  &  & \\
\hline
  &&  & $\mathbb{{Z}}^{  }_{  }$ & $\mathbb{{Z}}^{  }_{  }$ & $\mathbb{{Z}}^{ 3 }_{ 2 }  \ \mathbb{{Z}}^{  }_{ 4 }$ & $\mathbb{{Z}}^{  }_{ 2 }  \ \mathbb{{Z}}^{  }_{ 4 }$ &  &  &  &  & \\
\hline
  &&  & $\mathbb{{Z}}^{ 3 }_{  }  \ \mathbb{{Z}}^{  }_{ 2 }$ & $\mathbb{{Z}}^{ 4 }_{  }  \ \mathbb{{Z}}^{  }_{ 2 }$ & $\mathbb{{Z}}^{  }_{  }$ &  &  &  &  &  & \\
\hline
 &$\mathbb{{Z}}^{  }_{  }$ & $\mathbb{{Z}}^{  }_{  }$ & $\mathbb{{Z}}^{ 3 }_{ 2 }$ & $\mathbb{{Z}}^{  }_{ 2 }$ &  &  &  &  &  &  & \\
\hline
 $-24n+2$&$\mathbb{{Z}}^{ 3 }_{  }$ & $\mathbb{{Z}}^{ 3 }_{  }$ &  &  &  &  &  &  &  &  & \\
\hline
 $-24n$&$\mathbb{{Z}}^{ 2 }_{  }$ &  &  &  &  &  &  &  &  &  & \\
\hline
\end{tabular}

    }
        \caption{For $n\geq 7$ the unreduced Khovanov homology $\operatorname{Kh}^{i,j}(T(4,-4n))$ for homological degrees $i\leq-8n+10$.}
    \end{subfigure}
        \begin{subfigure}{\textwidth}

     \centering   
    \resizebox{\textwidth}{!}{

\footnotesize\begin{tabular}{|c||c|c|c|c|c|c|c|c|c|c|c|}
\hline
 & $-8n-4$ &&&&&&&&&& $-8n+6$\\
\hline
\hline
$-24n+10$ &  &  &  &  &  &  &  &  &  &  & $\mathbb{{Z}}^{ 2 }_{  }$\\
\hline
 &  &  &  &  &  &  &  &  &  & $\mathbb{{Z}}^{  }_{  }$ & $\mathbb{{Z}}^{ 2 }_{  } \ \mathbb{{Z}}^{ 2 }_{ 2 } \ \mathbb{{Z}}^{  }_{ 4 }$\\
\hline
 &  &  &  &  &  &  &  &  & $\mathbb{{Z}}^{ 2 }_{  }$ & $\mathbb{{Z}}^{ 3 }_{  } \ \mathbb{{Z}}^{ 2 }_{ 2 }$ & $\mathbb{{Z}}^{  }_{  } \ \mathbb{{Z}}^{ 2 }_{ 2 } \ \mathbb{{Z}}^{  }_{ 4 }$\\
\hline
 &  &  &  &  &  &  &  & $\mathbb{{Z}}^{  }_{  }$ & $\mathbb{{Z}}^{ 2 }_{  } \ \mathbb{{Z}}^{ 2 }_{ 2 } \ \mathbb{{Z}}^{  }_{ 4 }$ & $\mathbb{{Z}}^{ 2 }_{  } \ \mathbb{{Z}}^{ 2 }_{ 2 } \ \mathbb{{Z}}^{  }_{ 4 }$ & $\mathbb{{Z}}^{  }_{  } \ \mathbb{{Z}}^{  }_{ 2 }$\\
\hline
 &  &  &  &  &  &  & $\mathbb{{Z}}^{ 2 }_{  }$ & $\mathbb{{Z}}^{ 3 }_{  } \ \mathbb{{Z}}^{ 2 }_{ 2 }$ & $\mathbb{{Z}}^{  }_{  } \ \mathbb{{Z}}^{  }_{ 2 } \ \mathbb{{Z}}^{  }_{ 4 }$ & $\mathbb{{Z}}^{  }_{ 2 } \ \mathbb{{Z}}^{  }_{ 4 }$ & \\
\hline
&  &  &  &  &  & $\mathbb{{Z}}^{  }_{  }$ & $\mathbb{{Z}}^{  }_{  } \ \mathbb{{Z}}^{ 2 }_{ 2 } \ \mathbb{{Z}}^{  }_{ 4 }$ & $\mathbb{{Z}}^{  }_{  } \ \mathbb{{Z}}^{ 2 }_{ 2 } \ \mathbb{{Z}}^{  }_{ 4 }$ & $\mathbb{{Z}}^{  }_{  }$ &  & \\
\hline
 &  &  &  &  & $\mathbb{{Z}}^{ 2 }_{  }$ & $\mathbb{{Z}}^{ 3 }_{  } \ \mathbb{{Z}}^{ 2 }_{ 2 }$ & $\mathbb{{Z}}^{  }_{  } \ \mathbb{{Z}}^{  }_{ 2 }$ & $\mathbb{{Z}}^{  }_{ 2 }$ &  &  & \\
\hline
&  &  &  & $\mathbb{{Z}}^{  }_{  }$ & $\mathbb{{Z}}^{  }_{  } \ \mathbb{{Z}}^{  }_{ 2 } \ \mathbb{{Z}}^{  }_{ 4 }$ & $\mathbb{{Z}}^{  }_{ 2 } \ \mathbb{{Z}}^{  }_{ 4 }$ &  &  &  &  & \\
\hline
 &  &  & $\mathbb{{Z}}^{  }_{  }$ & $\mathbb{{Z}}^{ 2 }_{  } \ \mathbb{{Z}}^{ 2 }_{ 2 }$ & $\mathbb{{Z}}^{  }_{  }$ &  &  &  &  &  & \\
\hline
 &  & $\mathbb{{Z}}^{  }_{  }$ & $\mathbb{{Z}}^{  }_{  } \ \mathbb{{Z}}^{  }_{ 2 }$ & $\mathbb{{Z}}^{  }_{ 2 }$ &  &  &  &  &  &  & \\
\hline
$-24n-10$ & $\mathbb{{Z}}^{  }_{  }$ & $\mathbb{{Z}}^{  }_{  } \ \mathbb{{Z}}^{  }_{ 2 }$ &  &  &  &  &  &  &  &  & \\
\hline
$-24n-12$ & $\mathbb{{Z}}^{ 2 }_{  }$ &  &  &  &  &  &  &  &  &  & \\
\hline
\end{tabular}


    }
        \caption{For $n\geq 7$ the unreduced Khovanov homology $\operatorname{Kh}^{i,j}(T(4,-4n-2))$ for homological degrees $i\leq-8n+6$.}
    \end{subfigure}

    \begin{subfigure}{\textwidth}
        \centering

\footnotesize \begin{tabular}{|c||c|c|c|c|c|c|c|c|}
\hline
 & $-4n-1$ &$-4n$ &&&&&& $-4n+6$\\
\hline
\hline
$-12n+13$ &  &  &  &  &  &  &  & $\mathbb{{Z}}^{  }_{  }$\\
\hline
 &  &  &  &  &  &  &  & $\mathbb{{Z}}^{  }_{  }  \ \mathbb{{Z}}^{  }_{ 2 }  \ \mathbb{{Z}}^{  }_{ 4 }$\\
\hline
 &  &  &  &  &  & $\mathbb{{Z}}^{  }_{  }$ & $\mathbb{{Z}}^{ 2 }_{  }  \ \mathbb{{Z}}^{  }_{ 2 }$ & $\mathbb{{Z}}^{  }_{  }  \ \mathbb{{Z}}^{ 2 }_{ 2 }  \ \mathbb{{Z}}^{  }_{ 4 }$\\
\hline
 &  &  &  &  &  & $\mathbb{{Z}}^{  }_{  }  \ \mathbb{{Z}}^{  }_{ 2 }  \ \mathbb{{Z}}^{  }_{ 4 }$ & $\mathbb{{Z}}^{ 2 }_{  }  \ \mathbb{{Z}}^{ 2 }_{ 2 }  \ \mathbb{{Z}}^{  }_{ 4 }$ & $\mathbb{{Z}}^{  }_{  }  \ \mathbb{{Z}}^{  }_{ 2 }$\\
\hline
 &  &  &  & $\mathbb{{Z}}^{  }_{  }$ & $\mathbb{{Z}}^{ 2 }_{  }  \ \mathbb{{Z}}^{  }_{ 2 }$ & $\mathbb{{Z}}^{  }_{  }  \ \mathbb{{Z}}^{  }_{ 2 }  \ \mathbb{{Z}}^{  }_{ 4 }$ & $\mathbb{{Z}}^{  }_{ 2 }  \ \mathbb{{Z}}^{  }_{ 4 }$ & \\
\hline
 &  &  &  & $\mathbb{{Z}}^{  }_{ 2 }  \ \mathbb{{Z}}^{  }_{ 4 }$ & $\mathbb{{Z}}^{  }_{  }  \ \mathbb{{Z}}^{ 2 }_{ 2 }  \ \mathbb{{Z}}^{  }_{ 4 }$ & $\mathbb{{Z}}^{  }_{  }$ &  & \\
\hline
 &  & $\mathbb{{Z}}^{  }_{  }$ & $\mathbb{{Z}}^{ 2 }_{  }  \ \mathbb{{Z}}^{  }_{ 2 }$ & $\mathbb{{Z}}^{  }_{  }  \ \mathbb{{Z}}^{  }_{ 2 }$ & $\mathbb{{Z}}^{  }_{ 2 }$ &  &  & \\
\hline
 &  & $\mathbb{{Z}}^{  }_{ 4 }$ & $\mathbb{{Z}}^{  }_{ 2 }  \ \mathbb{{Z}}^{  }_{ 4 }$ &  &  &  &  & \\
\hline
$-12n-3$ & $\mathbb{{Z}}^{  }_{  }  \ \mathbb{{Z}}^{  }_{ 2 }$ & $\mathbb{{Z}}^{  }_{  }$ &  &  &  &  &  & \\
\hline
$-12n-5$ & $\mathbb{{Z}}^{  }_{ 2 }$ &  &  &  &  &  &  & \\
\hline
\end{tabular}

        \caption{For $n\geq 14$ the unreduced Khovanov homology $\operatorname{Kh}^{i,j}(T(4,-2n-1))$ for homological degrees $i\leq-4n+6$.}
    \end{subfigure}
    
    \caption{For $n\geq 28$ the unreduced integral Khovanov homology $\operatorname{Kh}^{i,j}(T(4,-n))$ in the lowest non-trivial homological degrees. Outside the marked entries the homology vanishes for $i\leq -8n +10$,  $i\leq -8n +6$ and $i\leq -4n +6$ respectively. }
    \label{Figure: T4 homology lowest degrees}
\end{figure}

\newpage

\begin{sidewaysfigure}[p]  
    \centering
    \begin{subfigure}{\textwidth}
        \centering
          \resizebox{\textwidth}{!}{
    
\begin{tabular}{|c||c|c|c|c|c|c|c|c|c|c|c|c|c|c|c|c|c|c|c|c|c|}
\hline
 & $-41$ & $-40$ & $-39$ & $-38$ & $-37$ & $-36$ & $-35$ & $-34$ & $-33$ & $-32$ & $-31$ & $-30$ & $-29$ & $-28$ & $-27$ & $-26$ & $-25$ & $-24$ & $-23$ & $-22$ & $-21$\\
\hline
\hline
$-3n-26$ &  &  &  &  &  &  &  &  &  &  &  &  &  &  &  &  &  &  &  & $\mathbb{{Z}}^{  }_{  }$ & \\
\hline
$-3n-28$ &  &  &  &  &  &  &  &  &  &  &  &  &  &  &  &  &  & $\mathbb{{Z}}^{  }_{  }$ &  & $\mathbb{{Z}}^{ 2 }_{  } \ \mathbb{{Z}}^{  }_{ 2 }$ & $\mathbb{{Z}}^{ 3 }_{  } \ \mathbb{{Z}}^{  }_{ 2 }$\\
\hline
 &  &  &  &  &  &  &  &  &  &  &  &  &  &  &  &  &  & $\mathbb{{Z}}^{ 2 }_{  } \ \mathbb{{Z}}^{  }_{ 2 }$ & $\mathbb{{Z}}^{ 2 }_{  }$ & $\mathbb{{Z}}^{ 2 }_{ 2 } \ \mathbb{{Z}}^{  }_{ 4 }$ & $\mathbb{{Z}}^{  }_{  } \ \mathbb{{Z}}^{ 2 }_{ 2 } \ \mathbb{{Z}}^{  }_{ 4 }$\\
\hline
&  &  &  &  &  &  &  &  &  &  &  &  &  &  &  & $\mathbb{{Z}}^{ 2 }_{  }$ & $\mathbb{{Z}}^{  }_{  }$ & $\mathbb{{Z}}^{  }_{  } \ \mathbb{{Z}}^{ 2 }_{ 2 } \ \mathbb{{Z}}^{  }_{ 4 }$ & $\mathbb{{Z}}^{ 3 }_{  } \ \mathbb{{Z}}^{ 2 }_{ 2 }$ & $\mathbb{{Z}}^{  }_{  } \ \mathbb{{Z}}^{  }_{ 2 }$ & $\mathbb{{Z}}^{  }_{ 2 }$\\
\hline
 &  &  &  &  &  &  &  &  &  &  &  &  &  & $\mathbb{{Z}}^{  }_{  }$ &  & $\mathbb{{Z}}^{  }_{  } \ \mathbb{{Z}}^{ 2 }_{ 2 } \ \mathbb{{Z}}^{  }_{ 4 }$ & $\mathbb{{Z}}^{ 3 }_{  } \ \mathbb{{Z}}^{  }_{ 2 }$ & $\mathbb{{Z}}^{  }_{  } \ \mathbb{{Z}}^{  }_{ 2 } \ \mathbb{{Z}}^{  }_{ 4 }$ & $\mathbb{{Z}}^{  }_{ 2 } \ \mathbb{{Z}}^{  }_{ 4 }$ &  & \\
\hline
 &  &  &  &  &  &  &  &  &  &  &  & $\mathbb{{Z}}^{  }_{  }$ &  & $\mathbb{{Z}}^{ 2 }_{  } \ \mathbb{{Z}}^{  }_{ 2 }$ & $\mathbb{{Z}}^{ 3 }_{  } \ \mathbb{{Z}}^{  }_{ 2 }$ & $\mathbb{{Z}}^{ 2 }_{ 2 } \ \mathbb{{Z}}^{  }_{ 4 }$ & $\mathbb{{Z}}^{  }_{  } \ \mathbb{{Z}}^{ 2 }_{ 2 } \ \mathbb{{Z}}^{  }_{ 4 }$ & $\mathbb{{Z}}^{  }_{  }$ &  &  & \\
\hline
 &  &  &  &  &  &  &  &  &  &  &  & $\mathbb{{Z}}^{ 2 }_{  } \ \mathbb{{Z}}^{  }_{ 2 }$ & $\mathbb{{Z}}^{ 2 }_{  }$ & $\mathbb{{Z}}^{  }_{  } \ \mathbb{{Z}}^{ 2 }_{ 2 } \ \mathbb{{Z}}^{  }_{ 4 }$ & $\mathbb{{Z}}^{ 2 }_{  } \ \mathbb{{Z}}^{ 2 }_{ 2 } \ \mathbb{{Z}}^{  }_{ 4 }$ & $\mathbb{{Z}}^{  }_{  } \ \mathbb{{Z}}^{  }_{ 2 }$ & $\mathbb{{Z}}^{  }_{ 2 }$ &  &  &  & \\
\hline
 &  &  &  &  &  &  &  &  &  & $\mathbb{{Z}}^{ 2 }_{  }$ & $\mathbb{{Z}}^{  }_{  }$ & $\mathbb{{Z}}^{  }_{  } \ \mathbb{{Z}}^{ 2 }_{ 2 } \ \mathbb{{Z}}^{  }_{ 4 }$ & $\mathbb{{Z}}^{ 3 }_{  } \ \mathbb{{Z}}^{ 2 }_{ 2 }$ & $\mathbb{{Z}}^{  }_{  } \ \mathbb{{Z}}^{  }_{ 2 } \ \mathbb{{Z}}^{  }_{ 4 }$ & $\mathbb{{Z}}^{  }_{ 2 } \ \mathbb{{Z}}^{  }_{ 4 }$ &  &  &  &  &  & \\
\hline
&  &  &  &  &  &  &  & $\mathbb{{Z}}^{  }_{  }$ &  & $\mathbb{{Z}}^{  }_{  } \ \mathbb{{Z}}^{ 2 }_{ 2 } \ \mathbb{{Z}}^{  }_{ 4 }$ & $\mathbb{{Z}}^{ 3 }_{  } \ \mathbb{{Z}}^{  }_{ 2 }$ & $\mathbb{{Z}}^{  }_{  } \ \mathbb{{Z}}^{ 2 }_{ 2 } \ \mathbb{{Z}}^{  }_{ 4 }$ & $\mathbb{{Z}}^{  }_{  } \ \mathbb{{Z}}^{ 2 }_{ 2 } \ \mathbb{{Z}}^{  }_{ 4 }$ & $\mathbb{{Z}}^{  }_{  }$ &  &  &  &  &  &  & \\
\hline
 &  &  &  &  &  & $\mathbb{{Z}}^{  }_{  }$ &  & $\mathbb{{Z}}^{ 2 }_{  } \ \mathbb{{Z}}^{  }_{ 2 }$ & $\mathbb{{Z}}^{ 3 }_{  } \ \mathbb{{Z}}^{  }_{ 2 }$ & $\mathbb{{Z}}^{  }_{  } \ \mathbb{{Z}}^{ 2 }_{ 2 } \ \mathbb{{Z}}^{  }_{ 4 }$ & $\mathbb{{Z}}^{ 2 }_{  } \ \mathbb{{Z}}^{ 2 }_{ 2 } \ \mathbb{{Z}}^{  }_{ 4 }$ & $\mathbb{{Z}}^{  }_{  } \ \mathbb{{Z}}^{  }_{ 2 }$ & $\mathbb{{Z}}^{  }_{ 2 }$ &  &  &  &  &  &  &  & \\
\hline
&  &  &  &  &  & $\mathbb{{Z}}^{ 2 }_{  } \ \mathbb{{Z}}^{  }_{ 2 }$ & $\mathbb{{Z}}^{ 2 }_{  }$ & $\mathbb{{Z}}^{  }_{  } \ \mathbb{{Z}}^{ 2 }_{ 2 } \ \mathbb{{Z}}^{  }_{ 4 }$ & $\mathbb{{Z}}^{ 2 }_{  } \ \mathbb{{Z}}^{ 2 }_{ 2 } \ \mathbb{{Z}}^{  }_{ 4 }$ & $\mathbb{{Z}}^{  }_{  } \ \mathbb{{Z}}^{  }_{ 2 } \ \mathbb{{Z}}^{  }_{ 4 }$ & $\mathbb{{Z}}^{  }_{ 2 } \ \mathbb{{Z}}^{  }_{ 4 }$ &  &  &  &  &  &  &  &  &  & \\
\hline
 &  &  &  & $\mathbb{{Z}}^{ 2 }_{  }$ & $\mathbb{{Z}}^{  }_{  }$ & $\mathbb{{Z}}^{  }_{  } \ \mathbb{{Z}}^{ 2 }_{ 2 } \ \mathbb{{Z}}^{  }_{ 4 }$ & $\mathbb{{Z}}^{ 3 }_{  } \ \mathbb{{Z}}^{ 2 }_{ 2 }$ & $\mathbb{{Z}}^{  }_{  } \ \mathbb{{Z}}^{ 2 }_{ 2 } \ \mathbb{{Z}}^{  }_{ 4 }$ & $\mathbb{{Z}}^{  }_{  } \ \mathbb{{Z}}^{ 2 }_{ 2 } \ \mathbb{{Z}}^{  }_{ 4 }$ & $\mathbb{{Z}}^{  }_{  }$ &  &  &  &  &  &  &  &  &  &  & \\
\hline
 &  & $\mathbb{{Z}}^{  }_{  }$ &  & $\mathbb{{Z}}^{  }_{  } \ \mathbb{{Z}}^{ 2 }_{ 2 } \ \mathbb{{Z}}^{  }_{ 4 }$ & $\mathbb{{Z}}^{ 3 }_{  } \ \mathbb{{Z}}^{  }_{ 2 }$ & $\mathbb{{Z}}^{ 2 }_{  } \ \mathbb{{Z}}^{ 2 }_{ 2 } \ \mathbb{{Z}}^{  }_{ 4 }$ & $\mathbb{{Z}}^{ 2 }_{  } \ \mathbb{{Z}}^{ 2 }_{ 2 } \ \mathbb{{Z}}^{  }_{ 4 }$ & $\mathbb{{Z}}^{  }_{  } \ \mathbb{{Z}}^{  }_{ 2 }$ & $\mathbb{{Z}}^{  }_{ 2 }$ &  &  &  &  &  &  &  &  &  &  &  & \\
\hline
&  & $\mathbb{{Z}}^{ 2 }_{  } \ \mathbb{{Z}}^{  }_{ 2 }$ & $\mathbb{{Z}}^{ 3 }_{  } \ \mathbb{{Z}}^{  }_{ 2 }$ & $\mathbb{{Z}}^{  }_{  } \ \mathbb{{Z}}^{ 2 }_{ 2 } \ \mathbb{{Z}}^{  }_{ 4 }$ & $\mathbb{{Z}}^{ 2 }_{  } \ \mathbb{{Z}}^{ 2 }_{ 2 } \ \mathbb{{Z}}^{  }_{ 4 }$ & $\mathbb{{Z}}^{  }_{  } \ \mathbb{{Z}}^{  }_{ 2 } \ \mathbb{{Z}}^{  }_{ 4 }$ & $\mathbb{{Z}}^{  }_{ 2 } \ \mathbb{{Z}}^{  }_{ 4 }$ &  &  &  &  &  &  &  &  &  &  &  &  &  & \\
\hline
 & $\mathbb{{Z}}^{ 2 }_{  }$ & $\mathbb{{Z}}^{  }_{  } \ \mathbb{{Z}}^{ 2 }_{ 2 } \ \mathbb{{Z}}^{  }_{ 4 }$ & $\mathbb{{Z}}^{ 2 }_{  } \ \mathbb{{Z}}^{ 2 }_{ 2 } \ \mathbb{{Z}}^{  }_{ 4 }$ & $\mathbb{{Z}}^{  }_{  } \ \mathbb{{Z}}^{ 2 }_{ 2 } \ \mathbb{{Z}}^{  }_{ 4 }$ & $\mathbb{{Z}}^{  }_{  } \ \mathbb{{Z}}^{ 2 }_{ 2 } \ \mathbb{{Z}}^{  }_{ 4 }$ & $\mathbb{{Z}}^{  }_{  }$ &  &  &  &  &  &  &  &  &  &  &  &  &  &  & \\
\hline
 & $\mathbb{{Z}}^{ 3 }_{  } \ \mathbb{{Z}}^{ 2 }_{ 2 }$ & $\mathbb{{Z}}^{ 2 }_{  } \ \mathbb{{Z}}^{ 2 }_{ 2 } \ \mathbb{{Z}}^{  }_{ 4 }$ & $\mathbb{{Z}}^{ 2 }_{  } \ \mathbb{{Z}}^{ 2 }_{ 2 } \ \mathbb{{Z}}^{  }_{ 4 }$ & $\mathbb{{Z}}^{  }_{  } \ \mathbb{{Z}}^{  }_{ 2 }$ & $\mathbb{{Z}}^{  }_{ 2 }$ &  &  &  &  &  &  &  &  &  &  &  &  &  &  &  & \\
\hline
 & $\mathbb{{Z}}^{ 2 }_{  } \ \mathbb{{Z}}^{ 2 }_{ 2 } \ \mathbb{{Z}}^{  }_{ 4 }$ & $\mathbb{{Z}}^{  }_{  } \ \mathbb{{Z}}^{  }_{ 2 } \ \mathbb{{Z}}^{  }_{ 4 }$ & $\mathbb{{Z}}^{  }_{ 2 } \ \mathbb{{Z}}^{  }_{ 4 }$ &  &  &  &  &  &  &  &  &  &  &  &  &  &  &  &  &  & \\
\hline
$-3n-60$ & $\mathbb{{Z}}^{  }_{  } \ \mathbb{{Z}}^{ 2 }_{ 2 } \ \mathbb{{Z}}^{  }_{ 4 }$ & $\mathbb{{Z}}^{  }_{  }$ &  &  &  &  &  &  &  &  &  &  &  &  &  &  &  &  &  &  & \\
\hline
$-3n-62$ & $\mathbb{{Z}}^{  }_{ 2 }$ &  &  &  &  &  &  &  &  &  &  &  &  &  &  &  &  &  &  &  & \\
\hline
\end{tabular}

    }
        \caption{For $n\geq 28$ the unreduced Khovanov homology $\operatorname{Kh}^{i,j}(T(4,-n))$ for homological degrees $i\in [-41, -21]$.}
    \end{subfigure}

    \begin{subfigure}{\textwidth}
        \centering
          \resizebox{\textwidth}{!}{

\begin{tabular}{|c||c|c|c|c|c|c|c|c|c|c|c|c|c|c|c|c|c|c|c|c|c|}
\hline
 & $-20$ & $-19$ & $-18$ & $-17$ & $-16$ & $-15$ & $-14$ & $-13$ & $-12$ & $-11$ & $-10$ & $-9$ & $-8$ & $-7$ & $-6$ & $-5$ & $-4$ & $-3$ & $-2$ & $-1$ & $0$\\
\hline
\hline
$-3n+4$ &  &  &  &  &  &  &  &  &  &  &  &  &  &  &  &  &  &  &  &  & $\mathbb{{Z}}^{  }_{  }$\\
\hline
$-3n+2$ &  &  &  &  &  &  &  &  &  &  &  &  &  &  &  &  &  &  &  &  & $\mathbb{{Z}}^{  }_{  }$\\
\hline
 &  &  &  &  &  &  &  &  &  &  &  &  &  &  &  &  &  &  & $\mathbb{{Z}}^{  }_{  }$ &  & \\
\hline
&  &  &  &  &  &  &  &  &  &  &  &  &  &  &  &  & $\mathbb{{Z}}^{  }_{  }$ &  & $\mathbb{{Z}}^{  }_{ 2 }$ &  & \\
\hline
 &  &  &  &  &  &  &  &  &  &  &  &  &  &  & $\mathbb{{Z}}^{  }_{  }$ &  & $\mathbb{{Z}}^{  }_{  }$ & $\mathbb{{Z}}^{  }_{  }$ &  &  & \\
\hline
 &  &  &  &  &  &  &  &  &  &  &  &  &  &  & $\mathbb{{Z}}^{  }_{  }  \  \mathbb{{Z}}^{  }_{ 2 }$ & $\mathbb{{Z}}^{  }_{  }$ &  &  &  &  & \\
\hline
 &  &  &  &  &  &  &  &  &  &  &  &  & $\mathbb{{Z}}^{ 2 }_{  }$ & $\mathbb{{Z}}^{  }_{  }$ & $\mathbb{{Z}}^{  }_{ 2 }$ & $\mathbb{{Z}}^{  }_{  }$ &  &  &  &  & \\
\hline
 &  &  &  &  &  &  &  &  &  &  & $\mathbb{{Z}}^{  }_{  }$ &  & $\mathbb{{Z}}^{  }_{ 2 }  \  \mathbb{{Z}}^{  }_{ 4 }$ & $\mathbb{{Z}}^{  }_{  }$ &  &  &  &  &  &  & \\
\hline
 &  &  &  &  &  &  &  &  & $\mathbb{{Z}}^{  }_{  }$ &  & $\mathbb{{Z}}^{  }_{  }  \  \mathbb{{Z}}^{  }_{ 2 }$ & $\mathbb{{Z}}^{ 2 }_{  }  \  \mathbb{{Z}}^{  }_{ 2 }$ &  &  &  &  &  &  &  &  & \\
\hline
 &  &  &  &  &  &  &  &  & $\mathbb{{Z}}^{ 2 }_{  }  \  \mathbb{{Z}}^{  }_{ 2 }$ & $\mathbb{{Z}}^{ 2 }_{  }$ & $\mathbb{{Z}}^{  }_{ 2 }$ & $\mathbb{{Z}}^{  }_{ 2 }$ &  &  &  &  &  &  &  &  & \\
\hline
 &  &  &  &  &  &  & $\mathbb{{Z}}^{ 2 }_{  }$ & $\mathbb{{Z}}^{  }_{  }$ & $\mathbb{{Z}}^{  }_{ 2 }  \  \mathbb{{Z}}^{  }_{ 4 }$ & $\mathbb{{Z}}^{  }_{  }  \  \mathbb{{Z}}^{  }_{ 2 }$ &  &  &  &  &  &  &  &  &  &  & \\
\hline
 &  &  &  &  & $\mathbb{{Z}}^{  }_{  }$ &  & $\mathbb{{Z}}^{ 2 }_{ 2 }  \  \mathbb{{Z}}^{  }_{ 4 }$ & $\mathbb{{Z}}^{ 2 }_{  }  \  \mathbb{{Z}}^{  }_{ 2 }$ & $\mathbb{{Z}}^{  }_{  }$ &  &  &  &  &  &  &  &  &  &  &  & \\
\hline
&  &  & $\mathbb{{Z}}^{  }_{  }$ &  & $\mathbb{{Z}}^{ 2 }_{  }  \  \mathbb{{Z}}^{  }_{ 2 }$ & $\mathbb{{Z}}^{ 3 }_{  }  \  \mathbb{{Z}}^{  }_{ 2 }$ & $\mathbb{{Z}}^{  }_{ 2 }$ & $\mathbb{{Z}}^{  }_{ 2 }$ &  &  &  &  &  &  &  &  &  &  &  &  & \\
\hline
&  &  & $\mathbb{{Z}}^{ 2 }_{  }  \  \mathbb{{Z}}^{  }_{ 2 }$ & $\mathbb{{Z}}^{ 2 }_{  }$ & $\mathbb{{Z}}^{  }_{ 2 }  \  \mathbb{{Z}}^{  }_{ 4 }$ & $\mathbb{{Z}}^{  }_{ 2 }  \  \mathbb{{Z}}^{  }_{ 4 }$ &  &  &  &  &  &  &  &  &  &  &  &  &  &  & \\
\hline
 & $\mathbb{{Z}}^{ 2 }_{  }$ & $\mathbb{{Z}}^{  }_{  }$ & $\mathbb{{Z}}^{ 2 }_{ 2 }  \  \mathbb{{Z}}^{  }_{ 4 }$ & $\mathbb{{Z}}^{ 2 }_{  }  \  \mathbb{{Z}}^{ 2 }_{ 2 }$ & $\mathbb{{Z}}^{  }_{  }$ &  &  &  &  &  &  &  &  &  &  &  &  &  &  &  & \\
\hline
 & $\mathbb{{Z}}^{  }_{  }  \  \mathbb{{Z}}^{ 2 }_{ 2 }  \  \mathbb{{Z}}^{  }_{ 4 }$ & $\mathbb{{Z}}^{ 3 }_{  }  \  \mathbb{{Z}}^{  }_{ 2 }$ & $\mathbb{{Z}}^{  }_{  }  \  \mathbb{{Z}}^{  }_{ 2 }$ & $\mathbb{{Z}}^{  }_{ 2 }$ &  &  &  &  &  &  &  &  &  &  &  &  &  &  &  &  & \\
\hline
$-3n-28$ & $\mathbb{{Z}}^{  }_{ 2 }  \  \mathbb{{Z}}^{  }_{ 4 }$ & $\mathbb{{Z}}^{  }_{ 2 }  \  \mathbb{{Z}}^{  }_{ 4 }$ &  &  &  &  &  &  &  &  &  &  &  &  &  &  &  &  &  &  & \\
\hline
$-3n-30$ & $\mathbb{{Z}}^{  }_{  }$ &  &  &  &  &  &  &  &  &  &  &  &  &  &  &  &  &  &  &  & \\
\hline
\end{tabular}

    }
        \caption{For $n\geq 28$ the unreduced Khovanov homology $\operatorname{Kh}^{i,j}(T(4,-n))$ for homological degrees $i\in [-20,0]$}
    \end{subfigure}
    
    \caption{For $n\geq 28$ the unreduced integral Khovanov homology $\operatorname{Kh}^{i,j}(T(4,-n))$ in the highest non-trivial homological degrees. Outside the marked entries the homology vanishes for $i\geq -41$.}
    \label{Figure: T4 homology highest degrees}
\end{sidewaysfigure}

\FloatBarrier
\printbibliography

\end{document}